\documentclass[a4paper,10pt,twoside]{article}

\usepackage[T1]{fontenc}
\usepackage[utf8]{inputenc}
\usepackage[english]{babel}

\usepackage[top=3cm,bottom=3.2cm,left=3.2cm,right=3.2cm]{geometry}

\usepackage{amscd}
\usepackage{amsthm}
\usepackage{bbm}
\usepackage{bm}
\usepackage{booktabs}
\usepackage{braket}
\usepackage{changepage}
\usepackage[babel]{csquotes}
\usepackage{enumitem}
\usepackage{fancyhdr}
\usepackage{float}
\usepackage[flushmargin,hang,symbol]{footmisc}
\usepackage{graphicx}
\usepackage{indentfirst}
\usepackage[intlimits]{mathtools}
\usepackage{mathrsfs}
\usepackage{microtype}
\usepackage{newtxtext,newtxmath}
\usepackage{subfig}
\usepackage[english]{varioref}
\usepackage{xcolor}
\usepackage{xurl}

\usepackage{hyperref}
    \hypersetup{
	colorlinks = true,
	citecolor = black,
	filecolor = black,
	linkcolor = black,
	urlcolor = black
    }

\selectfont

\makeatletter
\renewcommand{\@maketitle}{
	\begin{flushleft}\Large\@title\par\vspace{3ex}
	\large\@author\par\vspace{0ex}
	\end{flushleft}}
\makeatother

\makeatletter
\renewcommand{\footnoterule}{
  \kern -3pt
  \vspace{1.6\baselineskip}
  \hrule width 0.4\columnwidth
  \vspace{0.8\baselineskip} 
  \kern -0.4pt
}
\makeatother

\newcommand{\mail}[1]{\href{mailto:#1}{\texttt{#1}}}

\newcommand{\customlabel}[1]{\textup{\textbf{#1}}}

\newcommand{\N}{\mathbb{N}}

\newcommand{\R}{\mathbb{R}}

\newcommand{\cA}{\mathcal{A}}
\newcommand{\cB}{\mathcal{B}}

\newcommand{\cD}{\mathcal{D}}
\newcommand{\cE}{\mathcal{E}}
\newcommand{\cF}{\mathcal{F}}

\newcommand{\cI}{\mathcal{I}}

\newcommand{\cN}{\mathcal{N}}

\newcommand{\cR}{\mathcal{R}}
\newcommand{\cS}{\mathcal{S}}
\newcommand{\cT}{\mathcal{T}}

\newcommand{\cY}{\mathcal{Y}}

\newcommand{\ccA}{\mathscr{A}}
\newcommand{\ccB}{\mathscr{B}}

\newcommand{\ccP}{\mathscr{P}}

\newcommand{\bd}{\bm{d}}
\newcommand{\bS}{\bm{S}}
\newcommand{\by}{\bm{y}}
\newcommand{\bY}{\bm{Y}}
\newcommand{\bz}{\bm{z}}

\renewcommand{\d}{\mathrm{d}}

\renewcommand{\epsilon}{\varepsilon}
\renewcommand{\theta}{\vartheta}
\renewcommand{\rho}{\varrho}
\renewcommand{\phi}{\varphi}

\newcommand{\mo}[1]{\mathopen{#1}}
\newcommand{\mc}[1]{\mathclose{#1}}

\newcommand{\1}{\mathbbm{1}}

\DeclarePairedDelimiter{\rounds}{{\mathopen{(}}}{{\mathclose{)}}}
\DeclarePairedDelimiter{\bracks}{\lbrack}{\rbrack}
\DeclarePairedDelimiter{\braces}{\lbrace}{\rbrace}

\DeclarePairedDelimiter{\abs}{\lvert}{\rvert}
\DeclarePairedDelimiter{\norm}{\lVert}{\rVert}

\newcommand{\Lped}[3]{L_{#3}^{#1}\rounds{#2}}

\newcommand{\xto}[2]{\xrightarrow[#2]{#1}} 

\DeclareMathOperator{\bP}{\mathbf{P}}

\DeclareMathOperator{\bE}{\mathbf{E}}

\DeclareMathOperator{\blambda}{\bm{\lambda}}

\DeclareMathOperator{\supp}{supp}
\DeclareMathOperator{\tr}{tr}
\newcommand{\TV}{\mathrm{TV}}

\DeclareMathOperator*{\essinf}{ess\,inf}

\newenvironment{system}
{\left\lbrace\begin{array}{@{}l@{}}}
{\end{array}\right.}

\theoremstyle{plain}
    \newtheorem{theorem}{Theorem}[section]
	
    \newtheorem{proposition}{Proposition}[section]

	\newtheorem{corollary}{Corollary}[section]
	
	\newtheorem{algorithm}{Algorithm} 

\theoremstyle{definition}

	\newtheorem{example}{Example}[section]

\theoremstyle{remark}
	\newtheorem{remark}{Remark}[section]
	
\begin{document}


\title{\textbf{\Large{Compatibility-based Bayesian inference:
asymptotic limits, variational principles, and Wasserstein geometry}}}
\author{\textbf{Marco Tarsia $\bm{\cdot}$ Antonietta Mira $\bm{\cdot}$ Daniele Cassani}}
\date{}



\maketitle

{\let\thefootnote\relax\footnote{
    	\textbf{Marco Tarsia} \newline
	Department of Economics \newline
	University of Insubria, Varese 21100, Italy \newline
	E\:\!-mail: \mail{marcotarsia90@gmail.com} \newline
	Web page: \href{https://www.uninsubria.it/hpp/marco.tarsia}{\texttt{https://www.uninsubria.it/hpp/marco.tarsia}}
    \par\smallskip\noindent
	\textbf{Antonietta Mira} \newline
    Department of Science and High Technology \newline
    University of Insubria, Como 22100, Italy \newline
    E\:\!-mail: \mail{antonietta.mira@uninsubria.it} \newline
    Web page: \href{https://www.uninsubria.it/hpp/antonietta.mira}{\texttt{https://www.uninsubria.it/hpp/antonietta.mira}}
    \par\smallskip\noindent
    \textbf{Daniele Cassani} \newline
	Department of Science and High Technology \newline 
	University of Insubria, Como 22100, Italy \newline
	RISM--Riemann International School of Mathematics \newline
	Via G.\,B. Vico 46, Varese 21100, Italy \newline
	E\:\!-mail: \mail{daniele.cassani@uninsubria.it} \newline
	Web page: \href{https://www.uninsubria.it/hpp/daniele.cassani}{\texttt{https://www.uninsubria.it/hpp/daniele.cassani}}
	}
}

\thispagestyle{empty}

\par\noindent\rule{\linewidth}{0.8pt}\par

\bigskip\bigskip


\begin{description}

    \item[Abstract.]
    Approximate Bayesian Computation (ABC) replaces intractable likelihood evaluations with comparisons between observed and simulated data.
    We formulate discrepancy-based inference through nonnegative compatibility weights.
    Prior-integrated convergence of rescaled weights gives total variation convergence as the tolerance vanishes.
    For a fixed positive tolerance and increasing sample size, we give a general representation-based formulation of the conditional-prior limit known from ABC and coarsened Bayesian inference, and establish convergence in expected total variation under stabilization in probability and a prior-null boundary condition.
    We also derive acceptance-persistence bounds under partial stabilization.
    The Gibbs variational identity yields convergence of minimum values and stability of approximate minimizers; the fixed-tolerance limit is characterized separately by constrained relative-entropy minimization.
    For additional penalty-based functionals of the limiting Wasserstein discrepancy, we prove an exact reduction to distributions over model laws and, under a geodesic assumption, a representation by distances from a Wasserstein ball.
    These penalty-based updates are distinguished from the original fixed-tolerance ABC limit and illustrated through a Gaussian comparison.


    \item[Key words.]
    Approximate Bayesian Computation, compatibility filtering, entropy-regularized inference, optimal transport, Wasserstein distance.

    \item[2020 Mathematics Subject Classification.]
    49J45, 49Q22, 60B05, 62F15, 62F30.

\end{description}

\bigskip\bigskip

\noindent\hfil\rule{0.5\textwidth}{.8pt}


\section*{Introduction}
\label{sec:introduction}

Bayesian inference for generative models often requires comparing observed and simulated data because the likelihood cannot be evaluated.
Approximate Bayesian Computation (ABC) implements this comparison through a discrepancy and an acceptance tolerance; see, e.g.,~\cite{beaumont-zhang-balding_2002,marin-pudlo-robert-ryder_2012,sisson-fan-beaumont_2019}, and the applications in~\cite{dutta-mira-onnela_2018,ebert-dutta-mengersen-mira-ruggeri-wu_2018,chen-mira-onnela_2020,goyal-onnela_2023}.
Distributional comparisons include kernel embeddings, Kullback--Leibler discrepancies, and energy statistics~\cite{park-jitkrittum-sejdinovic_2016,jiang-wu-wong_2018,nguyen-arbel-lu-forbes_2019}.
Wasserstein discrepancies additionally use the geometry of the observation space~\cite{bernton-jacob-gerber-robert_2019,panaretos-zemel_2019,weed-bach_2019,lin-ho-chen-cuturi-jordan_2020}.

Two asymptotic questions must be distinguished: whether an ABC posterior recovers a likelihood-based posterior as the tolerance vanishes at fixed sample size, and what happens as the sample size grows at a fixed positive tolerance.
The second limit need not concentrate on a single parameter, even when the limiting discrepancy has a unique minimizer.
This phenomenon is already established in the literature.
Miller and Dunson~\cite{miller-dunson_2019} obtain a prior restricted to a discrepancy neighborhood as the deterministic-tolerance case of coarsened Bayesian inference.
Bernton et al.~\cite[Section~3.2]{bernton-jacob-gerber-robert_2019} explicitly state the corresponding conditional-prior limit for Wasserstein ABC and discuss its lack of point concentration.
These fixed-tolerance results complement the joint sample-size and tolerance asymptotics studied in~\cite{frazier-martin-robert-rousseau_2018,frazier-robert-rousseau_2020}.

Our purpose is to formulate these mechanisms through measurable compatibility weights and to relate them to entropy minimization and Wasserstein geometry.
The normalization principle in Proposition~\ref{prop:threshold-limit} is classical: $L^1(\bm{\pi})$ convergence of rescaled nonnegative weights implies total variation convergence of their normalized measures.
We make its hypotheses explicit and apply it to likelihood recovery, including a discussion of summary statistics and permutation-invariant empirical representations.

For the large-sample regime, Theorem~\ref{thm:compatibility-filtering} works with general observed and parameter-dependent simulated representations in a common pseudo-metric space of probability measures.
It assumes stabilization in probability rather than an almost-sure coupling of all sample sizes.
Outside a prior-null discrepancy boundary, the acceptance weights converge to an indicator in prior-integrated, observation-averaged $L^1$.
Corollary~\ref{cor:compatibility-filtering} then gives convergence in expected total variation to the conditional prior, together with an explicit normalization bound.
The result thus provides a precise formulation under these hypotheses of the established conditional-prior mechanism, rather than a new claim that fixed tolerances prevent point concentration.
We distinguish the prior-essential discrepancy threshold from the model-wide infimum, and derive quantitative acceptance-persistence bounds under weaker stabilization assumptions and local discrepancy-growth conditions.

The variational analysis uses the classical Gibbs--Donsker--Varadhan identity~\cite{donsker-varadhan_1975,dupuis-ellis_1997}.
For a compatibility weight $w_\epsilon$, the posterior minimizes its expected negative logarithm plus relative entropy from the prior.
Theorem~\ref{thm:inf_l-KL-asy} combines this identity with Pinsker's inequality to prove convergence of minimum values and stability of approximate minimizers under the weight convergence above.
It does not assert convergence of the objective functionals themselves or a $\Gamma$-convergence theorem.
For the fixed-tolerance large-sample limit, the exact variational problem is instead the minimization of relative entropy subject to assigning probability one to the compatibility region.

We also study additional penalty-based functionals of the limiting discrepancy.
These define generalized Bayesian updates~\cite{bissiri-holmes-walker_2016}, related in perspective to coarsening, scoring-rule inference, and discrepancy-based robust updating~\cite{miller-dunson_2019,pacchiardi-khoo-dutta_2024,matsubara-knoblauch-briol-oates_2022}.
Their penalties are specified separately from the finite-sample ABC acceptance rule; their Gibbs minimizers generally differ from the conditional-prior limit.
This distinction preserves the probabilistic interpretation of ABC acceptance probabilities, including model-error interpretations such as~\cite{wilkinson_2013}.

When the limiting discrepancy is Wasserstein, the parameter-to-distribution map sends a posterior candidate to a probability measure on the space of model laws.
Proposition~\ref{prop:KL-Phi} uses the data-processing inequality and an attaining lift to reduce relative-entropy minimization exactly to this space.
Corollary~\ref{cor:KL-Phi} applies the reduction to the additional penalty-based problem.
For power penalties and a geodesic Wasserstein space, Theorem~\ref{thm:R_phi,p^star-representation} expresses the excess-discrepancy cost as an expected power of the distance from the ball whose radius is the prior-essential minimum discrepancy.
This radius is distinct from the chosen ABC tolerance.

The Gaussian example illustrates fixed-tolerance filtering through reproducible finite-sample calculations and compares the limiting ABC posterior with penalty-based Gibbs updates.
Finite-state, mixed discrete--continuous, and deterministic models illustrate the scope of the distributional representation.

Section~\ref{sec:ABC} specifies the probabilistic framework.
Section~\ref{sec:threshold-limit} treats vanishing tolerances, and Section~\ref{sec:large-sample} treats fixed-tolerance filtering and partial stabilization.
Section~\ref{sec:variational-transport} develops the variational and transport constructions and the examples.
Section~\ref{sec:conclusions} concludes.

\section{Approximate Bayesian Computation}
\label{sec:ABC} 

This section introduces the observation and parameter spaces, the underlying probabilistic model, and the discrepancy-based construction of the corresponding posterior distributions.




\subsection{Probabilistic foundations of ABC}
\label{subsec:ABC-framework} 


Fix a probability space $\rounds{ \Omega \:\!, \;\!\! \cA , \;\!\! \bP }$.
For any topological space $X$\;\!\!, we denote by $\cB(X)$ its Borel $\sigma$-\:\!algebra and by $\ccP(X)$ the set of probability measures defined on $\cB(X)$. 
Given another topological space $Y$\;\!\!, $\ccB(X ,\;\!\! Y)$ stands for the class of Borel measurable mappings from $X$ to $Y$\;\!\!.

Let
\[
\cY \subseteq \R^{m_{\;\!\! \cY}} \:\!\! 
\equiv
\R^{m_{\;\!\! \cY} \times \:\! 1} \;\!\!
\]
be the observation space, assumed separable under the metric $\bd_{\;\!\! \cY}$, where 
\[
m_{\:\!\! \cY} \;\!\! \in \N^{\:\! *} \:\!\! \coloneqq \N \setminus \Set{\! 0 \!} 
\]
gives the dimension of a generic data point $\by \in \cY$\;\!\!.
We also specify the sample size $n \in \N^{\:\! *}$\:\!\!.


The observed random sample is represented by the random vector
\[
\bY^{1:\:\!n} \;\!\!
\equiv \;\!\!
\bracks[\big]{ \bY^1 \:\!\!, \dots ,\;\!\! \bY^n } \;\!\!
\colon
\rounds{ \Omega \:\!, \;\!\! \cA }
\longrightarrow
\big{(} \cY^{\:\! n} \:\!\!,\;\!\! \cB(\cY^{\:\! n}) \big{)} ,
\]
with distribution $\bP^{\:\! 1:\:\!n} \;\!\! \in \ccP \big{(} \cY^{\:\! n} \big{)}$.

For a fixed realization $\omega\in\Omega$, write $\by^{1:n}\coloneqq\bY^{1:n}(\omega)$ and $\by^i\coloneqq\bY^i(\omega)$.
The dependence on $\omega$ is suppressed when the observed sample is treated as given.


We next introduce the parameter space
\[
\Theta \subseteq \R^{m_{\:\! \Theta}} \;\!\!,
\]
governing the data-generating mechanism and assumed separable under the metric $\bd_{\;\! \Theta}$, where 
\[
m_{\:\! \Theta} \;\!\! \in \N^{\:\! *} \:\!\! 
\]
is the dimension of the parameter vectors $\bm{\theta} \in \Theta$.
The Borel structures and all finite products are understood with respect to these metric topologies.
Separability ensures that the Borel sigma-algebra of each finite product is the product of the corresponding Borel sigma-algebras, as needed for the kernel-integral measurability statements below.
This is an explicit metric-space assumption; set-theoretic inclusion in a Euclidean space alone does not imply it.
For probability measures $\bm{\rho},\bm{\eta}\in\ccP(\Theta)$, we use the convention
\[
\|\bm{\rho}-\bm{\eta}\|_{\:\!\TV}
\coloneqq \;\!\! \sup_{B\in\cB(\Theta)}|\bm{\rho}[B]-\bm{\eta}[B]|.
\]
If both measures have densities relative to a common dominating measure, this equals one half of the $L^1$\:\!\! distance between their densities.





We further consider a prior distribution
\[
\bm{\pi} \in \ccP(\Theta)
\]
on the parameter space $\Theta$, together with a Markov kernel from $\Theta$ to $\cY^{\:\! n}$\:\!\!,
\[
K_n \;\!\!
\colon
\Theta
\times
\cB \big{(} \cY^{\:\! n} \big{)} \;\!\!
\longrightarrow
\mo{[} 0 \:\!,\:\!\! 1 \mc{]} \:\! ,
\]
defining the associated family of sampling distributions: for every $\bm{\theta} \in \Theta$, the probability measure
\[
\bm{\mu}_{\bm{\theta}}^{\:\! n}
\coloneqq
K_n (\bm{\theta} , \bm{\cdot} \,)
\in
\ccP \big{(} \cY^{\:\! n} \big{)}
\]
represents the distribution of a pseudo-sample generated under the parameter $\bm{\theta}$.

In particular, for every $D \in \cB \big{(} \cY^{\:\! n} \big{)}$, the mapping
\[
\bm{\theta}
\longmapsto
\bm{\mu}_{\bm{\theta}}^{\:\! n} \bracks{\;\!\! D }
\]
belongs to $\ccB \big{(} \Theta \:\!,\;\!\! \mo{[} 0 \:\!,\:\!\! 1 \mc{]} \big{)}$ and is therefore $\bm{\pi}$-\:\!integrable (in the Lebesgue sense).


The first structural assumption underlying Approximate Bayesian Computation can be stated and will henceforth be maintained throughout the paper.

\begin{description}

\item[\textbf{($\ccA$-\:\!$1$)}] \hypertarget{item:A-1}{\label{item:A-1} --- \emph{Simulatability}. 
The statistical model $\Set{\! \bm{\mu}_{\bm{\theta}}^{\:\! n} \!}_{\bm{\theta} \:\! \in \:\! \Theta}$\;\!\! induced by $K_n$\;\!\! is generative: for every $\bm{\theta} \in \Theta$, one can simulate arbitrarily many independent pseudo-samples
\[
\bY_{\;\!\! \bm{\theta}}^{1:\:\!n} \;\!\!
\equiv \;\!\!
\bracks[\big]{ \bY_{\;\!\! \bm{\theta}}^1 , \dots ,\;\!\! \bY_{\;\!\! \bm{\theta}}^n } \;\!\!
\in \cY^{\:\! n} \:\!\! ,
\]
namely random vectors distributed according to $\bm{\mu}_{\bm{\theta}}^{\:\! n}$; that is,
\[
\bY_{\;\!\! \bm{\theta}}^{1:\:\!n} \!
\sim
\bm{\mu}_{\bm{\theta}}^{\:\! n} \:\! .
\]
} 

\end{description}

For every $\bm{\theta} \in \Theta$, any realization of the random pseudo-sample $\bY_{\;\!\! \bm{\theta}}^{1:\:\!n}$\;\!\! is written as
\[
\by_{\bm{\theta}}^{1:\:\!n} \:\!\!
\in \cY^{\:\! n} \:\!\! .
\]




We emphasize that the observed sample vector
$\by^{1:\:\!n}$\:\!\! constitutes the fixed empirical realization arising from the underlying statistical experiment and is regarded as given, whereas any pseudo-sample $\by_{\bm{\theta}}^{1:\:\!n}$\;\!\! is a realization of simulated data generated under a prescribed parameter value $\bm{\theta} \in \Theta$.


This distinction lies at the heart of ABC: in the absence of direct likelihood evaluation, inference compares observed and simulated data through a discrepancy measure
\[
\cD \colon \cY^{\:\! n} \:\!\! \times \cY^{\:\! n} \:\!\! \longrightarrow \R_{\:\! +} \;\!\! \coloneqq \mo{[} 0 \:\!,\:\!\! \infty \mc{[} \:\! ,
\]
defined as a continuous pseudo-metric on $\cY^{\:\! n}$\:\!\!, satisfying all the standard axioms of a metric except that distinct elements may have zero distance; see, e.g.,~\cite{park-jitkrittum-sejdinovic_2016,bernton-jacob-gerber-robert_2019,nguyen-arbel-lu-forbes_2019}.

To this end, for any observed sample $\by^{1:\:\!n} \:\!\! \in \cY^{\:\! n}$\:\!\! and radius $r \in \R_{\:\! +}$, we introduce the sample-dependent closed subset $D_{\;\!\! r \;\!\!} \big{(} \by^{1:\:\!n} \big{)}$ of the observation space $\cY^{\:\! n}$\:\!\! containing all elements whose deviation from $\by^{1:\:\!n}$\:\!\! does not exceed $r$ under the discrepancy measure $\cD$; namely,
\[
D_{\;\!\! r \;\!\!} \big{(} \by^{1:\:\!n} \big{)}
\coloneqq
\cD \big{(} \by^{1:\:\!n} \:\!\!, \bm{\cdot} \, \big{)}^{\;\!\! -1 \;\!\!} \big{(} \mo{[} 0 \:\! , \;\!\! r \mc{]} \big{)}
=
\Set{\! \widetilde{\by}^{\, 1:\:\!n} \;\!\! \in \cY^{\:\! n} | \cD \big{(} \by^{1:\:\!n} \:\!\!, \widetilde{\by}^{\, 1:\:\!n} \big{)} \le r \!} \;\!\!
\in \cB(\cY^{\:\! n}) \:\! .
\]
The family $\Set{\;\!\!\! D_{\;\!\! r \;\!\!} \big{(} \by^{1:\:\!n} \big{)} \;\!\!\!}_{r \:\! \ge \:\! 0}$\;\!\! is increasing in $r$ and, as $r \downarrow 0$,
\begin{equation*}
D_{\;\!\! r \;\!\!} \big{(} \by^{1:\:\!n} \big{)} \;\!\!
\downarrow
D_0 \big{(} \by^{1:\:\!n} \big{)} \;\!\!
\equiv
\cD \big{(} \by^{1:\:\!n} \:\!\!, \bm{\cdot} \, \big{)}^{\;\!\! -1 \;\!\!} \big{(}\! \Set{\! 0 \!} \!\big{)} .
\end{equation*}

For the classical ABC construction in this section and its application in Subsection~\ref{subsec:ABC-classical}, we impose the following nondegeneracy assumption.

\begin{description}

\item[\textbf{($\ccA$-\:\!$2$)}] \hypertarget{item:A-2}{\label{item:A-2} --- \emph{Nondegeneracy}. 
There exists $\epsilon^{\;\!\! \circ} \:\!\! \equiv \epsilon_{\;\!\! n}^{\;\!\! \circ} \;\!\! \in \mo{]}\:\! 0 \:\!,\:\!\! \infty \mc{[}$ such that, for every $\epsilon \in \mo{]}\:\! 0 \:\!,\;\!\! \epsilon^{\;\!\! \circ} \;\!\! \mc{[}$---\:\!hereafter referred to as an ABC threshold\:\!---\:\!and for $\bP^{\:\! 1:\:\!n}$\;\!\!-\:\!almost every $\by^{1:\:\!n} \:\!\! \in \cY^{\:\! n}$\:\!\!,
\begin{equation}
\label{eq:C_epsilon^n} 
C_{\:\!\! \epsilon}^{\;\! n} \big{(} \by^{1:\:\!n} \big{)}
\coloneqq \;\!\!
\int\nolimits_{\Theta} \bm{\mu}_{\bm{\theta}}^{\:\! n} \bracks[\big]{\;\!\! D_{\! \epsilon \;\!\!} \big{(} \by^{1:\:\!n} \big{)} \;\!\!} \;\! \bm{\pi}(\d\bm{\theta}) > 0 \:\! .
\end{equation}
} 

\end{description}

By monotonicity in $\epsilon$, the exceptional null set in~\hyperlink{item:A-2}{\customlabel{($\ccA$-\:\!$2$)}} can be chosen independently of the threshold, ensuring that the stated positivity holds simultaneously for all $\epsilon \in \mo{]}\:\! 0 \:\!,\;\!\! \epsilon^{\;\!\! \circ} \;\!\! \mc{[}$.

For each fixed $\by^{1:\:\!n} \:\!\! \in \cY^{\:\! n}$\:\!\!, ABC threshold $\epsilon$, and $\bm{\theta} \in \Theta$, the quantity $\bm{\mu}_{\bm{\theta}}^{\:\! n} \bracks[\big]{\;\!\! D_{\! \epsilon \;\!\!} \big{(} \by^{1:\:\!n} \big{)} \;\!\!}$ represents the acceptance probability under the model with parameter $\bm{\theta}$, namely the probability that simulated data generated under $\bm{\theta}$ fall within the $\epsilon$-neighborhood $D_{\! \epsilon \;\!\!} \big{(} \by^{1:\:\!n} \big{)}$ of the empirical realization $\by^{1:\:\!n}$.

As a function of the parameter $\bm{\theta}$, the acceptance probability
\begin{equation}
\label{eq:mu_theta^n} 
\bm{\theta}
\longmapsto
\bm{\mu}_{\bm{\theta}}^{\:\! n} \bracks[\big]{\;\!\! D_{\! \epsilon \;\!\!} \big{(} \by^{1:\:\!n} \big{)} \;\!\!}
\end{equation}
quantifies the compatibility between the model under $\bm{\theta}$ and the observed sample, and may be interpreted as a surrogate likelihood, in which the classical pointwise comparison between observed and model-implied data is replaced by a tolerance-based notion of proximity; see, e.g.,~\cite{marin-pudlo-robert-ryder_2012,sisson-fan-beaumont_2019}.
Large values indicate that the model under $\bm{\theta}$ assigns high probability to simulated realizations close to the observed sample, while small values reflect poor compatibility with the observed data.





\subsection{ABC posterior distribution}
\label{subsec:ABC-posterior} 

We now turn to the probabilistic formulation of the ABC posterior distribution.
When direct likelihood evaluation is unavailable or computationally intractable, but simulation from the statistical model remains feasible according to~\hyperlink{item:A-1}{\customlabel{($\ccA$-\:\!$1$)}}, approximate Bayesian inference can be performed by reweighting the prior distribution $\bm{\pi}$ according to the compatibility between observed and simulated data.



Given an ABC threshold $\epsilon$, each parameter value $\bm{\theta} \in \Theta$ is weighted according to the corresponding acceptance probability.
The resulting normalized measure defines an approximate Bayesian update and is formalized below as the ABC posterior distribution, in which the role traditionally played by the likelihood is taken precisely by the compatibility map~\eqref{eq:mu_theta^n}; see, e.g.,~\cite{marin-pudlo-robert-ryder_2012,sisson-fan-beaumont_2019}.


The following scheme summarizes the classical ABC rejection algorithm; see, e.g.,~\cite{beaumont-zhang-balding_2002,marin-pudlo-robert-ryder_2012}.

\begin{algorithm}[ABC rejection algorithm]
\label{alg:ABC-rej} 

Let $\by^{1:\:\!n} \:\!\! \in \cY^{\:\! n}$\:\!\! be an observed sample.


\begin{description} 

    \item[\textbf{Step 1}~\textmd{(Tolerance specification)}.]
    Choose an ABC threshold $\epsilon \in \mo{]}\:\! 0 \:\!,\;\!\! \epsilon^{\;\!\! \circ} \;\!\! \mc{[}$.

    \item[\textbf{Step 2}~\textmd{(Sampling and simulation)}.]
    Draw a parameter $\bm{\theta} \in \Theta$ from the prior distribution $\bm{\pi}$ and generate a pseudo-sample $\by_{\bm{\theta}}^{1:\:\!n}$\:\!\! according to $\bm{\mu}_{\bm{\theta}}^{\:\! n}$. 

    \item[\textbf{Step 3}~\textmd{(Acceptance rule)}.]
    Retain $\bm{\theta}$ if and only if $\by_{\bm{\theta}}^{1:\:\!n} \:\!\! \in D_{\! \epsilon \;\!\!} \big{(} \by^{1:\:\!n} \big{)}$.

\end{description}

\end{algorithm}




For every ABC threshold $\epsilon$, consider the ABC normalizing function defined on $\cY^{\:\! n}$\:\!\! through~\eqref{eq:C_epsilon^n}, 
\[
\by^{1:\:\!n} \:\!\!
\longmapsto
C_{\:\!\! \epsilon}^{\;\! n} \big{(} \by^{1:\:\!n} \big{)} ,
\quad
\by^{1:\:\!n} \:\!\! \in \cY^{\:\! n} \:\!\! .
\]
By~\hyperlink{item:A-2}{\customlabel{($\ccA$-\:\!$2$)}}, $C_{\:\!\! \epsilon}^{\;\! n} ( \, \bm{\cdot} \, ) > 0$ $\bP^{\:\! 1:\:\!n}$\;\!\!-\:\!a.s.
Moreover, $C_{\:\!\! \epsilon}^{\;\! n} ( \, \bm{\cdot} \, ) \in \ccB \big{(} \cY^{\:\! n} \:\!\!,\;\!\! \mo{[} 0 \:\!,\:\!\! 1 \mc{]} \big{)}$ since the Markov-kernel property of $K_n$, together with the assumed continuity of $\cD$\;\!\!, ensures the Borel measurability of the mapping 
\[
\big{(} \bm{\theta} ,\;\!\! \by^{1:\:\!n} \big{)}
\longmapsto
\bm{\mu}_{\bm{\theta}}^{\:\! n} \bracks[\big]{\;\!\! D_{\! \epsilon \;\!\!} \big{(} \by^{1:\:\!n} \big{)} \;\!\!} ,
\quad
\big{(} \bm{\theta} ,\;\!\! \by^{1:\:\!n} \big{)} \in \Theta \times \cY^{\:\! n} \:\!\! .
\]

It follows that the function $\bm{\Pi}_{\:\!\! \epsilon}^n \colon \cY^{\:\! n} \:\!\! \times \cB(\Theta) \to \mo{[} 0 \:\!,\:\!\! 1 \mc{]}$ given, for every $\by^{1:\:\!n} \:\!\! \in \cY^{\:\! n}$\:\!\! and $B \in \cB(\Theta)$, by
\begin{equation*}
\bm{\Pi}_{\:\!\! \epsilon}^n \big{(}\:\! \by^{1:\:\!n} \:\!\!,\;\!\! B \:\!\big{)} \;\!\!
\coloneqq
\begin{cases}
\displaystyle
\frac{1}{C_{\:\!\! \epsilon}^{\;\! n} \big{(} \by^{1:\:\!n} \big{)}} \;\!\!
\int\nolimits_{\;\!\! B} \bm{\mu}_{\bm{\theta}}^{\:\! n} \bracks[\big]{\;\!\! D_{\! \epsilon \;\!\!} \big{(} \by^{1:\:\!n} \big{)} \;\!\!} \;\! \bm{\pi}(\d\bm{\theta}) ,
& \text{if $C_{\:\!\! \epsilon}^{\;\! n} \big{(} \by^{1:\:\!n} \big{)} > 0$\:\!} ,
\\[2.25ex]
\:\! \bm{\pi} \bracks{ B \:\!} ,
& \text{otherwise} \:\! ,
\end{cases}
\end{equation*}
is a Markov kernel from $\cY^{\:\! n}$\:\!\! to $\Theta$, and
\[
\bm{\Pi}_{\:\!\! \epsilon}^n \big{(}\:\! \bY^{1:\:\!n} \:\!\!,\:\! \bm{\cdot} \:\!\,\big{)} \;\!\!
\]
is the random ABC posterior distribution induced by the observed sample.

For every fixed $\by^{1:\:\!n} \:\!\! \in \cY^{\:\! n}$\:\!\! with
$C_{\:\!\! \epsilon}^{\;\! n} \big{(} \by^{1:\:\!n} \big{)} > 0$, we write the corresponding posterior distribution as
\begin{equation}
\label{eq:pi_epsilon^n} 
\bm{\pi}_{\:\!\! \epsilon}^{\:\! n}
\equiv
\bm{\pi}_{\:\!\! \epsilon}^{\:\! \by^{1:\:\!n}} \:\!\!\!
\coloneqq
\bm{\Pi}_{\:\!\! \epsilon}^n \big{(}\:\! \by^{1:\:\!n} \:\!\!,\:\! \bm{\cdot} \:\!\,\big{)} \;\!\!
\in \ccP(\Theta) \:\! .
\end{equation}
By construction,
\[
\bm{\pi}_{\:\!\! \epsilon}^{\:\! n} \;\!\! \ll \;\!\! \bm{\pi}
\]
and its Radon--Nikodym derivative with respect to $\bm{\pi}$ is proportional to the mapping~\eqref{eq:mu_theta^n}.

The foregoing construction shows that repeated independent executions of steps~\textbf{2}\:\!–\:\!\textbf{3} of Algorithm~\ref{alg:ABC-rej} generate accepted parameter values distributed according to $\bm{\pi}_{\:\!\! \epsilon}^{\:\! n}$.

\section{On the vanishing-threshold limit}
\label{sec:threshold-limit} 

In this section, we study the asymptotic behavior of the posterior distributions~\eqref{eq:pi_epsilon^n} as the tolerance parameter $\epsilon$ vanishes, with the sample size $n$ held fixed.
We introduce a general formulation in terms of nonnegative compatibility weights and establish sufficient conditions for convergence in total variation.

The complementary large-sample regime is investigated separately in Section~\ref{sec:large-sample}.

\subsection{A general convergence principle}
\label{subsec:convergence-principle} 


Consider a general family of nonnegative measurable functions
\[
w_{\:\!\! \epsilon} \;\!\! \in \ccB(\Theta \:\!,\;\!\! \R_{\:\! +}) ,
\]
indexed by $\epsilon \in \mo{]}\:\! 0 \:\!,\:\!\! \infty \mc{[}$, intended to quantify the compatibility between the model under a parameter value $\bm{\theta} \in \Theta$ and the observed sample $\by^{1:\:\!n}$\:\!\!, in the same spirit as~\eqref{eq:mu_theta^n}; we refer to these functions as compatibility weights.
We assume that, for every sufficiently small $\epsilon$, $w_{\:\!\! \epsilon}$\;\!\! is $\bm{\pi}$-\:\!integrable and not equal to zero $\bm{\pi}$-\:\!almost everywhere.
Accordingly, these functions induce the family $\Set{\! \bm{\pi}_{\:\!\! \epsilon} \!}_{\;\!\! \epsilon}$\;\!\! in $\ccP(\Theta)$ defined by
\begin{equation}
\label{eq:pi_epsilon} 
\bm{\pi}_{\:\!\! \epsilon} \bracks{ B \:\!}
\coloneqq
\frac{
\int\nolimits_B w_{\:\!\! \epsilon}(\bm{\theta}) \;\! \bm{\pi}(\d\bm{\theta})
}
{
\int\nolimits_{\:\! \Theta} w_{\:\!\! \epsilon}(\bm{\theta}') \;\! \bm{\pi}(\d\bm{\theta}')
} ,
\quad
B \in \cB(\Theta) .
\end{equation}
This formulation encompasses the classical ABC posterior~\eqref{eq:pi_epsilon^n} as a special case; see, e.g.,~\cite{beaumont-zhang-balding_2002,marin-pudlo-robert-ryder_2012,sisson-fan-beaumont_2019}.

We next establish a fundamental convergence principle governed by the asymptotic behavior of the compatibility weights $w_{\:\!\! \epsilon}$\;\!\! in the vanishing-threshold regime.

\begin{proposition}
\label{prop:threshold-limit} 

Assume that there exist a family $\Set{\! \kappa_{\;\!\! \epsilon} \!}_{\;\!\! \epsilon}$\;\!\! of strictly positive constants and a function $\widebar{w} \in \ccB(\Theta \:\!,\;\!\! \R_{\:\! +})$, $\bm{\pi}$-\:\!integrable and not $\bm{\pi}$-\:\!almost everywhere equal to zero, such that
\begin{equation}
\label{eq:bar-w} 
\lim_{\epsilon \:\! \downarrow \:\! 0} \
\int\nolimits_{\Theta} \
\abs*{\:\!
\frac{1}{\kappa_{\;\!\! \epsilon}} \:\! w_{\:\!\! \epsilon}(\bm{\theta}) - \widebar{w}(\bm{\theta})
\:\!} \;\! \bm{\pi}(\d\bm{\theta})
= 0 \:\! .
\end{equation}
Then, as $\epsilon \downarrow 0$, the family of posterior distributions $\Set{\! \bm{\pi}_{\:\!\! \epsilon} \!}_{\;\!\! \epsilon}$ introduced in~\eqref{eq:pi_epsilon} converges in total variation to the probability measure $\widebar{\bm{\pi}} \in \ccP(\Theta)$ defined by
\begin{equation}
\label{eq:bar-pi} 
\widebar{\bm{\pi}} \bracks{ B \:\!}
\coloneqq
\frac{
\int\nolimits_B \widebar{w}(\bm{\theta}) \;\! \bm{\pi}(\d\bm{\theta})
}
{
\int\nolimits_{\:\! \Theta} \widebar{w}(\bm{\theta}') \;\! \bm{\pi}(\d\bm{\theta}')
} ,
\quad
B \in \cB(\Theta) ,
\end{equation}
that is,
\begin{equation}
\label{eq:threshold-limit} 
\lim_{\epsilon \:\! \downarrow \:\! 0} \;\!
\sup_{B \:\! \in \:\! \cB(\Theta)} \;\!
\abs[\big]{\:\! \bm{\pi}_{\:\!\! \epsilon} \bracks{ B \:\!} - \widebar{\bm{\pi}} \bracks{ B \:\!} \:\!}
= 0 \:\! .
\end{equation}

\end{proposition}

\begin{proof}

For every $B \in \cB(\Theta)$, by the definitions above and the triangle inequality for integrals,
\[
\begin{split}
&\hphantom{= \lvert}  \:\! \abs[\big]{\:\! \bm{\pi}_{\:\!\! \epsilon} \bracks{ B \:\!} - \widebar{\bm{\pi}} \bracks{ B \:\!} \:\!} \\[2.25ex]
&= \:\!
\abs*{\, \frac{\int\nolimits_B \;\! \bracks*{ \frac{1}{\kappa_{\;\!\! \epsilon}} \:\! w_{\:\!\! \epsilon}(\bm{\theta}) - \widebar{w}(\bm{\theta}) } \;\! \bm{\pi}(\d\bm{\theta})}{\int\nolimits_{\:\! \Theta} \frac{1}{\kappa_{\;\!\! \epsilon}} \:\! w_{\:\!\! \epsilon}(\bm{\theta}') \;\! \bm{\pi}(\d\bm{\theta}')} \:\! + \:\! \bracks*{ \frac{1}{\int\nolimits_{\:\! \Theta} \frac{1}{\kappa_{\;\!\! \epsilon}} \:\! w_{\:\!\! \epsilon}(\bm{\theta}') \;\! \bm{\pi}(\d\bm{\theta}')} - \frac{1}{\int\nolimits_{\:\! \Theta} \widebar{w}(\bm{\theta}') \;\! \bm{\pi}(\d\bm{\theta}')} } \:\!\! \int\nolimits_{\;\!\! B} \widebar{w}(\bm{\theta}) \;\! \bm{\pi}(\d\bm{\theta}) \,} \\[1.25ex]
&\le \:\!
\frac{\int\nolimits_{\:\! \Theta} \;\! \abs*{\:\! \frac{1}{\kappa_{\;\!\! \epsilon}} \:\! w_{\:\!\! \epsilon}(\bm{\theta}) - \widebar{w}(\bm{\theta}) \:\!} \;\! \bm{\pi}(\d\bm{\theta})}{\int\nolimits_{\:\! \Theta} \frac{1}{\kappa_{\;\!\! \epsilon}} \:\! w_{\:\!\! \epsilon}(\bm{\theta}') \;\! \bm{\pi}(\d\bm{\theta}')} \:\! + \:\! \abs*{\;\! \frac{1}{\int\nolimits_{\:\! \Theta} \frac{1}{\kappa_{\;\!\! \epsilon}} \:\! w_{\:\!\! \epsilon}(\bm{\theta}') \;\! \bm{\pi}(\d\bm{\theta}')} - \frac{1}{\int\nolimits_{\:\! \Theta} \widebar{w}(\bm{\theta}') \;\! \bm{\pi}(\d\bm{\theta}')} \;\!} \int\nolimits_{\Theta} \widebar{w}(\bm{\theta}) \;\! \bm{\pi}(\d\bm{\theta}) .
\end{split}
\]
Both terms on the right-hand side converge to zero as $\epsilon \downarrow 0$ directly from~\eqref{eq:bar-w}.
Since the bound is independent of $B \in \cB(\Theta)$, taking the supremum proves~\eqref{eq:threshold-limit}. \qedhere

\end{proof}

Within the present framework, Proposition~\ref{prop:threshold-limit} highlights that the asymptotic behavior of the induced posterior measures $\bm{\pi}_{\:\!\! \epsilon}$\;\!\! is essentially governed by the corresponding rescaled compatibility weights $w_{\:\!\! \epsilon}/\kappa_{\;\!\! \epsilon}$.
In particular, convergence in integral form of such weights\:\!---\:\!$\Lped{\;\!\! 1 \:\!\!}{\bm{\pi}}{}$-\:\!convergence\:\!---\:\!as specified in~\eqref{eq:bar-w} yields convergence in total variation toward a well-defined limiting posterior distribution $\widebar{\bm{\pi}}$.

This compatibility-weight formulation also underlies the variational and transport formulations of approximate Bayesian inference developed in Section~\ref{sec:variational-transport}.


\subsection{Bayesian posterior recovery in classical ABC}
\label{subsec:ABC-classical} 

The preceding approach applies to a broad class of approximate Bayesian constructions.
We first specialize it to the ABC setting based on the full sample, where the compatibility weights are induced by acceptance probabilities over shrinking neighborhoods of the observed realization; see, e.g.,~\cite{beaumont-zhang-balding_2002,marin-pudlo-robert-ryder_2012,sisson-fan-beaumont_2019}.

For this purpose, let $\by^{1:\:\!n} \:\!\! \in \cY^{\:\! n}$\:\!\! be the observed sample and assume that there exists a $\sigma$\;\!\!-\:\!finite measure $\blambda^{\;\!\! n}$\:\!\! on $\cB \big{(} \cY^{\:\! n} \big{)}$ such that, for $\bm{\pi}$-\:\!almost every $\bm{\theta} \in \Theta$, the probability measure $\bm{\mu}_{\bm{\theta}}^{\:\! n}$ is absolutely continuous with respect to $\blambda^{\;\!\! n}$\:\!\!, with Radon--Nikodym derivative
\[
f_{\;\!\! \bm{\theta}}^{\;\! n}
\coloneqq
\frac{\d\:\!\bm{\mu}_{\bm{\theta}}^{\:\! n}}{\d\:\!\!\blambda^{\;\!\! n} \:\!\!}
\]
admitting a version, still denoted by $f_{\;\!\! \bm{\theta}}^{\;\! n}$,
which is continuous at $\by^{1:\:\!n}$\:\!\!.
The generality of the dominating measure $\blambda^{\;\!\! n}$\;\!\! allows for continuous, discrete, lower-dimensional, and mixed data, including settings described by Lebesgue, counting, Hausdorff, or suitable combinations of such measures.


Moreover, define
\begin{equation}
\label{eq:f_theta^n} 
\widebar{w}(\bm{\theta})
\coloneqq
f_{\;\!\! \bm{\theta}}^{\;\! n} \big{(} \by^{1:\:\!n} \big{)}
\end{equation}
and suppose that it belongs to $\ccB \big{(} \Theta \:\!,\;\!\! \R_{\:\! +} \big{)}$ and is not $\bm{\pi}$-\:\!almost everywhere equal to zero.
The function~\eqref{eq:f_theta^n} then provides a version of the likelihood associated with $\by^{1:\:\!n}$\:\!\!.
Accordingly, for every ABC threshold $\epsilon$, set
\[
\kappa_{\;\!\! \epsilon} \;\!\!
\coloneqq
\blambda^{\;\!\! n} \bracks[\big]{\;\!\! D_{\! \epsilon \;\!\!} \big{(} \by^{1:\:\!n} \big{)} \;\!\!}
\] 
and, for $\bm{\pi}$-\:\!almost every $\bm{\theta} \in \Theta$,
\[
w_{\:\!\! \epsilon}(\bm{\theta})
\coloneqq
\bm{\mu}_{\bm{\theta}}^{\:\! n} \bracks[\big]{\;\!\! D_{\! \epsilon \;\!\!} \big{(} \by^{1:\:\!n} \big{)} \;\!\!}
= \;\!\!
\int\nolimits_{\;\!\! D_{\! \epsilon \;\!\!}(\by^{1:\:\!n})} \! f_{\;\!\! \bm{\theta}}^{\;\! n} \big{(} \bz^{1:\:\!n} \big{)} \:\! \blambda^{\;\!\! n} \! \big{(} \d\bz^{1:\:\!n} \big{)} .
\]

Assume, in addition, that $\cD$ is a genuine metric compatible with the topology of $\cY^{\:\! n}$\:\!\!, and that
\[
0 < \kappa_{\;\!\! \epsilon} \;\!\! < \infty
\]
for every sufficiently small $\epsilon$.
Then, for $\bm{\pi}$-\:\!almost every $\bm{\theta} \in \Theta$, continuity of $f_{\;\!\! \bm{\theta}}^{\;\! n}$ at $\by^{1:\:\!n}$\:\!\! yields, as $\epsilon \downarrow 0$,
\[
\abs*{\:\!
\frac{1}{\kappa_{\;\!\! \epsilon}} \:\! w_{\:\!\! \epsilon}(\bm{\theta}) - \widebar{w}(\bm{\theta})
\:\!}
\le \:\!\!
\sup_{\bz^{1:\:\!n} \in \:\! D_{\! \epsilon \;\!\!}(\by^{1:\:\!n})} \;\!
\abs*{\:\!
f_{\;\!\! \bm{\theta}}^{\;\! n} \big{(} \bz^{1:\:\!n} \big{)} \;\!\! - \;\!\! f_{\;\!\! \bm{\theta}}^{\;\! n} \big{(} \by^{1:\:\!n} \big{)}
}
\longrightarrow
0 \:\! .
\]



Therefore, if there exists a function $g \in \ccB(\Theta \:\!,\;\!\! \R_{\:\! +})$ which is $\bm{\pi}$-\:\!integrable and satisfies, for every sufficiently small $\epsilon$ and for $\bm{\pi}$-\:\!almost every $\bm{\theta} \in \Theta$,
\[
\frac{1}{\kappa_{\;\!\! \epsilon}} \:\! w_{\:\!\! \epsilon}(\bm{\theta})
\le
g(\bm{\theta}) \:\! ,
\]
then the integral convergence~\eqref{eq:bar-w} follows from the dominated convergence theorem.
Consequently, Proposition~\ref{prop:threshold-limit} implies convergence in total variation, as $\epsilon \downarrow 0$, of the ABC posterior $\bm{\pi}_{\:\!\! \epsilon} \;\!\! \equiv \bm{\pi}_{\:\!\! \epsilon}^{\:\! n}$ toward the probability measure~\eqref{eq:bar-pi}, here explicitly given by
\begin{equation*}
\widebar{\bm{\pi}} \bracks{ B \:\!}
=
\frac{
\int\nolimits_{B \:\!} f_{\;\!\! \bm{\theta}}^{\;\! n} \big{(} \by^{1:\:\!n} \big{)} \;\! \bm{\pi}(\d\bm{\theta})
}
{
\int\nolimits_{\:\! \Theta \:\!} f_{\bm{\theta}'}^{\;\! n} \big{(} \by^{1:\:\!n} \big{)} \;\! \bm{\pi}(\d\bm{\theta}')
} ,
\quad
B \in \cB(\Theta) .
\end{equation*}
From a measure-theoretic perspective, $\widebar{\bm{\pi}}$ provides a version of the usual likelihood-based posterior distribution given the full observed sample $\by^{1:\:\!n}$\:\!\!, thus recovering the classical Bayesian update.
Its density with respect to $\bm{\pi}$ is proportional to the likelihood function~\eqref{eq:f_theta^n}; see, e.g.,~\cite{parthasarathy_1967}. 

\begin{remark}

The preceding argument recovers the full-sample posterior under the stated metric and density assumptions.
More generally, let $\rounds{\cS ,\;\!\! \bd_{\cS}}$ be a metric space and suppose that the discrepancy is induced by a measurable summary map $S^{\:\! n} \;\!\! \colon \cY^{\:\! n} \:\!\! \to \cS$ through
\[
\cD \big{(} \by^{1:\:\!n} \:\!\!,\;\!\! \bz^{1:\:\!n} \big{)} \;\!\!
=
\bd_{\cS} \big{(}\:\! S^{\:\! n} \big{(} \by^{1:\:\!n} \big{)} ,\;\!\! S^{\:\! n} \big{(} \bz^{1:\:\!n} \big{)} \;\!\!\big{)} .
\]
Under analogous regularity and domination conditions for the summary distributions, the limiting posterior is determined by the likelihood associated with $S^{\:\! n}(\by^{1:\:\!n})$\:\!\!; see, e.g.,~\cite{fearnhead-prangle_2012,marin-pudlo-robert-ryder_2012}.
Hence, it coincides with the posterior distribution given the observed summary statistic, which generally differs from the posterior given the full sample; the two coincide, in particular, when $S^{\:\! n}$\:\!\! is sufficient.

When $\cD$ is only a pseudo-metric, observations in the equivalence class $\cD \big{(} \by^{1:\:\!n} \:\!\!, \bm{\cdot} \, \big{)}^{\;\!\! -1 \;\!\!} \big{(}\! \Set{\! 0 \!} \!\big{)}$ cannot be distinguished as $\epsilon \downarrow 0$, so that the limiting information need not identify the single realization $\by^{1:\:\!n}$\:\!\!.

An important case is $\cD(\by^{1:n},\bz^{1:n})=W_p(\widehat{\bm{\nu}}_n(\by^{1:n}),\widehat{\bm{\nu}}_n(\bz^{1:n}))$, where $\widehat{\bm{\nu}}_n$ is the empirical measure.
On ordered samples this is a pseudo-metric: its zero-discrepancy class consists of the permutations of the observed sample.
For exchangeable models admitting permutation-invariant likelihood versions, this loss of ordering need not lose information about the parameter.
Under the continuity and domination hypotheses for shrinking permutation neighborhoods, the full-sample posterior is recovered; see~\cite[Proposition~2(a)]{bernton-jacob-gerber-robert_2019}.
The genuine-metric argument above should therefore be distinguished from this permutation-invariant extension.

\end{remark}

The previous analysis shows that classical ABC convergence arises as a particular instance of Proposition~\ref{prop:threshold-limit}.
The fundamental mechanism underlying posterior convergence is not tied to the structure of ABC neighborhoods, but rather to integral control of suitably normalized compatibility weights.

Within this framework, uniform local bounds on the likelihood are not intrinsically required, as the necessary integral convergence may instead be obtained through an integrated domination condition with respect to the prior distribution.
Thus, regularity assumptions appearing in classical ABC convergence analyses (see, e.g.,~\cite{marin-pudlo-robert-ryder_2012,sisson-fan-beaumont_2019}) may be viewed as sufficient means of establishing such control, rather than as intrinsic requirements of the general convergence principle.





\section{The large-sample regime}
\label{sec:large-sample} 

We now study the established fixed-tolerance conditional-prior mechanism~\cite{miller-dunson_2019,bernton-jacob-gerber-robert_2019} under explicit assumptions on general probabilistic representations.
The observation and parameter spaces, prior, sampling kernels, and simulation interpretation from Section~\ref{sec:ABC} are retained.
However, the acceptance discrepancy is now the jointly measurable, possibly parameter-dependent map defined below; it need not be the continuous, parameter-independent pseudo-metric on ordered samples used in Section~\ref{sec:ABC}.
Assumption~\hyperlink{item:A-2}{\customlabel{($\ccA$-\:\!$2$)}} is not imposed on this generalized construction.
Its posterior is defined by normalization when possible and by a measurable probability-kernel extension otherwise.
We prove that the zero-normalization event has probability tending to zero in the supercritical regime.

Theorem~\ref{thm:compatibility-filtering} establishes prior-integrated, observation-averaged $L^1$\:\!\! convergence of the acceptance weights under stabilization in probability.
Corollary~\ref{cor:compatibility-filtering} strengthens the posterior conclusion to expected total variation convergence.
The remaining results quantify acceptance persistence under partial stabilization and describe compatibility regions near discrepancy minimizers.

\subsection{Probabilistic representations and asymptotic discrepancies}
\label{subsec:probabilistic-representations} 

The construction below separates the sampling mechanism of the ABC procedure from the asymptotic discrepancy geometry governing acceptance.
Rather than working with finite-dimensional sample vectors, we represent observed and simulated samples by probability measures.
This representation is well suited to the large-sample regime, since the ambient spaces vary with $n$, while these probability measures take values in a common space.
Acceptance is then expressed in terms of the discrepancy between them.

We equip $\ccP(\cY)$ with a generic pseudo-metric
\begin{equation}
\label{eq:T} 
\cT \colon \ccP(\cY) \times \ccP(\cY) \longrightarrow \R_{\:\! +} .
\end{equation}
This level of generality naturally encompasses transport-based discrepancies, whose role in statistical inference will be further explored in Section~\ref{sec:variational-transport}; see, e.g.,~\cite{bernton-jacob-gerber-robert_2019}.

For every sample size $n \in \N^{\:\! *}$\:\!\!, let
\[
\bm{\nu}_{n} \;\!\! 
\colon
\cY^{\:\! n} \:\!\!
\longrightarrow
\ccP(\cY)
\]
and, for every $\bm{\theta} \in \Theta$, 
\[
\bm{\nu}_{\bm{\theta} \:\!\!, n} \;\!\! 
\colon
\cY^{\:\! n} \:\!\!
\longrightarrow
\ccP(\cY)
\]
denote probabilistic representation maps for the observed and simulated samples, respectively.
More explicitly, any realization $\by^{1:\:\!n} \:\!\! \in \cY^{\:\! n}$\:\!\! of the observed random sample and any pseudo-sample $\by_{\bm{\theta}}^{1:\:\!n} \:\!\! \in \cY^{\:\! n}$ generated under parameter $\bm{\theta}$ are respectively associated with the probability measures
\[
\bm{\nu}_{n} \big{(} \by^{1:\:\!n} \big{)} \;\!\! \in \ccP(\cY)
\]
and
\[
\bm{\nu}_{\bm{\theta} \:\!\!, n} \big{(} \by_{\bm{\theta}}^{1:\:\!n} \big{)} \;\!\! \in \ccP(\cY) .
\]

A natural choice is given by empirical measures, although the subsequent analysis does not rely on this specific representation; see, e.g.,~\cite{bernton-jacob-gerber-robert_2019,panaretos-zemel_2019}.
What is essential is that observed and simulated samples are mapped into the common space $\ccP(\cY)$, so that the asymptotic behavior of the ABC acceptance mechanism can be characterized in terms of the geometry induced by $\cT$\:\!\!.

\begin{remark}

Since $\cT$ is only assumed to be a pseudo-metric, the topology induced on $\ccP(\cY)$ need not be Hausdorff (equivalently, $T_{\:\! 2}$).
In particular, distinct probability measures on $\cY$ may be topologically indistinguishable whenever their discrepancy with respect to $\cT$ vanishes.

\end{remark}

We next assume that, for every $n \in \N^{\:\! *}$\:\!\!, the mapping 
\begin{equation}
\label{eq:Delta_theta,n} 
\big{(} \bm{\theta} , \by^{1:\:\!n} \:\!\!, \by_{\bm{\theta}}^{1:\:\!n} \big{)} \;\!\!
\longmapsto
\Delta_{\;\! \bm{\theta} \:\!\!, n} \big{(} \by^{1:\:\!n} \:\!\!, \by_{\bm{\theta}}^{1:\:\!n} \big{)} \;\!\!
\coloneqq
\cT \big{(} \bm{\nu}_{n} \big{(} \by^{1:\:\!n} \big{)} , \;\!\! \bm{\nu}_{\bm{\theta} \:\!\!, n} \big{(} \by_{\bm{\theta}}^{1:\:\!n} \big{)} \;\!\! \big{)} \;\!\!
\end{equation}
belongs to $\ccB \big{(} \Theta \times \cY^{\:\! n} \:\!\! \times \cY^{\:\! n} \:\!\!, \R_{\:\! +} \big{)}$. 


The large-sample behavior of the ABC compatibility selection mechanism is ultimately governed by the discrepancy geometry emerging in $\ccP(\cY)$, rather than by finite-dimensional sample representations.
To identify this asymptotic structure, we assume the existence of a probability measure
\begin{equation}
\label{eq:mu^star} 
\bm{\mu}^{\;\!\! \star} \:\!\! \in \ccP(\cY)
\end{equation}
such that, for every $n \in \N^{\:\! *}$\:\!\!, the mapping
\begin{equation*}
\omega
\longmapsto
\cT_n^{\;\! \star} \;\!\! (\omega)
\coloneqq
\cT \big{(}\:\! \bm{\mu}^{\;\!\! \star} \:\!\!,\;\!\! \bm{\nu}_{n} \big{(} \bY^{1:\:\!n}(\omega) \big{)} \;\!\!\big{)}
\end{equation*}
is $\cA$-measurable and converges to zero in $\bP$-\:\!probability as $n \to \infty$, namely 
\begin{equation}
\label{eq:lim-T_n^star} 
\cT_n^{\;\! \star} \:\!\!
\xto{\bP}{\;\!\!\! n \:\! \to \:\! \infty \;\!\!\!}
\;\! 0 \:\! .
\end{equation}
This condition is fulfilled, in particular, whenever $\cT_n^{\;\! \star}$\:\!\! converges to zero $\bP$-\:\!almost surely.

Similarly, for every $\bm{\theta} \in \Theta$, we assume the existence of a probability measure 
\begin{equation}
\label{eq:eta_theta} 
\bm{\nu}_{\bm{\theta}} \;\!\! \in \ccP(\cY) 
\end{equation}
such that the function
\begin{equation}
\label{eq:T_theta^star} 
\bm{\theta}
\longmapsto
\cT_{\bm{\theta}}^{\;\! \star} \:\!\!
\coloneqq
\cT \big{(}\:\! \bm{\mu}^{\;\!\! \star} \:\!\!,\;\!\! \bm{\nu}_{\bm{\theta}} \big{)}
\end{equation}
belongs to $\ccB(\Theta \:\!,\;\!\! \R_{\:\! +})$ and, for $\bm{\pi}$-\:\!almost every $\bm{\theta} \in \Theta$ and every $n \in \N^{\:\! *}$\:\!\!, the mapping
\[
\by_{\bm{\theta}}^{1:\:\!n}
\longmapsto
\cT \big{(} \bm{\nu}_{\bm{\theta} \:\!\!, n} \big{(} \by_{\bm{\theta}}^{1:\:\!n} \big{)} , \;\!\! \bm{\nu}_{\bm{\theta}} \big{)}
\]
belongs to $\ccB \big{(} \cY^{\:\! n} \:\!\!,\;\!\! \R_{\:\! +} \big{)}$ and, for every $\xi > 0$,
\begin{equation}
\label{eq:lim-T_theta,n} 
\lim_{n \:\! \to \:\! \infty} \:\!
\bm{\mu}_{\bm{\theta}}^{\:\! n} \bracks*{\;\!\!
\Set{\:\!\!
\by_{\bm{\theta}}^{1:\:\!n} \:\!\! \in \cY^{\:\! n}
|
\cT \big{(} \bm{\nu}_{\bm{\theta} \:\!\!, n} \big{(} \by_{\bm{\theta}}^{1:\:\!n} \big{)} , \;\!\! \bm{\nu}_{\bm{\theta}} \big{)} \;\!\!
>
\xi
\:\!\!}
\;\!\!} \;\!\!
=
0 \:\! .
\end{equation}
This amounts to convergence in probability under the sequence of pseudo-sample laws~$\Set{\! \bm{\mu}_{\bm{\theta}}^{\:\! n} \!}_{n \:\! \in \:\! \N^{\:\! *}}$\;\!\!.

The quantity $\cT_{\bm{\theta}}^{\;\! \star}$\:\!\! represents the asymptotic discrepancy level of $\bm{\theta}$, comparing its model law with the limiting observational law.
As $\bm{\theta}$ varies over $\Theta$, the function~\eqref{eq:T_theta^star} induces the asymptotic compatibility geometry of the model class and determines the relative discrepancy levels across parameter values.

Although the limiting probability measures need not be uniquely determined under a pseudo-metric, the asymptotic discrepancy levels defined in~\eqref{eq:T_theta^star} are uniquely determined up to $\bm{\pi}$-\:\!null sets.


The asymptotic discrepancy function~\eqref{eq:T_theta^star} determines two thresholds, according to whether it is assessed relative to the prior or over the entire model class: the $\bm{\pi}$-\:\!essential compatibility threshold
\begin{equation}
\label{eq:epsilon_pi^star} 
\epsilon_{\;\!\! \bm{\pi}}^{\;\!\! \star} \;\!\!
\coloneqq
\essinf_{\bm{\theta} \:\! \in \:\! \Theta} \;\!
\cT_{\bm{\theta}}^{\;\! \star} \:\!\! ,
\end{equation}
where the essential infimum is taken with respect to $\bm{\pi}$, and the model\:\!-wide compatibility threshold
\begin{equation}
\label{eq:epsilon_Theta^star} 
\epsilon_{\Theta}^{\;\!\! \star} \;\!\!
\coloneqq \;\!\!
\inf_{\bm{\theta} \:\! \in \:\! \Theta}
\cT_{\bm{\theta}}^{\;\! \star} \:\!\!
\le
\epsilon_{\;\!\! \bm{\pi}}^{\;\!\! \star} \:\! .
\end{equation}

The threshold $\epsilon_{\;\!\! \bm{\pi}}^{\;\!\! \star}$\;\!\! represents the essential lower bound of the asymptotic discrepancy under $\bm{\pi}$, while $\epsilon_{\Theta}^{\;\!\! \star}$\;\!\! describes the best level available over the model class.
Strict inequality in~\eqref{eq:epsilon_Theta^star} thus reflects a restriction induced by the prior relative to the discrepancy geometry.
In particular, $\epsilon_{\Theta}^{\;\!\! \star} \;\!\! = 0$ means that the limiting observational law can be matched arbitrarily closely within the model class; if $\epsilon_{\Theta}^{\;\!\! \star} \;\!\! > 0$, an intrinsic asymptotic discrepancy remains, as may occur under model misspecification, structural incompatibility, or limited expressive capacity of the simulated probabilistic family; see, e.g.,~\cite{frazier-robert-rousseau_2020}.

For a fixed effective ABC tolerance $\tau \in \mo{]}\:\! 0 \:\!,\:\!\! \infty \mc{[}$, the asymptotic behavior of the compatibility selection mechanism depends on its relation to the threshold $\epsilon_{\;\!\! \bm{\pi}}^{\;\!\! \star}$.
We restrict attention to the supercritical regime
\begin{equation}
\label{eq:tau-supercritical} 
\tau
>
\epsilon_{\;\!\! \bm{\pi}}^{\;\!\! \star} \:\! .
\end{equation}
The boundary case $\tau = \epsilon_{\;\!\! \bm{\pi}}^{\;\!\! \star}$\;\!\! may depend on prior mass near this threshold and on sampling fluctuations, requiring a local asymptotic analysis.
In the subcritical regime $\tau < \epsilon_{\;\!\! \bm{\pi}}^{\;\!\! \star}$\;\!\!, acceptance probabilities vanish asymptotically for $\bm{\pi}$-\:\!almost every parameter.
The normalized posterior may then be governed by their relative decay rates, calling for large-deviation or Laplace-type analysis.

\begin{remark}

The assumptions introduced above can be satisfied in a broad range of standard Approximate Bayesian Computation settings.
This covers discrepancy constructions induced by Wasserstein distances and kernel-based metrics such as the maximum mean discrepancy; see, e.g.,~\cite{sriperumbudur-gretton-fukumizu-scholkopf-lanckriet_2010,park-jitkrittum-sejdinovic_2016,bernton-jacob-gerber-robert_2019}.
Functionals built from measurable summary statistics may also be accommodated, as may other probability metrics, including Prokhorov or total variation discrepancies, whenever the chosen probabilistic representations satisfy the required measurability and convergence properties.

The probability measures $\bm{\mu}^{\;\!\! \star}$\:\!\! and $\bm{\nu}_{\bm{\theta}}$\;\!\! may arise as limits under familiar probabilistic regimes, including i.i.d.~sampling, ergodic observations, and weak-convergence settings, provided that $\cT$ is compatible with the mode of convergence.
For empirical measures in Wasserstein distance, see, e.g.,~\cite{boissard-le.gouic_2014,fournier-guillin_2015}.

\end{remark}

\subsection{Acceptance geometry and posterior compatibility}
\label{subsec:acceptance-geometry} 

The discrepancy map~\eqref{eq:Delta_theta,n} determines the finite-sample acceptance mechanism.
Consistently with Algorithm~\ref{alg:ABC-rej}, a pseudo-sample $\by_{\bm{\theta}}^{1:\:\!n}$\:\!\! generated under $\bm{\theta} \in \Theta$ is retained whenever the discrepancy between $\bm{\nu}_{\bm{\theta} \:\!\!, n} \big{(} \by_{\bm{\theta}}^{1:\:\!n} \big{)}$\;\!\! and the observed representation $\bm{\nu}_{n} \big{(} \by^{1:\:\!n} \big{)}$\;\!\! does not exceed the tolerance $\tau$.

For every $n \in \N^{\:\! *}$\:\!\! and $\bm{\theta} \in \Theta$, define the data-dependent finite-sample acceptance region
\begin{equation}
\label{eq:A_theta^n} 
A_{\bm{\theta} \:\!\!, n \;\!\!} \big{(}\:\! \bY^{1:\:\!n} ;\;\!\! \tau \:\!\big{)} \;\!\!
\coloneqq
\Set{\:\!\!
\by_{\bm{\theta}}^{1:\:\!n} \:\!\! \in \cY^{\:\! n} \;\!\!
|
\Delta_{\;\! \bm{\theta} \:\!\!, n} \big{(}\:\! \bY^{1:\:\!n} \:\!\!, \by_{\bm{\theta}}^{1:\:\!n} \big{)} \;\!\!
\le
\tau
\:\!\!} \! .
\end{equation}
By measurability of $\Delta_{\;\! \bm{\theta} \:\!\!, n}$, for every $\omega \in \Omega$,
\[
A_{\bm{\theta} \:\!\!, n \;\!\!} \big{(}\:\! \bY^{1:\:\!n}(\omega) ;\;\!\! \tau \:\!\big{)} \;\!\!
\in
\cB \big{(} \cY^{\:\! n} \big{)} .
\]


Consistently with Section~\ref{sec:threshold-limit}, the finite-sample compatibility weight associated with $\bm{\theta} \in \Theta$ is defined as
\begin{equation}
\label{eq:w_tau^n} 
w_{\:\!\! \tau}^{\:\! n} \big{(} \bm{\theta} ;\;\!\! \bY^{1:\:\!n} \big{)} \;\!\!
\coloneqq
\bm{\mu}_{\bm{\theta}}^{\:\! n} \bracks[\big]{\;\!\! A_{\bm{\theta} \:\!\!, n \;\!\!} \big{(}\:\! \bY^{1:\:\!n} ;\;\!\! \tau \:\!\big{)} \;\!\!} .
\end{equation}
It represents the conditional acceptance probability of a pseudo-sample generated under $\bm{\theta}$, given the observed sample.
Under the preceding assumptions, joint Borel measurability holds for the mapping
\[
\big{(} \bm{\theta} ,\;\!\! \by^{1:\:\!n} \big{)}
\longmapsto
w_{\:\!\! \tau}^{\:\! n} \big{(} \bm{\theta} ;\;\!\! \by^{1:\:\!n} \big{)} .
\]


The corresponding normalization factor is
\begin{equation*}
C_{\;\!\! \tau}^{\;\! n} \big{(}\:\! \bY^{1:\:\!n} \big{)} \;\!\!
\coloneqq \;\!\!
\int\nolimits_{\Theta} w_{\:\!\! \tau}^{\:\! n} \big{(} \bm{\theta} ;\;\!\! \bY^{1:\:\!n} \big{)} \:\! \bm{\pi}(\d\bm{\theta}) .
\end{equation*}
For a given realization of the observed sample, it may vanish when the acceptance probability is zero for $\bm{\pi}$-\:\!almost every $\bm{\theta} \in \Theta$; whenever positive, the associated compatibility posterior is defined by 
\begin{equation}
\label{eq:Pi_tau^n} 
\bm{\Pi}_{\tau}^n \bracks*{ B ;\;\!\! \bY^{1:\:\!n} } \;\!\!
\coloneqq
\frac{1}{C_{\;\!\! \tau}^{\;\! n} \big{(}\:\! \bY^{1:\:\!n} \big{)} \;\!\!}
\int\nolimits_{\;\!\! B} w_{\:\!\! \tau}^{\:\! n} \big{(} \bm{\theta} ;\;\!\! \bY^{1:\:\!n} \big{)} \:\! \bm{\pi}(\d\bm{\theta}) \:\! ,
\quad
B \in \cB(\Theta) .
\end{equation}



In this representation, the posterior distribution~\eqref{eq:Pi_tau^n} may be interpreted as the prior reweighted by the compatibility weights~\eqref{eq:w_tau^n}.
Its large-sample behavior is therefore linked to the asymptotic discrepancy function~\eqref{eq:T_theta^star}, as made precise by the limiting behavior of these weights.


The limiting discrepancy profile~\eqref{eq:T_theta^star} induces a parameter-space counterpart of the finite-sample acceptance mechanism.
For the fixed tolerance $\tau$, define the asymptotic compatibility region
\begin{equation}
\label{eq:Theta_tau} 
\Theta_\tau \;\!\!
\coloneqq
\Set{\:\!\!
\bm{\theta} \in \Theta
|
\cT_{\bm{\theta}}^{\;\! \star} \:\!\!
<
\tau
\:\!\!} \;\!\!
\in \cB(\Theta) \:\! .
\end{equation}
In the supercritical regime~\eqref{eq:tau-supercritical}, definition~\eqref{eq:epsilon_pi^star} ensures that
\begin{equation*}
\bm{\pi} \bracks{ \Theta_{\tau \;\!\!} }
>
0 \:\! .
\end{equation*}
The next result identifies the asymptotic filtering behavior of the compatibility weights~\eqref{eq:w_tau^n} relative to $\Theta_\tau$.

\begin{theorem}
\label{thm:compatibility-filtering} 

Under the measurability assumptions of this section, suppose that
conditions~\eqref{eq:lim-T_n^star}, \eqref{eq:lim-T_theta,n},
and~\eqref{eq:tau-supercritical} hold, and that
\begin{equation}
\label{eq:tau-boundary} 
\bm{\pi}
\bracks*{\;\!\!
\Set{\:\!\!
\bm{\theta} \in \Theta
|
\cT_{\bm{\theta}}^{\;\! \star} \:\!\!
=
\tau
\:\!\!}
\;\!\!} \;\!\!
=
0 \:\! .
\end{equation}
Then the compatibility weights converge to the indicator of $\Theta_\tau$\;\!\! in the following integrated $L^{\;\!\! 1}$\:\!\! sense: 
\begin{equation}
\label{eq:w_tau^n-limit} 
\lim_{n \:\! \to \:\! \infty}
\bE_{\:\! \bP}
\bracks*{
\int\nolimits_{\Theta} \,
\abs*{\:\!
w_{\:\!\! \tau}^{\:\! n} \big{(} \bm{\theta} ;\;\!\! \bY^{1:\:\!n} \big{)} \;\!\! - \;\!\! \1_{\Theta_\tau}\:\!\!(\bm{\theta})
\:\!} \;\! \bm{\pi}(\d\bm{\theta})
} \:\!\!
=
0 \:\! .
\end{equation}

\end{theorem}

\begin{proof}

Fix $\bm{\theta} \in \Theta$ for which~\eqref{eq:lim-T_theta,n} holds, as is the case $\bm{\pi}$-\:\!almost everywhere.
For each $n \in \N^{\:\! *}$\:\!\!, consider the product probability measure $\bP \otimes \;\! \bm{\mu}_{\bm{\theta}}^{\:\! n}$ on $\big{(} \Omega \times \cY^{\:\! n} \:\!\!, \cA \otimes \cB(\cY^{\:\! n}) \;\!\!\big{)}$.

By the triangle inequality and symmetry of $\cT$\:\!\!, for every $\omega \in \Omega$ and $\by_{\bm{\theta}}^{1:\:\!n} \in \cY^{\:\! n}$\:\!\!,
\begin{equation}
\label{eq:triangle-inequality} 
\abs*{\:\!
\Delta_{\;\! \bm{\theta} \:\!\!, n} \big{(}\:\! \bY^{1:\:\!n}(\omega) ,\;\!\! \by_{\bm{\theta}}^{1:\:\!n} \big{)} \;\!\! - \cT_{\bm{\theta}}^{\;\! \star}
}
\le
\cT_n^{\;\! \star} \;\!\! (\omega)
+
\cT \big{(} \bm{\nu}_{\bm{\theta} \:\!\!, n} \big{(} \by_{\bm{\theta}}^{1:\:\!n} \big{)} , \;\!\! \bm{\nu}_{\bm{\theta}} \big{)} .
\end{equation}
Therefore, in view of~\eqref{eq:lim-T_n^star} and~\eqref{eq:lim-T_theta,n},
\[
\Delta_{\;\! \bm{\theta} \:\!\!, n} \big{(}\:\! \bY^{1:\:\!n}(\;\! \bm{\cdot} \;\!) \:\!, \text{\bf{-}} \, \big{)} \;\!\!
\longrightarrow
\cT_{\bm{\theta}}^{\;\! \star} \:\!\!
\]
in probability under the sequence of product measures~$\Set{\! \bP \otimes \;\! \bm{\mu}_{\bm{\theta}}^{\:\! n} \!}_{n \:\! \in \:\! \N^{\:\! *}}$\;\!\! as $n \to \infty$.

Condition~\eqref{eq:tau-boundary} allows us to further assume that $\cT_{\bm{\theta}}^{\;\! \star} \:\!\! \neq \tau$, so that the preceding convergence implies
\[
\1_{\:\!\! A_{\bm{\theta} \:\!\!, n \;\!\!} (\:\! \bY^{1:\:\!n}(\;\! \bm{\cdot} \;\!) ; \tau \:\!)}\;\!\!(\;\! \text{\bf{-}} \;\!)
\longrightarrow
\1_{\Theta_\tau}\:\!\!(\bm{\theta})
\]
in probability and, by boundedness, also in $L^{\;\!\! 1}$\:\!\! with respect to the same sequence of product measures.

In conclusion, combining~\eqref{eq:w_tau^n} with Tonelli's theorem and dominated convergence,
\begin{multline*}
\bE_{\:\! \bP}
\bracks*{
\int\nolimits_{\Theta} \,
\abs*{\:\!
w_{\:\!\! \tau}^{\:\! n} \big{(} \bm{\theta} ;\;\!\! \bY^{1:\:\!n} \big{)} \;\!\! - \;\!\! \1_{\Theta_\tau}\:\!\!(\bm{\theta})
\:\!} \;\! \bm{\pi}(\d\bm{\theta})
} \:\!\!
= \;\!\!
\int\nolimits_{\Theta}
\bE_{\:\! \bP}
\bracks*{\:\!
\abs*{\:\!
w_{\:\!\! \tau}^{\:\! n} \big{(} \bm{\theta} ;\;\!\! \bY^{1:\:\!n} \big{)} \;\!\! - \;\!\! \1_{\Theta_\tau}\:\!\!(\bm{\theta}) \:\!}
\:\!} \:\! \bm{\pi}(\d\bm{\theta}) \\
\le \;\!\!
\int\nolimits_{\Theta} \,
\bracks*{
\int\nolimits_{\:\! \Omega \;\! \times \cY^{\:\! n} \:\!\!} \,
\abs*{\:\!
\1_{\:\!\! A_{\bm{\theta} \:\!\!, n \;\!\!} (\:\! \bY^{1:\:\!n}(\omega) ; \tau \:\!)} \;\!\! \big{(} \by_{\bm{\theta}}^{1:\:\!n} \big{)} \;\!\! - \;\!\! \1_{\Theta_\tau}\:\!\!(\bm{\theta})
\:\!} \:\! \big{(}\;\!\! \bP \otimes \;\! \bm{\mu}_{\bm{\theta}}^{\:\! n} \big{)} \;\!\! \big{(} \d\omega \:\!,\:\!\! \d\by_{\bm{\theta}}^{1:\:\!n} \big{)}
\;\!\!} \:\! \bm{\pi}(\d\bm{\theta})
\longrightarrow
0
\end{multline*}
as $n \to \infty$, which proves~\eqref{eq:w_tau^n-limit} and hence the theorem. \qedhere

\end{proof}

The conditional-prior limit at fixed tolerance is established in Miller and Dunson~\cite{miller-dunson_2019} and, for Wasserstein ABC, in Bernton et al.~\cite[Section~3.2]{bernton-jacob-gerber-robert_2019}.
Here, an analogous filtering result is established for general probability-valued representations satisfying~\eqref{eq:lim-T_n^star} and~\eqref{eq:lim-T_theta,n}.
The integrated $L^1$\:\!\! convergence~\eqref{eq:w_tau^n-limit}, averaged over the observed data and the prior, provides quantitative control of the normalization factors and compatibility posteriors.

These consequences are made precise in the following corollary.

\begin{corollary}
\label{cor:compatibility-filtering}
Under the assumptions of Theorem~\ref{thm:compatibility-filtering},
\begin{equation}
\label{eq:C_tau^n-limit}
\lim_{n\to\infty}\bE_{\:\! \bP}
\left[\left|C_{\;\!\! \tau}^{\;\! n} \big{(}\:\! \bY^{1:\:\!n} \big{)}-\bm{\pi}[\Theta_\tau]\right|\right]=0\:\!,
\end{equation}
and
\begin{equation}
\label{eq:C_tau^n-zero}
\lim_{n\to\infty}\bP\left[C_{\;\!\! \tau}^{\;\! n} \big{(}\:\! \bY^{1:\:\!n} \big{)}=0\right]=0\:\!.
\end{equation}
Extend the compatibility posterior~\eqref{eq:Pi_tau^n} measurably to a probability kernel on the zero-normalization event, for example by assigning the prior there.
Then
\begin{equation}
\label{eq:posterior-mean-TV}
\lim_{n\to\infty}\bE_{\:\! \bP}
\left[\left\|\:\!\bm{\Pi}_\tau^n[\,\bm{\cdot} \,;\bY^{1:n}]
-\bm{\pi}[\,\bm{\cdot}\mid\Theta_\tau]\:\!\right\|_{\:\!\TV}\right]=0\:\!.
\end{equation}
In particular, the posteriors converge in total variation in $\bP$-probability, independently of the chosen measurable probability-kernel extension.
\end{corollary}

\begin{proof}
Set $a\coloneqq\bm{\pi}[\Theta_\tau]>0$ and
\[
D_n\coloneqq\int\nolimits_{\Theta}
\left|w_{\:\!\! \tau}^{\:\! n}(\bm{\theta};\bY^{1:n})
-\1_{\Theta_\tau}(\bm{\theta})\right|
\,\bm{\pi}(\d\bm{\theta}).
\]
By~\eqref{eq:w_tau^n-limit}, $\bE_{\:\! \bP}[D_n]\to0$.
Since
\[
|C_{\;\!\! \tau}^{\;\! n} \big{(}\:\! \bY^{1:\:\!n} \big{)}-a|\le D_n,
\]
we obtain~\eqref{eq:C_tau^n-limit}.
Moreover,
\[
\bP[C_{\;\!\! \tau}^{\;\! n} \big{(}\:\! \bY^{1:\:\!n} \big{)}=0]
\le \frac{\bE_{\:\! \bP}[D_n]}{a},
\]
which proves~\eqref{eq:C_tau^n-zero}.

On the event $\{C_{\;\!\! \tau}^{\;\! n} \big{(}\:\! \bY^{1:\:\!n} \big{)}>0\}$, the triangle inequality gives
\[
\left\|\:\!
\bm{\Pi}_\tau^n[\,\bm{\cdot}\,;\bY^{1:n}]
-\bm{\pi}[\,\bm{\cdot}\mid\Theta_\tau]
\:\!\right\|_{\:\!\TV}
\le
\frac{D_n+|C_{\;\!\! \tau}^{\;\! n} \big{(}\:\! \bY^{1:\:\!n} \big{)}-a|}{2a}
\le \frac{D_n}{a}.
\]
On the complementary event, $D_n=a$, so the same bound holds for any probability-kernel extension.
The total variation distance is measurable, since its defining supremum can be restricted to a countable algebra generating $\cB(\Theta)$.
Taking expectations proves~\eqref{eq:posterior-mean-TV}, and Markov's inequality yields convergence in probability.
\end{proof}

The proof of Theorem~\ref{thm:compatibility-filtering} also provides a direct geometric interpretation of the filtering mechanism.
Indeed, by~\eqref{eq:triangle-inequality}, for $\bm{\theta} \in \Theta_\tau$, $\omega \in \Omega$, and $\by_{\bm{\theta}}^{1:\:\!n} \in \cY^{\:\! n}$\:\!\!, the acceptance criterion~\eqref{eq:A_theta^n} is met whenever
\begin{equation}
\label{eq:triangle-inequality*} 
\cT_n^{\;\! \star} \;\!\! (\omega)
+
\cT \big{(} \bm{\nu}_{\bm{\theta} \:\!\!, n} \big{(} \by_{\bm{\theta}}^{1:\:\!n} \big{)} , \;\!\! \bm{\nu}_{\bm{\theta}} \big{)} \;\!\!
\le
\tau - \cT_{\bm{\theta}}^{\;\! \star} \:\!\! .
\end{equation}
The positive compatibility margin
\begin{equation*}
\tau_{\bm{\theta}}^{\star} \;\!\!
\coloneqq
\tau - \cT_{\bm{\theta}}^{\;\! \star} \:\!\! > 0 \:\! ,
\quad
\bm{\theta} \in \Theta_\tau ,
\end{equation*}
thus accommodates both observational and simulated finite-sample fluctuations.

\subsection{Compatibility under partial stabilization} 
\label{subsec:compatibility-persistence} 

The exact filtering results of the previous subsection rely on the asymptotic stabilization of both observed and simulated probabilistic representations, as specified in~\eqref{eq:lim-T_n^star} and~\eqref{eq:lim-T_theta,n}, respectively.
We now replace the latter condition with a weaker assumption, retaining the remaining hypotheses.
The compatibility margin in~\eqref{eq:triangle-inequality*} then provides quantitative bounds on asymptotic acceptance probabilities.

\begin{proposition}
\label{prop:partial-stabilization} 

Suppose there exist $B_{\;\!\! \tau} \;\!\! \in \cB(\Theta)$ with $\bm{\pi} \bracks{ B_{\;\!\! \tau \;\!\!} } > 0$ and $\bm{\pi} \bracks{ B_{\;\!\! \tau} \:\!\! \setminus \Theta_{\tau \;\!\!} } = 0$, and $\delta_{\;\!\! \tau} \;\!\! \in \mo{[} 0 \:\!,\:\!\! 1 \mc{[}$ such that, for $\bm{\pi}$-\:\!almost every $\bm{\theta} \in B_{\;\!\! \tau}$, some $\lambda_{\tau \:\!\!,\;\!\! \bm{\theta}} \;\!\! \in \mo{]}\:\! 0 \:\!,\:\!\! 1 \mc{[}$ satisfies
\begin{equation}
\label{eq:partial-stabilization} 
\limsup_{n \:\! \to \:\! \infty} \;\!
\bm{\mu}_{\bm{\theta}}^{\:\! n} \bracks*{\;\!\!
\Set{\:\!\!
\by_{\bm{\theta}}^{1:\:\!n} \:\!\! \in \cY^{\:\! n} \;\!\!
|
\cT \big{(} \bm{\nu}_{\bm{\theta} \:\!\!, n} \big{(} \by_{\bm{\theta}}^{1:\:\!n} \big{)} , \;\!\! \bm{\nu}_{\bm{\theta}} \big{)} \;\!\!
>
\lambda_{\tau \:\!\!,\;\!\! \bm{\theta}} \:\! \tau_{\bm{\theta}}^{\star}
\:\!\!}
\;\!\!} \;\!\!
\le
\delta_{\;\!\! \tau} .
\end{equation}
Then, for $\bm{\pi}$-\:\!almost every $\bm{\theta} \in B_{\;\!\! \tau}$,
\begin{equation}
\label{eq:partial-stabilization-weight}
\liminf_{n \:\! \to \:\! \infty} \;\!
\bE_{\:\! \bP}
\bracks*{
w_{\:\!\! \tau}^{\:\! n} \big{(} \bm{\theta} ;\;\!\! \bY^{1:\:\!n} \big{)}
} \;\!\!
\ge
1 - \delta_{\;\!\! \tau} .
\end{equation}

\end{proposition}

\begin{proof}

Fix $\bm{\theta} \in B_{\;\!\! \tau}$\;\!\! for which~\eqref{eq:partial-stabilization} holds.
By~\eqref{eq:triangle-inequality*}, for any $\omega \in \Omega$ and $\by_{\bm{\theta}}^{1:\:\!n} \:\!\! \in \cY^{\:\! n}$\:\!\!,
\[
\begin{system}
\:\! \cT_n^{\;\! \star} \;\!\! (\omega)
\le
\rounds{1 - \lambda_{\tau \:\!\!,\;\!\! \bm{\theta}}} \:\!
\tau_{\bm{\theta}}^{\star} \\[0.5ex]
\:\! \cT \big{(} \bm{\nu}_{\bm{\theta} \:\!\!, n} \big{(} \by_{\bm{\theta}}^{1:\:\!n} \big{)} , \;\!\! \bm{\nu}_{\bm{\theta}} \big{)} \;\!\!
\le
\lambda_{\tau \:\!\!,\;\!\! \bm{\theta}} \:\! \tau_{\bm{\theta}}^{\star}
\end{system} \:\!\!
\quad \Longrightarrow \quad
\by_{\bm{\theta}}^{1:\:\!n} \:\!\!
\in
A_{\bm{\theta} \:\!\!, n \;\!\!} \big{(}\:\! \bY^{1:\:\!n}(\omega) ;\;\!\! \tau \:\!\big{)} .
\]
Hence, by~\eqref{eq:w_tau^n} and the product structure,
\[
\bE_{\:\! \bP}
\bracks*{
w_{\:\!\! \tau}^{\:\! n} \big{(} \bm{\theta} ;\;\!\! \bY^{1:\:\!n} \big{)}
} \;\!\!
\ge
\bP \:\!
\bracks*{\:\!
\cT_n^{\;\! \star} \:\!\!
\le
\rounds{1 - \lambda_{\tau \:\!\!,\;\!\! \bm{\theta}}} \:\! \tau_{\bm{\theta}}^{\star}
\:\!}
\:\!
\bm{\mu}_{\bm{\theta}}^{\:\! n} \bracks*{\;\!\!
\Set{\:\!\!
\by_{\bm{\theta}}^{1:\:\!n} \:\!\! \in \cY^{\:\! n} \;\!\!
|
\cT \big{(} \bm{\nu}_{\bm{\theta} \:\!\!, n} \big{(} \by_{\bm{\theta}}^{1:\:\!n} \big{)} , \;\!\! \bm{\nu}_{\bm{\theta}} \big{)} \;\!\!
\le
\lambda_{\tau \:\!\!,\;\!\! \bm{\theta}} \:\! \tau_{\bm{\theta}}^{\star}
\:\!\!}
\;\!\!} .
\]
Since $\rounds{1 - \lambda_{\tau \:\!\!,\;\!\! \bm{\theta}}} \:\! \tau_{\bm{\theta}}^{\star} > 0$, relations~\eqref{eq:lim-T_n^star} and~\eqref{eq:partial-stabilization} imply that the first factor on the right-hand side converges to one and the lower limit of the second is at least $1 - \delta_{\;\!\! \tau}$, as claimed in~\eqref{eq:partial-stabilization-weight}. \qedhere

\end{proof}

An integrated version of Proposition~\ref{prop:partial-stabilization} follows directly from Tonelli's theorem and Fatou's lemma.

\begin{corollary}
\label{cor:partial-stabilization} 

Under the assumptions of Proposition~\ref{prop:partial-stabilization},
\begin{equation}
\label{eq:partial-stabilization-integrated} 
\liminf_{n \:\! \to \:\! \infty} \;\!
\bE_{\:\! \bP}
\bracks*{
\int\nolimits_{\;\!\! B_{\;\!\! \tau}} \! w_{\:\!\! \tau}^{\:\! n} \big{(} \bm{\theta} ;\;\!\! \bY^{1:\:\!n} \big{)} \:\! \bm{\pi}(\d\bm{\theta})
\;\!\!} \;\!\!
\ge
\rounds{1 - \delta_{\;\!\! \tau}} \:\! \bm{\pi} \bracks{ B_{\;\!\! \tau \;\!\!} } \:\! .
\end{equation}

\end{corollary}

Corollary~\ref{cor:partial-stabilization} shows that exact stabilization of the simulated representations is not necessary for compatible prior mass to persist asymptotically.
The residual level $\delta_{\;\!\! \tau}$\;\!\! quantifies the admissible loss, while $B_{\;\!\! \tau}$\;\!\! identifies the portion of the parameter space on which this control is retained.

We finally relate these persistence bounds to the geometry of the asymptotic discrepancy function~\eqref{eq:T_theta^star}.
We assume that the model\:\!-wide compatibility threshold~\eqref{eq:epsilon_Theta^star} is attained by at least one parameter, so that
\begin{equation}
\label{eq:epsilon_Theta^star-min} 
\epsilon_{\Theta}^{\;\!\! \star} \;\!\!
=
\min_{\bm{\theta} \:\! \in \:\! \Theta} \;\!
\cT_{\bm{\theta}}^{\;\! \star} \:\!\! ,
\end{equation}
and we then introduce the nonempty measurable set
\begin{equation}
\label{eq:Theta^star} 
\Theta^\star \:\!\!
\coloneqq \;\!\!
\Set{\:\!\!
\bm{\theta}^{\:\! \star} \! \in \Theta
|
\cT_{\bm{\theta}^{\:\! \star}}^{\:\! \star} \:\!\!
=
\epsilon_{\Theta}^{\;\!\! \star} \;\!\!
\:\!\!} \! ,
\end{equation}
the distance function to $\Theta^\star$\:\!\! defined, for every $\bm{\theta} \in \Theta$, by
\[
\bd_{\;\! \Theta^\star \:\!\!}(\bm{\theta})
\coloneqq \!
\inf_{\bm{\theta}^\star \:\!\! \in \:\! \Theta^\star \!}
\bd_{\;\! \Theta}(\bm{\theta} ,\:\!\! \bm{\theta}^{\:\! \star \;\!\!}) \:\! ,
\] 
and, for every $r \in \mo{]}\:\! 0 \:\!,\:\!\! \infty \mc{[}$, the open $r$-\:\!neighborhood of $\Theta^\star$\:\!\! in $\Theta$, 
\begin{equation}
\label{eq:Theta_r^star} 
\Theta_r^\star \;\!\!
\coloneqq \;\!\!
\Set{\:\!\!
\bm{\theta} \in \Theta
\:\!\ \big{|} \
\bd_{\;\! \Theta^\star \:\!\!}(\bm{\theta})
<
r
\:\!\!} \! .
\end{equation}


\begin{proposition}
\label{prop:compatibility-neighborhood} 

Suppose there exist $r^{\;\!\! \star} \:\!\! \in \mo{]}\:\! 0 \:\!,\:\!\! \infty \mc{[}$ and two non-decreasing functions
\[
\psi_{\;\!\! -} \:\!,\;\!\! \psi_{\;\!\! +} \;\!\!
\colon
\mo{[} 0 \:\!,\;\!\! r^{\;\!\! \star} \mc{]}
\longrightarrow
\R_{\:\! +} \;\!\!
\]
such that, for $\bm{\pi}$-\:\!almost every $\bm{\theta} \in \Theta_{r^{\;\!\! \star}}^\star$, 
\begin{equation}
\label{eq:compatibility-growth} 
\psi_{\;\!\! -} \big{(} \bd_{\;\! \Theta^\star \:\!\!}(\bm{\theta}) \;\!\!\big{)} \;\!\!
\le
\cT_{\bm{\theta}}^{\;\! \star} \:\!\!
-
\epsilon_{\Theta}^{\;\!\! \star} \;\!\!
\le
\psi_{\;\!\! +} \big{(} \bd_{\;\! \Theta^\star \:\!\!}(\bm{\theta}) \;\!\!\big{)} .
\end{equation}
Then, for any $r_{\;\!\! -} , r_{\;\!\! +} \;\!\! \in \mo{]} 0 \:\!,\;\!\! r^{\;\!\! \star} \;\!\!\mc{[}$ satisfying
\begin{equation}
\label{eq:compatibility-radii} 
\psi_{\;\!\! +}(r_{\;\!\! -})
<
\tau - \epsilon_{\Theta}^{\;\!\! \star} \;\!\!
\le
\psi_{\;\!\! -}(r_{\;\!\! +}) \:\! ,
\end{equation}
one has
\begin{equation}
\label{eq:compatibility-neighborhood}
\bm{\pi}
\bracks*{
\Theta_{r_{-}}^{\star} \!
\setminus
\Theta_\tau
\;\!\!} \;\!\!
=
\bm{\pi}
\bracks*{\;\!\!
\big{(}
\Theta_\tau \;\!\!
\cap
\Theta_{r^{\;\!\! \star}}^{\star} \:\!\!
\big{)} \;\!\!
\setminus
\Theta_{r_{+}}^{\star}
\;\!\!} \;\!\!
=
0 \:\! .
\end{equation}

\end{proposition}

\begin{proof}

For $\bm{\pi}$-\:\!almost every $\bm{\theta} \in \Theta_{r_{-}}^{\star} \:\!\! \subseteq \Theta_{r^{\;\!\! \star}}^{\star}$ (see~\eqref{eq:Theta_r^star}), in virtue of the upper bound in~\eqref{eq:compatibility-growth}, monotonicity of $\psi_{\;\!\! +}$, and the first inequality in~\eqref{eq:compatibility-radii}, we obtain
\[
\cT_{\bm{\theta}}^{\;\! \star} \:\!\!
\le
\epsilon_{\Theta}^{\;\!\! \star} \;\!\!
+
\psi_{\;\!\! +} \big{(} \bd_{\;\! \Theta^\star \:\!\!}(\bm{\theta}) \;\!\!\big{)} \;\!\!
\le
\epsilon_{\Theta}^{\;\!\! \star} \;\!\!
+
\psi_{\;\!\! +}(r_{\;\!\! -})
<
\tau ,
\]
hence $\bm{\theta} \in \Theta_\tau$\;\!\! by~\eqref{eq:Theta_tau}.
Similarly, for $\bm{\pi}$-\:\!almost every $\bm{\theta} \in \Theta_{r^{\;\!\! \star}}^{\star} \! \setminus \Theta_{r_{+}}^{\star}$, the lower bound in~\eqref{eq:compatibility-growth}, definition~\eqref{eq:Theta_r^star}, monotonicity of $\psi_{\;\!\! -}$, and the second inequality in~\eqref{eq:compatibility-radii} give
\[
\cT_{\bm{\theta}}^{\;\! \star} \:\!\!
\ge
\epsilon_{\Theta}^{\;\!\! \star} \;\!\!
+
\psi_{\;\!\! -} \big{(} \bd_{\;\! \Theta^\star \:\!\!}(\bm{\theta}) \;\!\!\big{)} \;\!\!
\ge
\epsilon_{\Theta}^{\;\!\! \star} \;\!\!
+
\psi_{\;\!\! -}(r_{\;\!\! +})
\ge
\tau ,
\]
which implies $\bm{\theta} \notin \Theta_\tau$.
The two conclusions together establish~\eqref{eq:compatibility-neighborhood}. \qedhere

\end{proof}

Proposition~\ref{prop:compatibility-neighborhood} offers a two-sided geometric description of the asymptotic compatibility region $\Theta_\tau$\;\!\! around $\Theta^\star$\:\!\!.
The function $\psi_{\;\!\! +}$\;\!\! identifies neighborhoods that remain asymptotically compatible, while $\psi_{\;\!\! -}$\;\!\! bounds the local extent of $\Theta_\tau$\;\!\! away from the set of model\:\!-wide discrepancy minimizers.

\begin{remark}

Whenever $\psi_{\;\!\! -}$\;\!\! and $\psi_{\;\!\! +}$\;\!\! are continuous and strictly increasing, the inner and outer geometry of the compatibility region is governed, respectively, by the inverse scales
\[
\psi_{\;\!\! +}^{\;\!\! -1} \;\!\!
\big{(} \tau - \epsilon_{\Theta}^{\;\!\! \star} \big{)} ,
\
\psi_{\;\!\! -}^{\;\!\! -1} \;\!\!
\big{(} \tau - \epsilon_{\Theta}^{\;\!\! \star} \big{)} ,
\]
provided that both are well defined and belong to $\mo{]} 0 \:\!,\;\!\! r^{\;\!\! \star} \;\!\!\mc{[}$.
Thus, flatter discrepancy profiles may produce compatibility regions extending farther from $\Theta^\star$\:\!\!, while faster growth leads to tighter localization.

\end{remark}

Under polynomial discrepancy growth, the same argument makes the compatibility scales explicit.

\begin{corollary}
\label{cor:polynomial-compatibility-growth} 

Let $r^{\;\!\! \star} \:\!\!, q , c_{-} , c_{+} \;\!\! \in \mo{]}\:\! 0 \:\!,\:\!\! \infty \mc{[}$ with $c_{-} \;\!\! \le c_{+}$\;\!\! be such that, for $\bm{\pi}$-\:\!almost every $\bm{\theta} \in \Theta_{r^{\;\!\! \star}}^\star$,
\begin{equation*}
c_{-} \:\!
\bd_{\:\! \Theta^\star \:\!\!}^{\;\! q}\;\!\!(\bm{\theta})
\le
\cT_{\bm{\theta}}^{\;\! \star} \:\!\!
-
\epsilon_{\Theta}^{\;\!\! \star} \;\!\!
\le
c_{+} \:\!
\bd_{\:\! \Theta^\star \:\!\!}^{\;\! q}\;\!\!(\bm{\theta}) \:\! .
\end{equation*}
Assuming that
\[
\bigg{(}
\frac{\tau - \epsilon_{\Theta}^{\;\!\! \star} \;\!\!}{c_{-}}
\bigg{)}^{\:\!\! 1/q} \:\!\!\!
<
r^{\;\!\! \star} \:\!\! ,
\]
the following two inclusions hold up to $\bm{\pi}$-\:\!null sets:
\begin{equation*}
\Theta_{\! \big{(}\:\!\! \frac{\tau \:\! - \:\! \epsilon_{\Theta}^{\;\!\! \star} \;\!\!}{c_{+}} \:\!\!\big{)}^{\:\!\! 1/q} \:\!\!\!}^\star
\subseteq
\Theta_\tau ,
\quad
\Theta_\tau \;\!\!
\cap
\Theta_{r^{\;\!\! \star}}^{\star} \:\!\!
\subseteq
\Theta_{\! \big{(}\:\!\! \frac{\tau \:\! - \:\! \epsilon_{\Theta}^{\;\!\! \star} \;\!\!}{c_{-}} \:\!\!\big{)}^{\:\!\! 1/q} \:\!\!\!}^\star \, .
\end{equation*}

\end{corollary}

The inner localization in Proposition~\ref{prop:compatibility-neighborhood} can also be coupled with Proposition~\ref{prop:partial-stabilization} and Corollary~\ref{cor:partial-stabilization}:
if $\bm{\pi} \bracks{ \Theta_{r_{-} \;\!\!}^{\star} } > 0$ and condition~\eqref{eq:partial-stabilization} holds with a uniform $\delta_{\;\!\! \tau}$ for $\bm{\pi}$-\:\!almost every $\bm{\theta} \in \Theta_{r_{-}}^{\star}$, we may take
\[
B_{\;\!\! \tau} \;\!\! = \Theta_{r_{-}}^{\star} \;\!\! .
\]
Consequently, inequality~\eqref{eq:partial-stabilization-integrated} specializes to
\begin{equation*}
\liminf_{n \:\! \to \:\! \infty} \;\!
\bE_{\:\! \bP}
\bracks*{
\int\nolimits_{\Theta_{\:\! r_{-}}^{\:\! \star}} \:\!\!\! w_{\:\!\! \tau}^{\:\! n} \big{(} \bm{\theta} ;\;\!\! \bY^{1:\:\!n} \big{)} \:\! \bm{\pi}(\d\bm{\theta})
\;\!\!} \;\!\!
\ge
\rounds{1 - \delta_{\;\!\! \tau}} \:\! \bm{\pi} \bracks{ \Theta_{r_{-} \;\!\!}^{\star} } \:\! .
\end{equation*}

Observe further that, in the polynomial case of Corollary~\ref{cor:polynomial-compatibility-growth}, one may choose any $r_{-} \;\!\! > 0$ satisfying
\[
r_{-} \;\!\!
<
\bigg{(}
\frac{\tau - \epsilon_{\Theta}^{\;\!\! \star} \;\!\!}{c_{+}}
\bigg{)}^{\:\!\! 1/q} \:\!\!\! .
\]

Taken together, the results of this section distinguish two complementary large-sample regimes.
Under full stabilization of the observed and simulated probabilistic representations, the posterior converges to the prior filtered by $\Theta_\tau$.
Under partial stabilization, exact filtering may fail, while positive prior-averaged acceptance can persist quantitatively within neighborhoods controlled by the local growth of the asymptotic discrepancy.
This analysis provides the foundation for the next section.

\section{Variational and transport structures in ABC}
\label{sec:variational-transport} 

We first characterize compatibility-weight posteriors through the classical Gibbs identity and prove stability of minimum values and approximate minimizers.
We then identify the constrained entropy problem solved by the fixed-tolerance large-sample limit.
Finally, we introduce additional penalty-based functionals of the limiting discrepancy and study their reduction to distributions over model laws and their Wasserstein geometry.
These last functionals are separate generalized Bayesian constructions; their minimizers are not asserted to be limits of the original fixed-tolerance ABC procedure.

\subsection{A compatibility-induced variational principle}
\label{subsec:compatibility-variational} 

Throughout the subsection, let $\epsilon > 0$ be given, and let $w_{\:\!\! \epsilon} \colon \Theta \to \R_{\:\! +}$\;\!\! be a compatibility weight as introduced in Subsection~\ref{subsec:convergence-principle}.
From a variational perspective, rather than viewing $w_{\:\!\! \epsilon}$\;\!\! solely as quantifying the compatibility of each parameter value $\bm{\theta} \in \Theta$ with the observed sample $\by^{1:\:\!n} \:\!\! \in \cY^{\:\! n}$\:\!\!, we regard it as defining a data-dependent loss structure on the parameter space.

Henceforth, we restrict our attention to compatibility weights satisfying, $\bm{\pi}$-\:\!a.s.,
\begin{equation}
\label{eq:w_epsilon-bound} 
w_{\:\!\! \epsilon} \le 1 .
\end{equation}
Up to a modification on a $\bm{\pi}$-\:\!null set, we may assume that~\eqref{eq:w_epsilon-bound} holds everywhere on $\Theta$, and introduce the corresponding extended-valued compatibility loss $\ell_{\:\!\! \epsilon}$\;\!\! by setting, for any $\bm{\theta} \in \Theta$, 
\begin{equation}
\label{eq:l_epsilon} 
\ell_{\:\!\! \epsilon}(\bm{\theta}) \;\!\!
\coloneqq
\begin{cases}
\displaystyle
\:\! - \log w_{\:\!\! \epsilon}(\bm{\theta}) ,
& \text{if $w_{\:\!\! \epsilon}(\bm{\theta}) > 0$} \:\! ,
\\[0.25ex]
\:\! \infty ,
& \text{otherwise} \:\! .
\end{cases}
\end{equation}

This logarithmic transformation converts the multiplicative compatibility scale into an additive one, playing a role analogous to that of the negative log-likelihood in classical statistical inference.
Larger values of $w_{\:\!\! \epsilon}$\;\!\! correspond to smaller values of $\ell_{\:\!\! \epsilon}$, so that minimizing the compatibility loss is equivalent to maximizing compatibility with the observed sample.
Under the present assumptions, the induced compatibility loss is nonnegative, consistent with the classical ABC setting (see also Subsection~\ref{subsec:ABC-classical}).
This additive formulation is particularly well suited for variational analysis, as it naturally gives rise to entropy-regularized problems in approximate Bayesian inference.

A further connection arises with loss-based forms of generalized Bayesian updating; see, e.g.,~\cite{bissiri-holmes-walker_2016}.
The compatibility loss~\eqref{eq:l_epsilon} plays the role of an empirical loss, while the compatibility weights act as likelihood surrogates.
The resulting inferential structure is induced by the discrepancy geometry relating observed and simulated data, without requiring explicit evaluation of the observation likelihood.

At the level of posterior candidates, this pointwise loss induces the extended-valued functional
\begin{equation}
\label{eq:R-epsilon} 
\bm{\rho}
\longmapsto
\;\!\!
\int\nolimits_{\Theta} \ell_{\:\!\! \epsilon}(\bm{\theta}) \;\! \bm{\rho}\:\!(\d\bm{\theta}) \:\! ,
\quad
\bm{\rho} \in \ccP(\Theta) \:\! ,
\end{equation}
which represents the compatibility loss averaged with respect to the candidate posterior distribution $\bm{\rho}$.
Thus, posterior candidates are no longer evaluated only through the compatibility of individual parameter values, but through the aggregate loss associated with their distribution over the parameter space.

This formulation extends to a general variational framework suited for the subsequent developments.
It is based on the interplay between a data-fitting functional
\[
\cR_{\:\!\! \epsilon} \colon \ccP(\Theta) \longrightarrow \mo{[} 0 \:\!,\:\!\! \infty \mc{]}
\]
and a regularization functional
\[
\cI(\;\! \bm{\cdot} \;\!\! \mid \:\!\! \bm{\pi}) \colon \ccP(\Theta) \longrightarrow \mo{[} 0 \:\!,\:\!\! \infty \mc{]} \:\! ,
\]
with the former representing the compatibility component associated with the observed sample $\by^{1:\:\!n}$\:\!\!, and the latter acting as a prior-based informational regularization term.
This yields the variational problem
\begin{equation}
\label{eq:inf_R-I} 
\inf_{\bm{\rho} \:\! \in \:\! \ccP(\Theta)} \;\!\!
\braces[\big]{
\cR_{\:\!\! \epsilon} (\bm{\rho} \:\!)
+
\cI(\bm{\rho} \:\!\! \mid \:\!\! \bm{\pi})
} .
\end{equation}


In the canonical compatibility-loss setting introduced above, the functional $\cR_{\:\!\!\epsilon}$ reduces precisely to~\eqref{eq:R-epsilon}.
In this case, it admits a natural interpretation as a compatibility-induced risk functional, thereby motivating its notation.
The regularization term then acts to prevent the minimization from collapsing into a purely compatibility-driven concentration on the most compatible regions of the parameter space, while preserving the role of the prior as an informational reference measure.

Among the possible choices for the functional $\cI(\;\! \bm{\cdot} \;\!\! \mid \:\!\! \bm{\pi})$, the Kullback--Leibler divergence, also known as relative entropy, emerges as a particularly natural candidate.
It is defined on $\ccP(\Theta)$ by
\begin{equation}
\label{eq:KL} 
\mathrm{KL}\:\!(\bm{\rho} \:\!\! \mid \:\!\! \bm{\pi})
\coloneqq
\begin{cases}
\displaystyle
\:\! \int\nolimits_{\Theta} \log \:\!
\bracks*{\;\!\!
\frac{\d\bm{\rho}}{\d\bm{\pi}}(\bm{\theta})
\;\!\!} \:\! \bm{\rho}\:\!(\d\bm{\theta}) ,
& \text{if $\bm{\rho} \ll \;\!\! \bm{\pi}$} ,
\\[1.5ex]
\:\! \infty ,
& \text{otherwise} \:\! .
\end{cases}
\end{equation}
The relative entropy quantifies the informational divergence of the posterior candidate $\bm{\rho}$ from the prior distribution $\bm{\pi}$.
Within this framework, it acts as an informational regularization term, discouraging excessive departures from the prior while balancing the influence of the compatibility component.

Within this compatibility-loss setting, and taking the Kullback--Leibler divergence as the regularization functional, the variational problem~\eqref{eq:inf_R-I} specializes to
\begin{equation}
\label{eq:inf_l-KL} 
\inf_{\bm{\rho} \:\! \in \:\! \ccP(\Theta)}
\braces*{ 
\:\! \int\nolimits_{\Theta} \ell_{\:\!\! \epsilon}(\bm{\theta}) \;\! \bm{\rho}\:\!(\d\bm{\theta}) \;\!\!
+
\mathrm{KL}\:\!(\bm{\rho} \:\!\! \mid \:\!\! \bm{\pi})
} .
\end{equation}


This is an entropy-regularized risk minimization problem over $\ccP(\Theta)$, balancing compatibility with the observed data against informational divergence from the prior.


The compatibility-weight posterior distribution $\bm{\pi}_{\:\!\! \epsilon} \;\!\! \in \ccP(\Theta)$ introduced in~\eqref{eq:pi_epsilon} admits a variational characterization in terms of the objective functional $\cF_{\:\!\! \epsilon}$\;\!\! on $\ccP(\Theta)$ given by 
\begin{equation}
\label{eq:F_epsilon} 
\cF_{\:\!\! \epsilon} (\bm{\rho} \:\!)
\coloneqq \;\!\!
\int\nolimits_{\Theta} \ell_{\:\!\! \epsilon}(\bm{\theta}) \;\! \bm{\rho}\:\!(\d\bm{\theta}) \;\!\!
+
\mathrm{KL}\:\!(\bm{\rho} \:\!\! \mid \:\!\! \bm{\pi}) \:\! .
\end{equation}
This representation is a compatibility-based instance of the classical Gibbs--Donsker--Varadhan principle; see, e.g.,~\cite{donsker-varadhan_1975,dembo-zeitouni_1998}.

\begin{proposition}
\label{prop:inf_l-KL} 

Under the standing assumptions of this subsection, the probability measure $\bm{\pi}_{\:\!\! \epsilon}$, defined as in~\eqref{eq:pi_epsilon}, is the unique minimizer of $\cF_{\:\!\! \epsilon}$\;\!\! on $\ccP(\Theta)$, with minimum value 
\begin{equation*}
\min_{\bm{\rho} \:\! \in \:\! \ccP(\Theta)} \:\!\!
\cF_{\:\!\! \epsilon} (\bm{\rho} \:\!)
\equiv
\cF_{\:\!\! \epsilon}(\bm{\pi}_{\:\!\! \epsilon})
=
- \log \int\nolimits_{\Theta} w_{\:\!\! \epsilon}(\bm{\theta}) \;\! \bm{\pi}(\d\bm{\theta}) \:\! .
\end{equation*}

\end{proposition}

\begin{proof}
Let $C_{\;\!\! \epsilon}$ denote the normalization constant associated with $\bm{\pi}_{\:\!\! \epsilon}$, namely,
\begin{equation}
\label{eq:C_epsilon}
C_{\;\!\! \epsilon}
\coloneqq \;\!\!
\int\nolimits_{\Theta} w_{\:\!\! \epsilon}(\bm{\theta}) \;\! \bm{\pi}(\d\bm{\theta}) \:\! .
\end{equation}
By the Gibbs variational identity, understood in the extended sense,
\begin{equation}
\label{eq:KL_epsilon}
\cF_{\:\!\! \epsilon}(\bm{\rho})
=
-\log C_{\;\!\! \epsilon}
+
\mathrm{KL}\:\!(\bm{\rho} \:\!\! \mid \:\!\! \bm{\pi}_{\:\!\! \epsilon}) \:\! ,
\quad
\bm{\rho} \in \ccP(\Theta) \:\! .
\end{equation}
The conclusion follows from the nonnegativity of relative entropy and its equality case. \qedhere
\end{proof}

Proposition~\ref{prop:inf_l-KL} identifies the compatibility-weight posterior as the unique entropy-regularized minimizer, with the normalization constant~\eqref{eq:C_epsilon} acting as its partition function.





We next examine the stability of this variational characterization under the vanishing-threshold behavior established in Proposition~\ref{prop:threshold-limit}.

\begin{theorem}
\label{thm:inf_l-KL-asy} 

Under the standing assumptions of this subsection and of Proposition~\ref{prop:threshold-limit}, define 
\begin{equation}
\label{eq:F_epsilon-nor} 
\widetilde{\cF}_{\:\!\! \epsilon} (\bm{\rho} \:\!)
\coloneqq
\cF_{\:\!\! \epsilon} (\bm{\rho} \:\!)
+
\log \kappa_{\;\!\! \epsilon} \:\! ,
\quad
\bm{\rho} \in \ccP(\Theta) \:\! ,
\end{equation}
with $\cF_{\:\!\! \epsilon}$\;\!\! as in~\eqref{eq:F_epsilon}, and
\begin{equation}
\label{eq:F-bar} 
\widebar{\cF} (\bm{\rho} \:\!)
\coloneqq
- \log \int\nolimits_{\Theta} \widebar{w}(\bm{\theta}) \;\! \bm{\pi}(\d\bm{\theta}) \;\!\!
+
\mathrm{KL}\:\!(\bm{\rho} \:\!\! \mid \:\!\! \widebar{\bm{\pi}}) \:\! ,
\quad
\bm{\rho} \in \ccP(\Theta) \:\! ,
\end{equation}
where $\widebar{\bm{\pi}}$ is given by~\eqref{eq:bar-pi}.
Then $\bm{\pi}_{\:\!\! \epsilon}$\;\!\! is the unique minimizer of $\widetilde{\cF}_{\:\!\! \epsilon}$\;\!\! on $\ccP(\Theta)$, with minimum value
\begin{equation}
\label{eq:min_F_epsilon-nor} 
\min_{\bm{\rho} \:\! \in \:\! \ccP(\Theta)} \:\!\!
\widetilde{\cF}_{\:\!\! \epsilon} (\bm{\rho} \:\!)
\equiv
\widetilde{\cF}_{\:\!\! \epsilon}(\bm{\pi}_{\:\!\! \epsilon})
=
- \log \int\nolimits_{\Theta} \frac{1}{\kappa_{\;\!\! \epsilon}} \:\! w_{\:\!\! \epsilon}(\bm{\theta}) \;\! \bm{\pi}(\d\bm{\theta}) \:\! .
\end{equation}
Likewise, the probability measure $\widebar{\bm{\pi}}$ is the unique minimizer of $\widebar{\cF}$ on $\ccP(\Theta)$, and
\begin{equation}
\label{eq:min_F-bar} 
\min_{\bm{\rho} \:\! \in \:\! \ccP(\Theta)} \:\!\!
\widebar{\cF} (\bm{\rho} \:\!)
\equiv
\widebar{\cF}(\widebar{\bm{\pi}})
=
- \log \int\nolimits_{\Theta} \widebar{w}(\bm{\theta}) \;\! \bm{\pi}(\d\bm{\theta}) \:\! .
\end{equation}
Furthermore,
\begin{equation}
\label{eq:lim_F_epsilon-nor_F-bar} 
\lim_{\epsilon \:\! \downarrow \:\! 0} \;\!
\widetilde{\cF}_{\:\!\! \epsilon}(\bm{\pi}_{\:\!\! \epsilon}) \;\!\!
=
\widebar{\cF}(\widebar{\bm{\pi}}) \:\! .
\end{equation}
Finally, if $\Set{\! \bm{\rho}_{\;\!\! \epsilon} \!}_{\;\!\! \epsilon \:\! > \:\! 0}$\;\!\! is any family in $\ccP(\Theta)$ of asymptotic minimizers, that is,
\begin{equation}
\label{eq:lim_F_epsilon-nor} 
\lim_{\epsilon \:\! \downarrow \:\! 0} \,
\bracks*{
\widetilde{\cF}_{\:\!\! \epsilon}(\bm{\rho}_{\;\!\! \epsilon})
- \widetilde{\cF}_{\:\!\! \epsilon}(\bm{\pi}_{\:\!\! \epsilon})
} \;\!\!
=
0 \:\! ,
\end{equation}
then this family converges in total variation to $\widebar{\bm{\pi}}$.

\end{theorem}

\begin{proof}

The characterization of $\bm{\pi}_{\:\!\! \epsilon}$\;\!\! as the unique minimizer of $\widetilde{\cF}_{\:\!\! \epsilon}$, together with identity~\eqref{eq:min_F_epsilon-nor}, follows directly from Proposition~\ref{prop:inf_l-KL}.
The corresponding statement for $\widebar{\cF}$ is immediate from definition~\eqref{eq:F-bar}.
Identity~\eqref{eq:lim_F_epsilon-nor_F-bar} is then obtained by combining~\eqref{eq:min_F_epsilon-nor} and~\eqref{eq:min_F-bar} with the $\Lped{\;\!\! 1 \:\!\!}{\bm{\pi}}{}$-\:\!convergence~\eqref{eq:bar-w} as $\epsilon \downarrow 0$.


It remains to demonstrate the stability of asymptotic minimizers $\Set{\! \bm{\rho}_{\;\!\! \epsilon} \!}_{\;\!\! \epsilon \:\! > \:\! 0}$.
By~\eqref{eq:F_epsilon-nor} and~\eqref{eq:KL_epsilon}, one has
\[
\widetilde{\cF}_{\:\!\! \epsilon}(\bm{\rho}_{\;\!\! \epsilon})
- \widetilde{\cF}_{\:\!\! \epsilon}(\bm{\pi}_{\:\!\! \epsilon})
=
\mathrm{KL}\:\!(\bm{\rho}_{\;\!\! \epsilon} \:\!\! \mid \:\!\! \bm{\pi}_{\:\!\! \epsilon}) \:\! ,
\]
and, in particular, condition~\eqref{eq:lim_F_epsilon-nor} may be rewritten as
\begin{equation}
\label{eq:lim_KL_epsilon} 
\lim_{\epsilon \:\! \downarrow \:\! 0} \,
\mathrm{KL}\:\!(\bm{\rho}_{\;\!\! \epsilon} \:\!\! \mid \:\!\! \bm{\pi}_{\:\!\! \epsilon})
=
0 \:\! .
\end{equation}
On the other hand, Pinsker's inequality gives
\begin{equation}
\label{eq:pinsker-inequality} 
\sup_{B \:\! \in \:\! \cB(\Theta)} \;\!
\abs[\big]{\:\! \bm{\rho}_{\;\!\! \epsilon} \bracks{ B \:\!} - \bm{\pi}_{\:\!\! \epsilon} \bracks{ B \:\!} \:\!}
\le
\sqrt{
\:\! \frac{1}{2} \:\!
\mathrm{KL}\:\!(\bm{\rho}_{\;\!\! \epsilon} \:\!\! \mid \:\!\! \bm{\pi}_{\:\!\! \epsilon})
} \;\! ,
\end{equation}
which, combined with~\eqref{eq:lim_KL_epsilon}, shows that the left-hand side in~\eqref{eq:pinsker-inequality} also vanishes as $\epsilon \downarrow 0$.

In conclusion, since Proposition~\ref{prop:threshold-limit} ensures that $\Set{\! \bm{\pi}_{\:\!\! \epsilon} \!}_{\;\!\! \epsilon}$ converges to $\widebar{\bm{\pi}}$ in total variation, an application of the triangle inequality implies the same convergence for $\Set{\! \bm{\rho}_{\;\!\! \epsilon} \!}_{\;\!\! \epsilon \:\! > \:\! 0}$. \qedhere

\end{proof}

\begin{remark}
The theorem does not establish pointwise convergence or $\Gamma$\;\!\!-\:\!convergence of the functionals.
A $\Gamma$\;\!\!-\:\!convergence statement would require separate liminf and recovery arguments in a specified topology; see~\cite{dal.maso_1993}.
Extending the stability conclusion to other regularizers likewise requires additional assumptions controlling the distance from approximate minimizers to the relevant posterior.
\end{remark}

\paragraph{The fixed-tolerance limit as an entropy projection.}
\label{par:hard-filter}
Under the assumptions of Theorem~\ref{thm:compatibility-filtering}, let $a=\bm{\pi}[\Theta_\tau]>0$ and $\bm{\pi}_\tau^\infty=\bm{\pi}[\,\bm{\cdot}\mid\Theta_\tau]$.
The exact variational characterization of this limiting posterior is
\begin{equation}
\label{eq:hard-filter-variational}
\min_{\substack{\bm{\rho}\in\ccP(\Theta)\\\bm{\rho}[\Theta_\tau]=1}}
\mathrm{KL}(\bm{\rho}\mid\bm{\pi})=-\log a,
\qquad\text{with unique minimizer }\bm{\pi}_\tau^\infty.
\end{equation}
Indeed, for every feasible $\bm{\rho}$,
\[
\mathrm{KL}(\bm{\rho}\mid\bm{\pi})
=-\log a+\mathrm{KL}(\bm{\rho}\mid\bm{\pi}_\tau^\infty).
\]
For $\bm{\rho}\ll\bm{\pi}$ this follows from
$\d\bm{\pi}_\tau^\infty/\d\bm{\pi}=1/a$ $\bm{\rho}$-almost everywhere.
If $\bm{\rho}\not\ll\bm{\pi}$, both relative entropies are infinite.
Nonnegativity and the equality case of relative entropy prove~\eqref{eq:hard-filter-variational}.
Equivalently, the limit minimizes expected loss plus relative entropy for the loss equal to zero on $\Theta_\tau$ and infinity on its complement.
This characterizes the limit itself; no convergence of the finite-sample loss functionals is claimed.

\paragraph{Additional penalties on the limiting discrepancy.}
\label{par:additional-penalty}
The profile $\bm{\theta}\mapsto\cT_{\bm{\theta}}^{\;\! \star}$\;\!\! also permits a separate construction in which departures from the prior-essential optimum are assigned a finite penalty.
Unlike the model-wide threshold~\eqref{eq:epsilon_Theta^star}, the prior-essential threshold~\eqref{eq:epsilon_pi^star} is invariant under modifications of the discrepancy profile on $\bm{\pi}$-null sets.
We adopt the latter as the reference level for the entropy-regularized update.
The penalty function, including its scale, is prescribed independently of the fixed-tolerance ABC procedure.

For this additional construction, define the asymptotic compatibility functional
\begin{equation}
\label{eq:R_phi^star} 
\cR_{\;\!\! \phi}^{\star \;\!\!}(\bm{\rho} \:\!)
\coloneqq \;\!\!
\int\nolimits_{\Theta}
\phi \:\! \big{(}
\bracks*{\;\!\! \cT_{\bm{\theta}}^{\;\! \star}\:\!\! - \epsilon_{\;\!\! \bm{\pi}}^{\;\!\! \star} }_{\;\!\! +}
\big{)} \;\!
\bm{\rho}\:\!(\d\bm{\theta}) \:\! ,
\quad
\bm{\rho} \in \ccP(\Theta) \:\! , \
\bm{\rho} \ll \;\!\! \bm{\pi}
\end{equation}
(understood in the extended sense), with $\bracks{\, \bm{\cdot} \,}_{+}$\;\!\! denoting the positive-part operator on~$\R$.
The prescribed penalty function
\[
\phi \colon \R_{\:\! +} \;\!\! \longrightarrow \R_{\:\! +}
\]
is strictly increasing, continuous at zero, and satisfies $\phi(0) = 0$, while $\epsilon_{\;\!\! \bm{\pi}}^{\;\!\! \star}$\;\!\! is the prior-essential compatibility threshold~\eqref{eq:epsilon_pi^star}.

The profile $\bm{\theta} \mapsto \cT_{\bm{\theta}}^{\;\! \star}$\;\!\! determines the asymptotic compatibility structure on $\Theta$, while $\phi$ quantifies deviations from the best compatibility level accessible under the prior.
The resulting functional~$\cR_{\;\!\! \phi}^\star$ evaluates the average penalized asymptotic incompatibility of posterior candidates, separating the asymptotic discrepancy profile from the penalization mechanism.

\begin{remark}
\label{rem:R_phi^star}

The functional~$\cR_{\;\!\! \phi}^\star$ depends only on the $\bm{\pi}$-\:\!almost everywhere equivalence class of the asymptotic discrepancy profile.
Moreover,
\[
\inf_{\substack{
\bm{\rho} \:\! \in \:\! \ccP(\Theta)
\\
\bm{\rho} \:\! \ll \bm{\pi}
}} \:\!\!
\cR_{\;\!\! \phi}^{\star \;\!\!}(\bm{\rho} \:\!)
=
0 \:\! .
\]
Indeed, the defining property of $\epsilon_{\;\!\! \bm{\pi}}^{\;\!\! \star}$\;\!\! and the continuity of $\phi$ at zero allow the penalty to be made arbitrarily small.

The infimum is attained if and only if
\[
\bm{\pi}\bracks*{\;\!\!
\Set{\!
\bm{\theta} \in \Theta
|
\cT_{\bm{\theta}}^{\;\! \star} \:\!\!
=
\epsilon_{\;\!\! \bm{\pi}}^{\;\!\! \star}
\!}
\;\!\!}
>
0 \:\! .
\]
In that case, the minimizers are precisely the probability measures $\bm{\rho} \ll \;\!\! \bm{\pi}$ concentrated on this set.

The entropy-regularized minimizer associated with~\eqref{eq:R_phi^star} is a separate update, rather than the conditional prior characterized by~\eqref{eq:hard-filter-variational}.

\end{remark}

Among the possible geometric extensions, an optimal-transport approach constitutes a natural specialization of the framework introduced so far. 

As anticipated in Section~\ref{sec:large-sample}, transport-based choices of the discrepancy map~\eqref{eq:T} induce compatibility mechanisms in Approximate Bayesian Computation.
The same geometric principles also allow posterior candidates to be represented through probability distributions on the space of model laws.

These considerations motivate the formulation introduced in the next subsection.




\subsection{Optimal transport structures}
\label{subsec:optimal-transport} 

In this subsection, compatibility is characterized via transport discrepancies on the data space, while the associated variational structure is represented through the transport geometry of the induced model laws.
The analysis is carried out within the theory of optimal transport; see, e.g., the foundational contribution~\cite{kantorovich_1942}, the classical monographs~\cite{villani_2003,villani_2009}, and the more recent treatments~\cite{ambrosio-gigli-savare_2008,santambrogio_2015,peyre-cuturi_2019}. 

In this subsection, we additionally assume that the metric $\bd_{\;\!\! \cY}$ is complete.
Together with the separability imposed in Subsection~\ref{subsec:ABC-framework}, this makes $(\cY,\bd_{\;\!\! \cY})$ a Polish metric space.
This is the case, for instance, when $\cY$ is a closed subset of $\R^{m_{\;\!\! \cY}}$\:\!\! endowed with the induced Euclidean metric.
Besides encompassing most statistical models encountered in practice, this provides the regularity needed for the transport-theoretic constructions.

We restrict attention to the classical Wasserstein setting generated by the ground metric $\bd_{\;\!\! \cY}$, although several of the constructions below admit extensions to more general transport costs.

To this end, we fix once and for all an exponent
\[
p \in \mo{[}\;\!\! 1 ,\:\!\! \infty \mc{[} .
\]
Let $\ccP_{\;\!\! p}(\cY)$ denote the subclass of $\ccP(\cY)$ consisting of all probability measures $\bm{\eta}$ on $\cY$ admitting a finite $p$-th moment; i.e., for some (equivalently, every) $\by_{\:\! 0} \;\!\! \in \cY$\;\!\!, 
\[
\int\nolimits_{\;\!\! \cY} \bd_{\:\!\! \cY}^{\;\! p}(\by \:\!,\;\!\! \by_{\:\! 0}) \;\! \bm{\eta}(\d\by)
<
\infty \:\! . 
\] 


We equip $\ccP_{\;\!\! p}(\cY)$ with the $p$\:\!-Wasserstein distance associated with $\bd_{\;\!\! \cY}$, defined by
\[
W_{\:\!\! p}(\bm{\mu} \:\!,\;\!\! \bm{\nu})
\coloneqq
\bracks*{
\inf_{\bm{\gamma} \:\! \in \:\! \Pi(\bm{\mu},\bm{\nu})}
\int\nolimits_{\;\!\! \cY \times \cY} \;\!\!
\bd_{\:\!\! \cY}^{\;\! p}(\by \:\!,\;\!\! \bz) \;\! \bm{\gamma}\:\!(\d\by \:\!,\;\!\! \d\bz)
}^{1/p} \! ,
\quad
\bm{\mu} \:\!,\;\!\! \bm{\nu} \in \ccP_{\;\!\! p}(\cY) \:\! ,
\]
where $\Pi(\bm{\mu} \:\!,\;\!\! \bm{\nu})$ denotes the set of all couplings of $\bm{\mu}$ and $\bm{\nu}$, that is, all probability measures on $\cY \;\!\! \times \cY$ whose first and second marginals are $\bm{\mu}$ and $\bm{\nu}$, respectively.

The above definition coincides with the Kantorovich formulation of the optimal transport problem corresponding to the cost function $\bd_{\:\!\! \cY}^{\;\! p}$.
Under the present assumptions, this problem admits a minimizer, namely an optimal transport plan. 
Moreover, the Wasserstein space $\rounds{ \ccP_{\;\!\! p}(\cY) \:\!,\:\!\! W_{\:\!\! p} }$ is again Polish.

Accordingly, the framework introduced in Subsection~\ref{subsec:probabilistic-representations} is henceforth understood upon replacing $\ccP(\cY)$ with $\ccP_{\;\!\! p}(\cY)$.
Specifically, for every $\bm{\mu} , \bm{\nu} \in \ccP_{\;\!\! p}(\cY)$, we set
\[
\cT(\bm{\mu} \:\!,\;\!\! \bm{\nu})
\coloneqq
W_{\:\!\! p}(\bm{\mu} \:\!,\;\!\! \bm{\nu}) \:\! .
\]

The statistical model associates with every parameter value $\bm{\theta} \in \Theta$ the corresponding law $\bm{\nu}_{\bm{\theta}} \in \ccP_{\;\!\! p}(\cY)$, introduced in~\eqref{eq:eta_theta}.
The family $\Set{\! \bm{\nu}_{\bm{\theta}} \!}_{\;\!\! \bm{\theta} \:\! \in \:\! \Theta}$\;\!\! provides a representation of the model class in $\ccP_{\;\!\! p}(\cY)$, while the limiting measure $\bm{\mu}^{\;\!\! \star} \:\!\! \in \ccP_{\;\!\! p}(\cY)$, specified in~\eqref{eq:mu^star}, encodes the observational information.

In particular, the asymptotic compatibility profile~\eqref{eq:T_theta^star} is given, for every $\bm{\theta} \in \Theta$, by
\begin{equation}
\label{eq:T_theta^star-p} 
\cT_{\bm{\theta}}^{\;\! \star} \:\!\!
=
W_{\:\!\! p} \big{(}\:\! \bm{\mu}^{\;\!\! \star} \:\!\!,\;\!\! \bm{\nu}_{\bm{\theta}} \big{)} .
\end{equation}
Under the attainment assumption~\eqref{eq:epsilon_Theta^star-min}, the model\:\!-wide compatibility threshold~\eqref{eq:epsilon_Theta^star} becomes
\[
\epsilon_{\Theta}^{\;\!\! \star} \;\!\!
=
\min_{\bm{\theta} \:\! \in \:\! \Theta} \:\!
W_{\:\!\! p} \big{(}\:\! \bm{\mu}^{\;\!\! \star} \:\!\!,\;\!\! \bm{\nu}_{\bm{\theta}} \big{)} .
\]

The asymptotic compatibility functional~\eqref{eq:R_phi^star} then takes the explicit Wasserstein form
\begin{equation}
\label{eq:R_phi,p^star} 
\cR_{\;\!\! \phi \;\!\!, p \;\!\!}^{\star \;\!\!}(\bm{\rho} \:\!)
\coloneqq \;\!\!
\int\nolimits_{\Theta}
\phi \:\! \big{(}
\bracks*{\;\!\! W_{\:\!\! p} \big{(}\:\! \bm{\mu}^{\;\!\! \star} \:\!\!,\;\!\! \bm{\nu}_{\bm{\theta}} \big{)} \;\!\! - \epsilon_{\;\!\! \bm{\pi}}^{\;\!\! \star} }_{\;\!\! +}
\big{)} \:\!
\bm{\rho}\:\!(\d\bm{\theta}) \:\! ,
\quad
\bm{\rho} \in \ccP(\Theta) \:\! , \
\bm{\rho} \ll \;\!\! \bm{\pi} ,
\end{equation}
where $\phi \colon \R_{\:\! +} \;\!\! \to \R_{\:\! +}$\;\!\!
is the strictly increasing penalty function continuous at zero with $\phi(0) = 0$.

Remark~\ref{rem:R_phi^star} implies that, for every $\bm{\rho} \in \ccP(\Theta)$ with $\bm{\rho} \ll \;\!\! \bm{\pi}$,
\[
\cR_{\;\!\! \phi \;\!\!, p \;\!\!}^{\star \;\!\!}(\bm{\rho} \:\!)
=
0
\quad
\Longleftrightarrow
\quad
\bm{\rho}
\bracks*{\;\!\!
\Set{\!
\bm{\theta} \in \Theta
|
W_{\:\!\! p} \big{(}\:\! \bm{\mu}^{\;\!\! \star} \:\!\!,\;\!\! \bm{\nu}_{\bm{\theta}} \big{)} \;\!\!
=
\epsilon_{\;\!\! \bm{\pi}}^{\;\!\! \star}
\!}
\;\!\!}
=
1 .
\]


Altogether, these observations show that asymptotic compatibility depends on the parameters only through the corresponding model laws.
It is therefore natural to represent asymptotic compatibility directly on the Wasserstein space $\ccP_{\;\!\! p}(\cY)$ of model laws.

To formalize this idea, we assume that the parameter-to-distribution map $\Phi \colon \Theta \to \ccP_{\;\!\! p}(\cY)$ defined by
\begin{equation}
\label{eq:Phi} 
\Phi(\bm{\theta})
=
\bm{\nu}_{\bm{\theta}} \:\! ,
\quad
\bm{\theta} \in \Theta \:\! ,
\end{equation}
is Borel measurable.
In this way, every $\bm{\rho} \in \ccP(\Theta)$, $\bm{\rho} \ll \;\!\! \bm{\pi}$, induces the push-forward (or image) measure
\begin{equation}
\label{eq:Lambda_rho} 
\bm{\Lambda}_{\bm{\rho}} \;\!\!
\coloneqq
\Phi_{\:\! \#} \:\! \bm{\rho}
\in
\ccP \big{(} \ccP_{\;\!\! p}(\cY) \big{)} ,
\end{equation}
which represents the distribution of the random probability measure $\Phi(\bm{\theta}) = \bm{\nu}_{\bm{\theta}}$ when the parameter $\bm{\theta}$ is distributed according to $\bm{\rho}$.
Indeed, by definition, 
\[
\Phi_{\:\! \#} \:\! \bm{\rho} \:\! \bracks*{ \cE }
=
\bm{\rho} \:\! \bracks*{\;\!\! \Phi^{-1}(\cE) \;\!\!} ,
\quad
\cE \in \cB \big{(} \ccP_{\;\!\! p}(\cY) \;\!\!\big{)} .
\]
Moreover, it immediately follows that $\bm{\Lambda}_{\bm{\rho}} \;\!\! \ll \bm{\Lambda}_{\bm{\pi}}$, so that
\begin{equation}
\label{eq:Lambda_rho<<Lambda_pi} 
W_{\:\!\! p} \big{(}\:\! \bm{\mu}^{\;\!\! \star} \:\!\!,\;\!\! \bm{\nu} \big{)} \;\!\!
\ge
\epsilon_{\;\!\! \bm{\pi}}^{\;\!\! \star} \;\!\!
\quad
\text{$\bm{\Lambda}_{\bm{\rho}}$-\;\!a.s.}
\end{equation}

Applying the standard integration formula for image measures to~\eqref{eq:R_phi,p^star}, we obtain 
\begin{equation}
\label{eq:R_phi,p^star-Lambda} 
\cR_{\;\!\! \phi \;\!\!, p \;\!\!}^{\star \;\!\!}(\bm{\rho} \:\!)
= \;\!\!
\int\nolimits_{\ccP_{\;\!\! p}(\cY)} \!
\phi \:\! \big{(}
\bracks*{\;\!\! W_{\:\!\! p} \big{(}\:\! \bm{\mu}^{\;\!\! \star} \:\!\!,\;\!\! \bm{\nu} \big{)} \;\!\! - \epsilon_{\;\!\! \bm{\pi}}^{\;\!\! \star} }_{\;\!\! +}
\big{)} \;\!
\bm{\Lambda}_{\bm{\rho}}(\d\bm{\nu}) \:\! ,
\quad
\bm{\rho} \in \ccP(\Theta) \:\! , \
\bm{\rho} \ll \;\!\! \bm{\pi} .
\end{equation}
Hence, the asymptotic compatibility functional depends on posterior candidates only through their distributions of model laws.
Statistical information is carried by the transport discrepancy from the observational law~$\bm{\mu}^{\;\!\! \star}$\:\!\!, while $\bm{\rho}$ enters only through the push-forward measure~$\bm{\Lambda}_{\bm{\rho}}$.

In particular, any two probability measures $\bm{\rho} , \widetilde{\bm{\rho}} \in \ccP(\Theta)$ with $\bm{\rho} , \widetilde{\bm{\rho}} \ll \;\!\! \bm{\pi}$ are indistinguishable with respect to asymptotic transport compatibility whenever they share the same induced law; namely,
\[
\bm{\Lambda}_{\bm{\rho}} \;\!\!
=
\bm{\Lambda}_{\widetilde{\bm{\rho}}} \;\!\!
\quad
\Longrightarrow
\quad
\cR_{\;\!\! \phi \;\!\!, p \;\!\!}^{\star \;\!\!}(\bm{\rho} \:\!)
= \:\!
\cR_{\;\!\! \phi \;\!\!, p \;\!\!}^{\star \;\!\!}(\widetilde{\bm{\rho}}\;\!) .
\]
The measure~\eqref{eq:Lambda_rho} thus serves as a natural transport-geometric representation of each posterior candidate.

The preceding representation reveals a direct connection between the entropy-regularized variational structure developed in Subsection~\ref{subsec:compatibility-variational} and the induced distributions of model laws.
The next result makes this connection precise by establishing an exact reduction of relative entropy~\eqref{eq:KL} under the map~$\Phi$.

\begin{proposition}
\label{prop:KL-Phi} 

Under the standing assumptions of this subsection, for every $\bm{\Lambda} \in \ccP \big{(} \ccP_{\;\!\! p}(\cY) \big{)}$,
\begin{equation}
\label{eq:KL-Phi}
\inf_{\substack{
\bm{\rho} \:\! \in \:\! \ccP(\Theta)
\\
\Phi_{\:\! \#} \bm{\rho} \;\! = \:\! \bm{\Lambda}
}} \;\!\!
\mathrm{KL}\:\!(\bm{\rho} \:\!\! \mid \:\!\! \bm{\pi})
=
\mathrm{KL}\:\!(\bm{\Lambda} \:\!\! \mid \:\!\! \bm{\Lambda}_{\bm{\pi}}) \:\! ,
\end{equation}
where the divergence on the right-hand side is understood between probability measures on $\ccP_{\;\!\! p}(\cY)$.
Whenever $\bm{\Lambda} \ll \;\!\! \bm{\Lambda}_{\bm{\pi}}$, the infimum is attained by the measure $\bm{\rho}_{\bm{\Lambda}} \;\!\! \in \ccP(\Theta)$ satisfying $\bm{\rho}_{\bm{\Lambda}} \;\!\! \ll \;\!\! \bm{\pi}$ and 
\begin{equation}
\label{eq:rho_lambda}
\frac{\d\bm{\rho}_{\bm{\Lambda}}}{\d\bm{\pi}}
=
\frac{\d\bm{\Lambda}}{\d\bm{\Lambda}_{\bm{\pi}} \;\!\!} \;\!\!
\circ \Phi
\quad
\text{$\bm{\pi}$-\:\!a.s.}
\end{equation}

\end{proposition}

\begin{proof}

Let $\bm{\Lambda} \in \ccP \big{(} \ccP_{\;\!\! p}(\cY) \big{)}$.
For every $\bm{\rho} \in \ccP(\Theta)$ such that
$\Phi_{\:\! \#} \:\! \bm{\rho} = \;\!\! \bm{\Lambda}$, i.e., $\bm{\Lambda}_{\bm{\rho}} \;\!\! = \bm{\Lambda}$ in the notation of~\eqref{eq:Lambda_rho}, the data-processing inequality for relative entropy yields
\[
\mathrm{KL}\:\!(\bm{\Lambda} \:\!\! \mid \:\!\! \bm{\Lambda}_{\bm{\pi}})
\le
\mathrm{KL}\:\!(\bm{\rho} \:\!\! \mid \:\!\! \bm{\pi})
\]
(see, e.g.,~\cite{dupuis-ellis_1997}).
Minimizing over these measures establishes the ``\:\!$\ge$''~inequality in~\eqref{eq:KL-Phi}.

For the reverse inequality, if $\bm{\Lambda} \not\ll \;\!\! \bm{\Lambda}_{\bm{\pi}}$, then $\mathrm{KL}\:\!(\bm{\Lambda} \:\!\! \mid \:\!\! \bm{\Lambda}_{\bm{\pi}}) = \infty$, and the same bound in turn implies that
\[
\mathrm{KL}\:\!(\bm{\rho} \:\!\! \mid \:\!\! \bm{\pi}) = \infty 
\]
for all $\bm{\rho} \in \ccP(\Theta)$ with
$\Phi_{\:\! \#} \:\! \bm{\rho} = \;\!\! \bm{\Lambda}$.
Thus,~\eqref{eq:KL-Phi} follows directly in this case.

Suppose now that $\bm{\Lambda} \ll \;\!\! \bm{\Lambda}_{\bm{\pi}}$.
By the integration formula for image measures, 
it is readily seen that~\eqref{eq:rho_lambda} defines a probability measure $\bm{\rho}_{\bm{\Lambda}} \;\!\! \in \ccP(\Theta)$;
moreover, $\Phi_{\:\! \#} \:\! \bm{\rho}_{\bm{\Lambda}} \;\!\! = \;\!\! \bm{\Lambda}$ since, for every
$\cE \in \cB \big{(} \ccP_{\;\!\! p}(\cY) \;\!\!\big{)}$,
\[
\Phi_{\:\! \#} \:\! \bm{\rho}_{\bm{\Lambda}} \bracks*{ \cE }
= \;\!\!
\int\nolimits_{\Theta} \;\!\!
\1_{\;\!\! \cE} \big{(} \Phi(\bm{\theta}) \;\!\!\big{)} \:\!
\frac{\d\bm{\Lambda}}{\d\bm{\Lambda}_{\bm{\pi}} \;\!\!} \big{(} \Phi(\bm{\theta}) \;\!\!\big{)}
\;\! \bm{\pi}(\d\bm{\theta})
= \;\!\!
\int\nolimits_{\;\!\! \cE}
\frac{\d\bm{\Lambda}}{\d\bm{\Lambda}_{\bm{\pi}} \;\!\!}(\bm{\nu})
\;\! \bm{\Lambda}_{\bm{\pi}}(\d\bm{\nu})
=
\bm{\Lambda} \bracks{ \cE } .
\]
Finally, a simple calculation based on~\eqref{eq:rho_lambda} and this identity gives
\[
\mathrm{KL}\:\!(\bm{\rho}_{\bm{\Lambda}} \! \mid \:\!\! \bm{\pi})
=
\mathrm{KL}\:\!(\bm{\Lambda} \:\!\! \mid \:\!\! \bm{\Lambda}_{\bm{\pi}}) \:\! ,
\]
which concludes the proof. \qedhere

\end{proof}

In conjunction with the push-forward representation~\eqref{eq:R_phi,p^star-Lambda}, Proposition~\ref{prop:KL-Phi} provides an exact reduction of an entropy-regularized asymptotic compatibility problem, analogous to~\eqref{eq:inf_l-KL}.
For this purpose, set 
\begin{equation}
\label{eq:c_phi,p^star} 
c_{\:\!\! \phi \;\!\!, p \;\!\!}^{\star \;\!\!}(\bm{\nu})
\coloneqq
\phi \:\! \big{(}
\bracks*{\;\!\! W_{\:\!\! p} \big{(}\:\! \bm{\mu}^{\;\!\! \star} \:\!\!,\;\!\! \bm{\nu} \big{)} \;\!\! - \epsilon_{\;\!\! \bm{\pi}}^{\;\!\! \star} }_{\;\!\! +}
\big{)} \:\! ,
\quad
\bm{\nu} \in \ccP_{\;\!\! p}(\cY) \:\! .
\end{equation}

\begin{corollary}
\label{cor:KL-Phi} 

Under the standing assumptions of this subsection,
\begin{equation}
\label{eq:variational-pushforward} 
\inf_{\substack{
\bm{\rho} \:\! \in \:\! \ccP(\Theta)
\\
\bm{\rho} \:\! \ll \bm{\pi}
}} \;\!\!
\braces[\big]{
\cR_{\;\!\! \phi \;\!\!, p \;\!\!}^{\star \;\!\!}(\bm{\rho} \:\!)
+
\mathrm{KL}\:\!(\bm{\rho} \:\!\! \mid \:\!\! \bm{\pi})
}
= \:\!\!
\inf_{\bm{\Lambda} \in \:\! \ccP ( \ccP_{\;\!\! p}\;\!\!(\cY) \:\!\!)}
\braces*{
\int\nolimits_{\ccP_{\;\!\! p}(\cY)} \!
c_{\:\!\! \phi \;\!\!, p \;\!\!}^{\star \;\!\!}(\bm{\nu}) \;\!
\bm{\Lambda}(\d\bm{\nu}) \;\!\!
+
\mathrm{KL}\:\!(\bm{\Lambda} \:\!\! \mid \:\!\! \bm{\Lambda}_{\bm{\pi}}) \;\!\!
} .
\end{equation}
The respective problems admit unique minimizers $\bm{\rho}_{\phi \;\!\!, p \;\!\!}^{\star \;\!\!} \:\!\! \ll \;\!\! \bm{\pi}$ and $\bm{\Lambda}_{\phi \;\!\!, p \;\!\!}^{\! \star \;\!\!} \:\!\! = \Phi_{\:\! \#} \:\! \bm{\rho}_{\phi \;\!\!, p \;\!\!}^{\star \;\!\!} \:\!\! \ll \;\!\! \bm{\Lambda}_{\bm{\pi}}$ with
\begin{equation}
\label{eq:rho-Lambda_phi,p^star} 
\frac{\d\bm{\Lambda}_{\phi \;\!\!, p \;\!\!}^{\! \star \;\!\!}}{\d\bm{\Lambda}_{\bm{\pi}} \;\!\!}
\propto
e^{\;\!\! - \:\! c_{\:\!\! \phi \;\!\!, p \;\!\!}^{\star \;\!\!} } \!
\quad
\text{$\bm{\Lambda}_{\bm{\pi}}$-\:\!a.s.},
\qquad
\frac{\d\bm{\rho}_{\phi \;\!\!, p \;\!\!}^{\star \;\!\!}}{\d\bm{\pi}}
=
\frac{\d\bm{\Lambda}_{\phi \;\!\!, p \;\!\!}^{\! \star \;\!\!}}{\d\bm{\Lambda}_{\bm{\pi}} \;\!\!} \;\!\!
\circ \Phi
\quad
\text{$\bm{\pi}$-\:\!a.s.}
\end{equation}
\end{corollary}

\begin{proof}

For every $\bm{\rho} \in \ccP(\Theta)$ with $\bm{\rho} \ll \;\!\! \bm{\pi}$, by~\eqref{eq:R_phi,p^star-Lambda} and~\eqref{eq:c_phi,p^star},
\[
\cR_{\;\!\! \phi \;\!\!, p \;\!\!}^{\star \;\!\!}(\bm{\rho} \:\!)
= \;\!\!
\int\nolimits_{\ccP_{\;\!\! p}(\cY)} \!
c_{\:\!\! \phi \;\!\!, p \;\!\!}^{\star \;\!\!}(\bm{\nu}) \;\!
\bm{\Lambda}_{\bm{\rho}}(\d\bm{\nu}) \:\! .
\]
Proposition~\ref{prop:KL-Phi} identifies, for each prescribed $\bm{\Lambda} \in \ccP \big{(} \ccP_{\;\!\! p}(\cY) \big{)}$, the least relative entropy $\mathrm{KL}\:\!(\bm{\Lambda} \:\!\! \mid \:\!\! \bm{\Lambda}_{\bm{\pi}})$ among all $\bm{\rho} \in \ccP(\Theta)$ such that $\Phi_{\:\! \#} \:\! \bm{\rho} \equiv \bm{\Lambda}_{\bm{\rho}} \;\!\! = \;\!\! \bm{\Lambda}$, and thus~\eqref{eq:variational-pushforward} follows upon taking the infimum over $\bm{\Lambda}$.

The Gibbs variational principle on $\ccP_{\;\!\! p}(\cY)$, with reference measure $\bm{\Lambda}_{\bm{\pi}}$\;\!\! and cost $c_{\:\!\! \phi \;\!\!, p \;\!\!}^{\star \;\!\!}$, selects the unique minimizer $\bm{\Lambda}_{\phi \;\!\!, p \;\!\!}^{\! \star \;\!\!} \:\!\! \ll \;\!\! \bm{\Lambda}_{\bm{\pi}}$, whose density satisfies the first relation in~\eqref{eq:rho-Lambda_phi,p^star} (see, e.g.,~\cite[Proposition~1.4.2]{dupuis-ellis_1997}).
Applying~\eqref{eq:rho_lambda} to $\bm{\Lambda}_{\phi \;\!\!, p \;\!\!}^{\! \star \;\!\!}$\;\!\! then recovers the unique minimizer $\bm{\rho}_{\phi \;\!\!, p \;\!\!}^{\star \;\!\!} \:\!\! \ll \;\!\! \bm{\pi}$, with density as specified in~\eqref{eq:rho-Lambda_phi,p^star}. \qedhere

\end{proof}

In parameter space, the minimizer in Corollary~\ref{cor:KL-Phi} has density
\[
\frac{\d\bm{\rho}_{\phi \;\!\!, p \;\!\!}^{\star \;\!\!}}{\d\bm{\pi}}(\bm{\theta})
=
\frac{
\exp\:\!\braces*{-\:\!\phi \:\! \big{(}
\bracks*{\;\!\! \cT_{\bm{\theta}}^{\;\! \star}\:\!\! - \epsilon_{\;\!\! \bm{\pi}}^{\;\!\! \star} }_{\;\!\! +}
\big{)}\;\!\!}
}{
\int\nolimits_{\:\! \Theta}
\exp\:\!\braces*{-\:\!\phi \:\! \big{(}
\bracks*{\;\!\! \cT_{\bm{\theta}}^{\;\! \star}\:\!\! - \epsilon_{\;\!\! \bm{\pi}}^{\;\!\! \star} }_{\;\!\! +}
\big{)}\;\!\!}
\;\!\bm{\pi}(\d\bm{\theta})
}.
\]
Replacing $\phi$ by $\lambda\phi$, $\lambda>0$, controls the strength of penalization and generally changes the resulting Gibbs posterior.
Because $\phi$ and the discrepancy are finite, this density is strictly positive $\bm{\pi}$-a.s.; it does not depend on the chosen ABC tolerance $\tau$.
If $\bm{\pi}[\Theta\setminus\Theta_\tau]>0$, it therefore differs from the conditional prior, whose density is $\1_{\Theta_\tau}/\bm{\pi}[\Theta_\tau]$.
The exact reduction in~\eqref{eq:variational-pushforward} concerns this additional penalty-based optimization problem.
Since this Gibbs update depends on the limiting observational law $\bm{\mu}^{\;\!\! \star}$\:\!\!, it is a population-level construction; an implementable version would require a data-dependent approximation, whose properties are not investigated here.

A further consequence of Corollary~\ref{cor:KL-Phi} is that the variational update does not distinguish between parameters inducing the same model law.
Indeed, by~\eqref{eq:rho-Lambda_phi,p^star}, the density of $\bm{\rho}_{\phi \;\!\!, p \;\!\!}^{\star \;\!\!}$\;\!\! with respect to $\bm{\pi}$ depends on $\bm{\theta}$ only through $\Phi(\bm{\theta})$.
Whenever regular conditional distributions with respect to $\Phi$ exist, the posterior minimizer therefore preserves the prior conditional distribution given $\Phi$.

\begin{remark}

If $\phi$ is lower semicontinuous (i.e., left-continuous), then, by standard arguments, $\cR_{\;\!\! \phi \;\!\!, p \;\!\!}^{\star \;\!\!}$ is lower semicontinuous with respect to weak convergence of the induced measures $\bm{\Lambda}_{\bm{\rho}}$.
Thus, the limiting incompatibility cost cannot exceed the lower limit of the approximating values.

Furthermore, linearity of integration implies that $\cR_{\;\!\! \phi \;\!\!, p \;\!\!}^{\star \;\!\!}$ is affine with respect to posterior mixtures.
This reflects the fact that posterior randomization is evaluated by averaging transport-based compatibility.

\end{remark}

The reformulation~\eqref{eq:R_phi,p^star-Lambda} naturally calls for considering the second-order Wasserstein geometry of
\begin{equation}
\label{eq:P_p-P_p-Y} 
\ccP_{\;\!\! p} \big{(} \ccP_{\;\!\! p}(\cY) \big{)} .
\end{equation}
We equip this space with the $p$\:\!-Wasserstein distance $\mathbb{W}_{\! p}$\;\!\! induced by the ground metric~$W_{\:\!\! p}$\;\!\! on $\ccP_{\;\!\! p}(\cY)$: 
\begin{equation*}
\mathbb{W}_{\! p}(\bm{\Lambda}_1 ,\;\!\! \bm{\Lambda}_2)
\coloneqq
\bracks*{
\inf_{\bm{\Gamma} \:\! \in \;\! \Pi(\bm{\Lambda}_1 \;\!\!,\:\! \bm{\Lambda}_2)}
\int\nolimits_{\ccP_{\;\!\! p}(\cY) \times \:\! \ccP_{\;\!\! p}(\cY)} \;\!\!\!
W_{\:\!\! p}^{\:\! p \;\!\!}(\bm{\mu} \:\!,\;\!\! \bm{\nu}) \, \bm{\Gamma}(\d\bm{\mu} \:\!,\;\!\! \d\bm{\nu})
}^{1/p} \! ,
\quad
\bm{\Lambda}_1 ,\;\!\! \bm{\Lambda}_2 \in \ccP_{\;\!\! p} \big{(} \ccP_{\;\!\! p}(\cY) \big{)} .
\end{equation*}

This lifted geometry admits a natural stability estimate for asymptotic compatibility, offering a first indication of its relevance to the present variational analysis.
Specifically, whenever the penalty function is Lipschitz continuous, the associated asymptotic compatibility functional inherits the same regularity with respect to the second-order $p$\:\!-Wasserstein metric between the induced laws.
An analogous conclusion can be drawn for Hölder-continuous penalty functions.

\begin{proposition}
\label{prop:R_phi,p^star-stability} 

Under the standing assumptions of this subsection, suppose that $\phi$ is $L_{\phi}$\;\!\!-\:\!Lipschitz continuous for some $L_{\phi} \;\!\! \in \mo{]}\:\! 0 \:\!,\:\!\! \infty \mc{[}$.
Then, for every $\bm{\rho} , \widetilde{\bm{\rho}} \in \ccP(\Theta)$ with $\bm{\rho} , \widetilde{\bm{\rho}} \ll \;\!\! \bm{\pi}$ and
\begin{equation}
\label{eq:Lambda_rho-P_p} 
\bm{\Lambda}_{\bm{\rho}} \:\!,\;\!\! \bm{\Lambda}_{\widetilde{\bm{\rho}}} \;\!\!
\in
\ccP_{\;\!\! p} \big{(} \ccP_{\;\!\! p}(\cY) \big{)} ,
\end{equation}
one has
\[
\abs*{\:\!
\cR_{\;\!\! \phi \;\!\!, p \;\!\!}^{\star \;\!\!}(\bm{\rho} \:\!)
- \:\!
\cR_{\;\!\! \phi \;\!\!, p \;\!\!}^{\star \;\!\!}(\widetilde{\bm{\rho}}\;\!)
} \:\!
\le
L_{\phi} \;\!\!
\mathbb{W}_{\! p} \big{(} \bm{\Lambda}_{\bm{\rho}} \:\!,\;\!\! \bm{\Lambda}_{\widetilde{\bm{\rho}} \:\!} \big{)} .
\]

\end{proposition}

\begin{proof}

For any $\bm{\Gamma} \;\!\! \in \Pi \big{(} \bm{\Lambda}_{\bm{\rho}} \:\!,\;\!\! \bm{\Lambda}_{\widetilde{\bm{\rho}} \:\!} \big{)}$, by combining~\eqref{eq:R_phi,p^star-Lambda}, the $L_{\phi}$\;\!\!-\:\!Lipschitz continuity of~$\phi$, the (reverse) triangle inequality for $W_{\:\!\! p}$, and Hölder's inequality, we deduce
\[
\begin{split}
\abs*{\:\!
\cR_{\;\!\! \phi \;\!\!, p \;\!\!}^{\star \;\!\!}(\bm{\rho} \:\!)
- \:\!
\cR_{\;\!\! \phi \;\!\!, p \;\!\!}^{\star \;\!\!}(\widetilde{\bm{\rho}}\;\!)
} \:\!
&\le L_{\phi} \;\!\! \int\nolimits_{\ccP_{\;\!\! p}(\cY) \times \:\! \ccP_{\;\!\! p}(\cY)} \;\!\!\!
W_{\:\!\! p}(\bm{\mu} \:\!,\;\!\! \bm{\nu}) \, \bm{\Gamma}(\d\bm{\mu} \:\!,\;\!\! \d\bm{\nu}) \\[0.25ex]
&\le L_{\phi} \:\! \bracks*{
\int\nolimits_{\ccP_{\;\!\! p}(\cY) \times \:\! \ccP_{\;\!\! p}(\cY)} \;\!\!\!
W_{\:\!\! p}^{\:\! p \;\!\!}(\bm{\mu} \:\!,\;\!\! \bm{\nu}) \, \bm{\Gamma}(\d\bm{\mu} \:\!,\;\!\! \d\bm{\nu})
}^{1/p} \! .
\end{split}
\]
Since $\bm{\Gamma}$ was arbitrary, taking the infimum over all couplings in $\Pi \big{(} \bm{\Lambda}_{\bm{\rho}} \:\!,\;\!\! \bm{\Lambda}_{\widetilde{\bm{\rho}} \:\!} \big{)}$ yields the desired estimate. \qedhere

\end{proof}

\begin{remark}

The $p$-moment requirement~\eqref{eq:Lambda_rho-P_p} is fairly mild.
In particular, whenever it holds for the prior~$\bm{\pi}$, it remains valid for every posterior whose density with respect to~$\bm{\pi}$ is essentially bounded.
More generally, this property follows under suitable Hölder-type conditions on the Radon--Nikodym derivative.

\end{remark}

A particularly transparent interpretation emerges when the minimal compatibility threshold vanishes, under a distinguished power-type specialization.
In this case, the transport-induced compatibility functional is exactly the $p$-th power of the second-order $p$\:\!-Wasserstein distance between the Dirac mass concentrated at the observational law and the induced law of model distributions.

\begin{proposition}
\label{prop:R_phi,p^star-representation} 

Under the standing assumptions of this subsection, suppose that
\begin{equation*}
\epsilon_{\;\!\! \bm{\pi}}^{\;\!\! \star} \;\!\!
=
0
\end{equation*}
and
\[
\phi(x) = x^{\:\! p} \;\!\! ,
\quad
x \in \R_{\:\! +} .
\]
Then, for every $\bm{\rho} \in \ccP(\Theta)$ with $\bm{\rho} \ll \;\!\! \bm{\pi}$ and $\bm{\Lambda}_{\bm{\rho}} \;\!\! \in \ccP_{\;\!\! p} \big{(} \ccP_{\;\!\! p}(\cY) \big{)}$, 
\[
\cR_{\;\!\! \phi \;\!\!, p \;\!\!}^{\star \;\!\!}(\bm{\rho} \:\!)
=
\mathbb{W}_{\! p}^{\:\! p \;\!\!} \big{(} \bm{\delta}_{\bm{\mu}^{\;\!\! \star} \;\!\!} ,\;\!\! \bm{\Lambda}_{\bm{\rho}} \big{)} .
\]

\end{proposition}

\begin{proof}

Let $\bm{\rho} \in \ccP(\Theta)$, $\bm{\rho} \ll \;\!\! \bm{\pi}$, satisfy $\bm{\Lambda}_{\bm{\rho}} \;\!\! \in \ccP_{\;\!\! p} \big{(} \ccP_{\;\!\! p}(\cY) \big{)}$.
Since the only coupling between $\bm{\delta}_{\bm{\mu}^{\;\!\! \star} \:\!\!}$ and $\bm{\Lambda}_{\bm{\rho}}$\;\!\! is
\[
\bm{\delta}_{\bm{\mu}^{\;\!\! \star} \:\!\!} \;\!\!
\otimes
\bm{\Lambda}_{\bm{\rho}} \:\! ,
\]
we directly obtain
\[
\mathbb{W}_{\! p}^{\:\! p \;\!\!} \big{(} \bm{\delta}_{\bm{\mu}^{\;\!\! \star} \;\!\!} ,\;\!\! \bm{\Lambda}_{\bm{\rho}} \big{)}
= \;\!\!
\int\nolimits_{\ccP_{\;\!\! p}(\cY)} \;\!\!\!
W_{\:\!\! p}^{\:\! p \;\!\!} \big{(}\:\! \bm{\mu}^{\;\!\! \star} \:\!\!,\;\!\! \bm{\nu} \big{)} \;\!
\bm{\Lambda}_{\bm{\rho}}(\d\bm{\nu})
\equiv \:\!
\cR_{\;\!\! \phi \;\!\!, p \;\!\!}^{\star \;\!\!}(\bm{\rho} \:\!) \:\! ,
\]
as asserted. \qedhere

\end{proof}

Proposition~\ref{prop:R_phi,p^star-representation} identifies the asymptotic incompatibility of a posterior distribution with the $p$-th power of the second-order $p$\:\!-Wasserstein distance between $\bm{\delta}_{\bm{\mu}^{\;\!\! \star} \:\!\!}$ and its transport-geometric representative.
In this sense, asymptotic compatibility is governed by the geometry of the space~\eqref{eq:P_p-P_p-Y}.

The same geometric principle extends naturally to arbitrary compatibility thresholds, with the reference observational law replaced by a closed Wasserstein ball.
With this aim, let
\begin{equation}
\label{eq:B^star} 
\widebar{B}_{\! \epsilon_{\;\!\! \bm{\pi}}^{\;\!\! \star} }(\bm{\mu}^{\;\!\! \star})
\coloneqq \;\!\!
\Set{\:\!\!
\bm{\eta} \in \ccP_{\;\!\! p}(\cY)
|
W_{\:\!\! p} \big{(}\:\! \bm{\mu}^{\;\!\! \star} \:\!\!,\;\!\! \bm{\eta} \big{)} \;\!\!
\le
\epsilon_{\;\!\! \bm{\pi}}^{\;\!\! \star} \;\!\!
\:\!\!} \:\!\! ,
\end{equation}
and let $\bd_{\widebar{B}_{\;\!\!\! \epsilon_{\;\!\! \bm{\pi}}^{\;\!\! \star} \;\!\!}(\bm{\mu}^{\;\!\! \star}\:\!\!)} \;\!\! \colon \ccP_{\;\!\! p}(\cY) \to \R_{\:\! +}$\;\!\! be the associated distance function, given by
\begin{equation}
\label{eq:d_B^star} 
\bd_{\widebar{B}_{\;\!\!\! \epsilon_{\;\!\! \bm{\pi}}^{\;\!\! \star} \;\!\!}(\bm{\mu}^{\;\!\! \star}\:\!\!)}\;\!\!(\;\! \bm{\cdot} \;\!)
\coloneqq \!
\inf_{\bm{\eta} \:\! \in \:\! \widebar{B}_{\;\!\!\! \epsilon_{\;\!\! \bm{\pi}}^{\;\!\! \star} \;\!\!}(\bm{\mu}^{\;\!\! \star}\:\!\!)} \;\!\!\!
W_{\:\!\! p} \big{(}\, \bm{\cdot} \,,\;\!\! \bm{\eta} \big{)} \:\! .
\end{equation}

\begin{theorem}
\label{thm:R_phi,p^star-representation} 

Under the standing assumptions of this subsection, suppose that the space $\rounds{ \ccP_{\;\!\! p}(\cY) \:\!,\:\!\! W_{\:\!\! p} }$ is geodesic and that $\phi(x) = x^{\:\! p}$\;\!\!, $x \in \R_{\:\! +}$.
Then, for every $\bm{\rho} \in \ccP(\Theta)$ with $\bm{\rho} \ll \;\!\! \bm{\pi}$, 
\[
\cR_{\;\!\! \phi \;\!\!, p \;\!\!}^{\star \;\!\!}(\bm{\rho} \:\!)
= \;\!\!
\int\nolimits_{\ccP_{\;\!\! p}(\cY)} \!
\bd_{\widebar{B}_{\;\!\!\! \epsilon_{\;\!\! \bm{\pi}}^{\;\!\! \star} \;\!\!}(\bm{\mu}^{\;\!\! \star}\:\!\!)}^{\;\! p}\:\!\!(\bm{\nu}) \;\!
\bm{\Lambda}_{\bm{\rho}}(\d\bm{\nu}) \:\! .
\]

\end{theorem}

\begin{proof}

In view of~\eqref{eq:Lambda_rho<<Lambda_pi} and~\eqref{eq:R_phi,p^star-Lambda}, it suffices to prove that, for every $\bm{\nu} \in \ccP_{\;\!\! p}(\cY)$,
\begin{equation}
\label{eq:thm-bb-W} 
\bd_{\widebar{B}_{\;\!\!\! \epsilon_{\;\!\! \bm{\pi}}^{\;\!\! \star} \;\!\!}(\bm{\mu}^{\;\!\! \star}\:\!\!)}\:\!\!(\bm{\nu})
=
\bracks*{\;\!\! W_{\:\!\! p} \big{(}\:\! \bm{\mu}^{\;\!\! \star} \:\!\!,\;\!\! \bm{\nu} \big{)} \;\!\! - \epsilon_{\;\!\! \bm{\pi}}^{\;\!\! \star} }_{+} .
\end{equation}

The inequality ``\:\!$\ge$'' in~\eqref{eq:thm-bb-W} is an immediate consequence of the metric structure of $\ccP_{\;\!\! p}(\cY)$.
Indeed, by the triangle inequality for $W_{\:\!\! p}$, for every $\bm{\nu} \in \ccP_{\;\!\! p}(\cY)$ and $\bm{\eta} \in \widebar{B}_{\! \epsilon_{\;\!\! \bm{\pi}}^{\;\!\! \star} }(\bm{\mu}^{\;\!\! \star})$ (see~\eqref{eq:B^star}),
\[
W_{\:\!\! p}(\bm{\nu} \:\!,\;\!\! \bm{\eta})
\ge
W_{\:\!\! p} \big{(}\:\! \bm{\mu}^{\;\!\! \star} \:\!\!,\;\!\! \bm{\nu} \big{)} \;\!\!
-
W_{\:\!\! p} \big{(}\:\! \bm{\mu}^{\;\!\! \star} \:\!\!,\;\!\! \bm{\eta} \big{)} \;\!\!
\ge
W_{\:\!\! p} \big{(}\:\! \bm{\mu}^{\;\!\! \star} \:\!\!,\;\!\! \bm{\nu} \big{)} \;\!\!
-
\epsilon_{\;\!\! \bm{\pi}}^{\;\!\! \star} \:\! ,
\]
whence the claim follows by~\eqref{eq:d_B^star}.

The converse inequality is precisely where the assumed geodesic structure of $\ccP_{\;\!\! p}(\cY)$ comes into play. 
Fix $\bm{\nu} \in \ccP_{\;\!\! p}(\cY)$.
If $W_{\:\!\! p} \big{(}\:\! \bm{\mu}^{\;\!\! \star} \:\!\!,\;\!\! \bm{\nu} \big{)} \;\!\! \le \epsilon_{\;\!\! \bm{\pi}}^{\;\!\! \star}$, then $\bm{\nu} \in \widebar{B}_{\! \epsilon_{\;\!\! \bm{\pi}}^{\;\!\! \star} }(\bm{\mu}^{\;\!\! \star})$, and both sides of~\eqref{eq:thm-bb-W} vanish.
Otherwise, set
\[
s_{\bm{\nu}}^{\star}
\coloneqq
\frac{\epsilon_{\;\!\! \bm{\pi}}^{\;\!\! \star} \;\!\!}{W_{\:\!\! p} \big{(}\:\! \bm{\mu}^{\;\!\! \star} \:\!\!,\;\!\! \bm{\nu} \big{)}}
\in
\mo{[} 0 \:\!,\:\!\! 1 \mc{[} \:\! ,
\]
and let $\rounds{ \bm{\eta}_{\:\! s} }_{s \:\! \in \:\! \mo{[}\;\!\! 0 \:\!, 1 \:\!\!\mc{]}}$\;\!\! be a constant-speed Wasserstein geodesic in $\ccP_{\;\!\! p}(\cY)$ connecting $\bm{\mu}^{\;\!\! \star}$\:\!\! and $\bm{\nu}$.
Then 
\[
W_{\:\!\! p} \big{(}\:\! \bm{\mu}^{\;\!\! \star} \:\!\!,\;\!\! \bm{\eta}_{\:\! s_{\bm{\nu}}^{\star}} \big{)} \;\!\!
=
s_{\bm{\nu}}^{\star} \:\! W_{\:\!\! p} \big{(}\:\! \bm{\mu}^{\;\!\! \star} \:\!\!,\;\!\! \bm{\nu} \big{)} \;\!\!
\equiv
\epsilon_{\;\!\! \bm{\pi}}^{\;\!\! \star} \;\!\!
\]
and
\[
W_{\:\!\! p} \big{(}\:\! \bm{\nu} \:\!,\;\!\! \bm{\eta}_{\:\! s_{\bm{\nu}}^{\star}} \big{)}
=
\rounds{ 1 - s_{\bm{\nu}}^{\star} \:\!} \:\! W_{\:\!\! p} \big{(}\:\! \bm{\mu}^{\;\!\! \star} \:\!\!,\;\!\! \bm{\nu} \big{)} \;\!\!
\equiv
W_{\:\!\! p} \big{(}\:\! \bm{\mu}^{\;\!\! \star} \:\!\!,\;\!\! \bm{\nu} \big{)} \;\!\!
-
\epsilon_{\;\!\! \bm{\pi}}^{\;\!\! \star} \:\! .
\]
Therefore, by construction, $\bm{\eta}_{\:\! s_{\bm{\nu}}^{\star}} \;\!\! \in \widebar{B}_{\! \epsilon_{\;\!\! \bm{\pi}}^{\;\!\! \star} }(\bm{\mu}^{\;\!\! \star})$ and
\[
\bd_{\widebar{B}_{\;\!\!\! \epsilon_{\;\!\! \bm{\pi}}^{\;\!\! \star} \;\!\!}(\bm{\mu}^{\;\!\! \star}\:\!\!)}\:\!\!(\bm{\nu})
\le
W_{\:\!\! p} \big{(}\:\! \bm{\nu} \:\!,\;\!\! \bm{\eta}_{\:\! s_{\bm{\nu}}^{\star}} \big{)} \;\!\!
=
W_{\:\!\! p} \big{(}\:\! \bm{\mu}^{\;\!\! \star} \:\!\!,\;\!\! \bm{\nu} \big{)} \;\!\!
-
\epsilon_{\;\!\! \bm{\pi}}^{\;\!\! \star} \:\! .
\]
This completes the proof. \qedhere

\end{proof}

\begin{remark}

The geodesic assumption in Theorem~\ref{thm:R_phi,p^star-representation} holds, for example, whenever the underlying space $\rounds{\cY \;\!\!, \bd_{\;\!\! \cY}}$ is geodesic (as the Polish assumption is already in place); see, e.g.,~\cite{mccann_1997,villani_2009,agueh-carlier_2011}.
In particular, this includes the classical Euclidean setting and complete Riemannian manifolds.

\end{remark}

The preceding results provide a natural transport-geometric characterization of asymptotic compatibility.
The functional~\eqref{eq:R_phi,p^star} factors through the push-forward map induced by~$\Phi$ in~\eqref{eq:Phi} and thus depends on posterior candidates only through their induced distributions of model laws in Wasserstein space, while the prior-essential compatibility threshold determines the reference Wasserstein ball.


We conclude by illustrating these constructions through several classes of statistical models.
The transport representations in the first two examples are classical; see, e.g.,~\cite{villani_2009,santambrogio_2015}.

\begin{example}[Gaussian model laws]
\label{ex:gaussian-models} 

Let $\cY = \R^{m_{\;\!\! \cY}}$\:\!\! with its Euclidean metric, and suppose that
\[
\bm{\mu}^{\;\!\! \star} \;\!\!
=
\cN
\big{(}
M^{\star}
\:\!\!,\;\!\!
\Sigma^{\star}
\big{)}
\]
and, for every $\bm{\theta} \in \Theta$,
\[
\bm{\nu}_{\bm{\theta}}
=
\cN
\big{(}
M_{\;\!\! \bm{\theta}}
\:\!,\;\!\!
\Sigma_{\:\! \bm{\theta}}
\big{)} ,
\]
where $M_{\;\!\! \bm{\theta}} , M^{\star} \:\!\! \in \R^{m_{\;\!\! \cY}}$\:\!\! are mean vectors and $\Sigma_{\:\! \bm{\theta}} , \Sigma^{\star} \:\!\! \in \R^{m_{\;\!\! \cY} \times \:\! m_{\;\!\! \cY}}$\;\!\! are covariance matrices.
Then, for $p = 2$, 
\[
\big{(} \cT_{\bm{\theta}}^{\;\! \star} \;\!\!\big{)}^{2} \:\!\!
\equiv
W_{\;\!\! 2}^{\;\! 2 \;\!\!}
\big{(}\:\! \bm{\mu}^{\;\!\! \star} \:\!\!,\;\!\! \bm{\nu}_{\bm{\theta}} \big{)} \;\!\!
=
\norm*{M^{\star} \:\!\! - M_{\bm{\;\!\! \theta}}}^{\;\! 2} \:\!\!
+
\tr
\bracks*{
\Sigma^{\star} \:\!\!
+
\Sigma_{\:\! \bm{\theta}} \;\!\!
-
2
\bracks*{
\big{(} \Sigma^{\star} \;\!\!\big{)}^{\:\!\! 1/2} \;\! 
\Sigma_{\:\! \bm{\theta}}
\big{(} \Sigma^{\star} \;\!\!\big{)}^{\:\!\! 1/2}
}^{\;\!\! 1/2}
}
\]
(Bures\:\!--\;\!\!Wasserstein formula; see, e.g.,~\cite{gelbrich_1990}).
Consequently, both the asymptotic compatibility profile~\eqref{eq:T_theta^star-p} and the integrand of the associated transport-based compatibility functional~\eqref{eq:R_phi,p^star} admit explicit expressions in terms of the corresponding mean vectors and covariance matrices.

The univariate case further elucidates the large-sample filtering mechanism developed in Section~\ref{sec:large-sample}.
Let $m_{\:\!\! \cY} \;\!\! = 1$, $\Theta = \R$, $\bm{\theta}_0 \;\!\! \in \R$, and consider
\[
\bm{\mu}^{\;\!\! \star} \;\!\!
=
\cN
\big{(}
\bm{\theta}_0
\:\!,\;\!\!
\sigma_{0}^{\:\! 2}
\big{)}
\]
and, for every $\bm{\theta} \in \R$,
\[
\bm{\nu}_{\bm{\theta}}
=
\cN
\big{(}
\bm{\theta}
,\;\!\!
\sigma^{\:\! 2}
\big{)} ,
\]
where $\sigma_{\:\! 0} , \sigma \in \mo{]}\:\! 0 \:\!,\:\!\! \infty \mc{[}$, not necessarily equal.
Then, still for $p = 2$,
\[
W_{\;\!\! 2}^{\;\! 2 \;\!\!}
\big{(}\:\! \bm{\mu}^{\;\!\! \star} \:\!\!,\;\!\! \bm{\nu}_{\bm{\theta}} \big{)} \;\!\!
=
\big{(} \bm{\theta} - \bm{\theta}_0 \big{)}^{2} \:\!\!
+
\big{(} \sigma - \sigma_{\:\! 0} \big{)}^{2} \:\!\! ,
\]
and hence, from~\eqref{eq:epsilon_Theta^star-min} and~\eqref{eq:Theta^star},
\[
\epsilon_{\Theta}^{\;\!\! \star} \;\!\!
=
\abs[\big]{\:\!
\sigma - \sigma_{\:\! 0}
\:\!} \:\! ,
\quad
\Theta^\star \:\!\!
= \;\!\!
\Set{\!
\bm{\theta}_0
\!} \! .
\]

Suppose further that $\bm{\theta}_0 \;\!\! \in \supp \bm{\pi}$.
Then $\epsilon_{\;\!\! \bm{\pi}}^{\;\!\! \star} \;\!\! = \epsilon_{\Theta}^{\;\!\! \star}$\;\!\! (see~\eqref{eq:epsilon_pi^star} and~\eqref{eq:epsilon_Theta^star}) and, for every supercritical tolerance $\tau > \epsilon_{\Theta}^{\;\!\! \star}$, the corresponding asymptotic compatibility region~\eqref{eq:Theta_tau} is
\[
\Theta_\tau \;\!\!
=
\Set{\!
\bm{\theta} \in \R
|
\abs[\big]{\:\!
\bm{\theta} - \bm{\theta}_0
\:\!}
<
\sqrt{\tau^{\:\! 2} \;\!\! - \big{(} \epsilon_{\Theta}^{\;\!\! \star} \:\!\big{)}^{2} \;\!\!}
\:\!\!} \:\!\! .
\]

If $\sigma = \sigma_{\:\! 0}$, then $\epsilon_{\Theta}^{\;\!\! \star} \;\!\! = 0$ and the model is well specified, in the sense that the limiting observational measure $\bm{\mu}^{\;\!\! \star}$\:\!\! belongs to the model class $\Set{\! \bm{\nu}_{\bm{\theta}} \!}_{\;\!\! \bm{\theta} \:\! \in \:\! \R}$.
In this case, the compatibility region reduces to
\[
\Theta_\tau \;\!\!
=
\Set{\!
\bm{\theta} \in \R
|
\abs[\big]{\:\!
\bm{\theta} - \bm{\theta}_0
\:\!}
<
\tau
\:\!\!} \! .
\]
If instead $\sigma \ne \sigma_{\:\! 0}$, i.e., $\epsilon_{\Theta}^{\;\!\! \star} \;\!\! > 0$, the model is misspecified:
no parameter is asymptotically compatible with a subcritical tolerance $\tau < \epsilon_{\Theta}^{\;\!\! \star}$, while at the critical value $\tau = \epsilon_{\Theta}^{\;\!\! \star}$ the compatibility region is empty.

Under the assumptions of Theorem~\ref{thm:compatibility-filtering}, as $n \to \infty$, the finite-sample compatibility weights~\eqref{eq:w_tau^n} approach the indicator $\1_{\Theta_\tau}$\;\!\! of $\Theta_\tau$\;\!\! in the integrated $L^{\;\!\! 1}$\:\!\! sense~\eqref{eq:w_tau^n-limit}.
The associated compatibility posteriors~\eqref{eq:Pi_tau^n} then converge in total variation, in $\bP$-\:\!probability, to the conditional prior distribution
\[
\bm{\pi} \bracks{\;\! \bm{\cdot} \:\!\! \mid \:\!\! \Theta_{\tau \;\!\!} } \:\! .
\]
Thus, even though $\bm{\theta}_0$\;\!\! is the unique asymptotically optimal parameter, a fixed supercritical tolerance alone does not in general lead to posterior concentration at $\bm{\theta}_0$\;\!\! itself.

For a finite-sample counterpart to this behavior, fix $\bm{\theta}_0 \;\!\! = 0$, $\sigma_{\:\! 0} \;\!\! = 1$, $\sigma = 1.5$, $\tau = 1$, and the prior
\[
\bm{\pi} = \cN (0 \:\!,\! 1) \:\! .
\]
For each $n \in \N^{\:\! *}$\:\!\!, given the observed sample $\bY^{1:\:\!n} \;\!\! = \;\!\! \bracks[\big]{ \bY^1 \:\!\!, \dots ,\;\!\! \bY^n } \;\!\!$, define, on $\Omega$,
\[
\widebar{\bY}_{\;\!\!\! n} \;\!\!
\coloneqq
\frac{1}{n}
\sum_{i \:\! = \:\! 1}^{n} \bY^{\:\! i} \:\!\! ,
\quad
\bS_{\:\! \bY \:\!\!\!,\:\! n}^{\;\! 2} \:\!\!
\coloneqq
\frac{1}{n}
\sum_{i \:\! = \:\! 1}^{n}
\big{(} \bY^{\:\! i} \:\!\! - \widebar{\bY}_{\;\!\!\! n} \big{)}^{2} \:\!\! .
\]
The corresponding probabilistic representation, as in Subsection~\ref{subsec:probabilistic-representations}, is
\[
\bm{\nu}_{n} \big{(} \bY^{1:\:\!n} \big{)} \;\!\!
\coloneqq
\cN \big{(} \widebar{\bY}_{\;\!\!\! n} \:\!,\;\!\! \bS_{\:\! \bY \:\!\!\!,\:\! n}^{\;\! 2} \big{)} \:\! .
\]
For the simulated sample $\bY_{\;\!\! \bm{\theta}}^{1:\:\!n}$\:\!\!, $\bm{\theta} \in \R$, introduce $\widebar{\bY}_{\;\!\!\! \bm{\theta} \:\!\!, n}$\;\!\! and $\bS_{\:\! \bY_{\;\!\! \bm{\theta}} \;\!\!, n}^{\;\! 2}$\;\!\! analogously, and set
\[
\bm{\nu}_{\bm{\theta} \:\!\!, n} \big{(} \bY_{\;\!\! \bm{\theta}}^{1:\:\!n} \big{)} \;\!\!
\coloneqq
\cN \big{(} \widebar{\bY}_{\;\!\!\! \bm{\theta} \:\!\!, n} \:\!,\;\!\! \bS_{\:\! \bY_{\;\!\! \bm{\theta}} \;\!\!, n}^{\;\! 2} \big{)} .
\]
For every $\bm{\theta} \in \R$, the resulting squared finite-sample discrepancy takes the form 
\[
W_{\;\!\! 2}^{\;\! 2 \;\!\!}
\big{(}\:\! \bm{\nu}_{n} \big{(} \bY^{1:\:\!n} \big{)} ,\;\!\! \bm{\nu}_{\bm{\theta} \:\!\!, n} \big{(} \bY_{\;\!\! \bm{\theta}}^{1:\:\!n} \big{)} \;\!\!\big{)} \;\!\!
=
\big{(}
\widebar{\bY}_{\;\!\!\! n} \;\!\!
-
\widebar{\bY}_{\;\!\!\! \bm{\theta} \:\!\!, n}
\big{)}^{2} \:\!\!
+
\big{(}
\bS_{\:\! \bY \:\!\!\!,\:\! n} \;\!\!
-
\bS_{\:\! \bY_{\;\!\! \bm{\theta}} \;\!\!, n}
\big{)}^{2} \:\!\! ,
\]
with $\bS_{\:\! \bY \:\!\!\!,\:\! n}$\;\!\! and $\bS_{\:\! \bY_{\;\!\! \bm{\theta}} \;\!\!, n}$\;\!\! denoting the respective empirical standard deviations.
In particular,
\[
\epsilon_{\Theta}^{\;\!\! \star} \;\!\!
=
\epsilon_{\;\!\! \bm{\pi}}^{\;\!\! \star} \;\!\!
=
1/2 \;\! ,
\quad
\Theta_\tau \;\!\!
=
\mo{]}
-\sqrt{3}/2
\:\!,\!
\sqrt{3}/2
\mc{[} .
\]

For a reproducible finite-sample calculation, generate $1000$ independent $\cN(0,1)$ observations using NumPy's \texttt{Generator(PCG64(20261010))}, and use the first $n$ observations for $n\in\{50,200,1000\}$.
The calculation uses the Gaussian fitted representations above, not the empirical measures themselves.
Write $m_n$ and $s_n$ for the realized sample mean and standard deviation, with divisor $n$ in the sample variance.
For the simulated standard deviation $S$, one has $nS^2/\sigma^2\;\!\!\sim\chi^2_{n-1}$, independently of the simulated sample mean, whose conditional law is $\cN(\theta,\sigma^2/n)$.
Let $f_{n,\sigma}$ be the density of $S$, let $\Phi_0$ denote the standard normal distribution function, and put $r_n(s)=\sqrt{\tau^2-(s-s_n)^2}$ on $I_n=[\max(0,s_n-\tau),s_n+\tau]$.
The compatibility weight therefore has the one-dimensional representation
\begin{equation}
\label{eq:gaussian-quadrature}
w_{\:\!\! \tau}^{\:\! n}(\theta;\by^{1:n})
=\int\nolimits_{\;\!\! I_n}\!\bigg[
\Phi_0\!\left(\frac{m_n+r_n(s)-\theta}{\sigma/\sqrt n}\right) \;\!\! - \Phi_0\!\left(\frac{m_n-r_n(s)-\theta}{\sigma/\sqrt n}\right)
\bigg] f_{n,\sigma}(s)\,\d s. 
\end{equation}
Integrating over the prior $\cN (0 \:\!,\! 1)$ gives the normalizing constant by replacing the two arguments in~\eqref{eq:gaussian-quadrature} with $(m_n\pm r_n(s))/\sqrt{1+\sigma^2/n}$.
Thus no simulated sample means or Monte Carlo estimates of the acceptance probabilities are needed for this example.

We evaluate these integrals by $600$-node Gauss--Legendre quadrature.
Figure~\ref{fig:gaussian-filtering} uses a grid of $1001$ parameter values from $-2.5$ to $2.5$, with spacing $0.005$.
The limiting normalizer is $\bm{\pi}[\Theta_\tau]=0.6135237692$, and the limiting density at zero is $0.6502474727$.
Table~\ref{tab:gaussian-errors} reports the total variation error, computed as one half of the density $L^1$\:\!\! error over $[-8,8]$, with the quadrature split at both compatibility boundaries.
The omitted tail contribution is below $1.1\times10^{-15}$ for these realizations, using $w_{\:\!\! \tau}^{\:\! n}\le1$ and the Gaussian prior tail.
Doubling the quadrature order to $1200$ changes the reported total variation errors by less than $6.6\times10^{-8}$ and the normalizers by less than $1.1\times10^{-13}$.
The observed errors decrease along this simulated sequence; this illustration does not assert monotonicity for arbitrary data realizations or a convergence rate.

\begin{figure}[H]
\centering
\includegraphics[width=\textwidth]{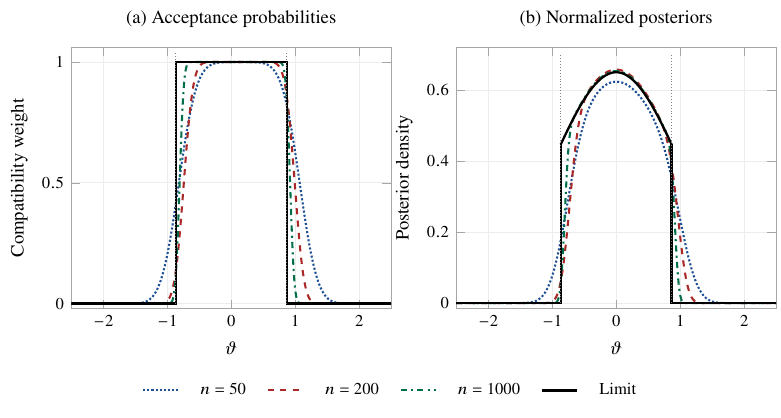}
\caption{Finite-sample compatibility filtering using Gaussian fitted representations in the misspecified model. Left: conditional acceptance weights. Right: normalized posterior densities. The black limits are the indicator of $\Theta_\tau$ and the prior conditioned on $\Theta_\tau$, respectively. Dotted vertical lines mark $\pm\sqrt{3}/2$. The finite-sample curves use deterministic quadrature conditional on the observed samples.}
\label{fig:gaussian-filtering}
\end{figure}

\begin{table}[H]
\centering
\begin{tabular}{r r r r r}
\toprule
$n$ & $m_n$ & $s_n$ & $C_\tau^n$ & Total variation error \\
\midrule
50 & 0.125410 & 1.184947 & 0.639483 & 0.104930 \\
200 & 0.125055 & 1.001095 & 0.606501 & 0.065867 \\
1000 & 0.056368 & 0.994954 & 0.610417 & 0.030386 \\
\bottomrule
\end{tabular}

\caption{Observed summaries, normalizing constants, and total variation errors relative to the conditional prior for the reproducible Gaussian example.}
\label{tab:gaussian-errors}
\end{table}

The supplementary script \texttt{numerics/generate\_gaussian.py} records the generator, software versions, quadrature settings, observation values, and numerical checks, and exports the plotted data and this table.
The supplied figure source uses PGFPlots to reproduce the vector graphic.

The same Gaussian model also illustrates the additional penalty-based update introduced in~\eqref{eq:variational-pushforward}.
For $\lambda>0$, consider the quadratic penalty
\[
\phi_\lambda(t)
=
\lambda t^2 \:\! ,
\quad
t \in \R_{\:\! +} .
\]
By Corollary~\ref{cor:KL-Phi}, the corresponding Gibbs minimizer has density
\[
q_\lambda(\bm{\theta})
=
\frac{1}{Z_\lambda\sqrt{2\pi}}
\exp\Bigg\{
-\frac{\bm{\theta}^{\:\! 2}}{2}
-
\lambda
\bigg(
\sqrt{\bm{\theta}^{\:\! 2}+\frac14}-\frac12
\bigg)^2
\Bigg\} ,
\]
where $Z_\lambda>0$ is the normalizing constant.

Unlike the limiting ABC posterior, which is the prior conditioned on $\Theta_\tau$, the Gibbs posterior assigns positive density to every parameter value.

For $\lambda\in\{1,5,20\}$, Figure~\ref{fig:gaussian-gibbs} and Table~\ref{tab:gaussian-gibbs} compare the resulting distributions with the prior and the limiting ABC posterior.
The reported quantities, computed under each distribution, are the standard deviation, the probability assigned to the complement of $\Theta_\tau$, and the expected excess discrepancy
\[
D(\bm{\theta})
\coloneqq
\sqrt{\bm{\theta}^{\:\! 2}+\frac14}-\frac12 .
\]
Since $D(\bm{\theta}) = \bm{\theta}^{\:\! 2}\:\!\! + O(\bm{\theta}^{\:\! 4})$ near zero, the penalty $\lambda D(\bm{\theta})^2$\:\!\! is locally quartic in the parameter, despite being quadratic in the discrepancy.

These calculations use Gauss--Legendre quadrature on $[-8,8]$, split at the compatibility boundaries, with $600$ nodes per subinterval.
Doubling the quadrature order changes the reported quantities by less than $1.3\times10^{-13}$\:\!\!.
The supplementary script \texttt{numerics/generate\_gaussian.py} also generates the numerical data underlying Figure~\ref{fig:gaussian-gibbs} and Table~\ref{tab:gaussian-gibbs}.

Increasing $\lambda$ reduces both posterior dispersion and expected excess discrepancy, while the Gibbs density remains strictly positive.
In particular, the penalized update need not have greater dispersion than the conditional prior.

The geometric representation in Theorem~\ref{thm:R_phi,p^star-representation} gives a direct interpretation of the penalty:
\[
D(\bm{\theta})^2 \;\!\!
=
\bd_{\widebar{B}_{\;\!\!\! \epsilon_{\;\!\! \bm{\pi}}^{\;\!\! \star} \;\!\!}(\bm{\mu}^{\;\!\! \star}\:\!\!)}^{\;\! 2} \:\!\!
(\bm{\nu}_{\bm{\theta}}) .
\]
Thus, the Gibbs weights are determined by the squared Wasserstein distance of each model law from the reference ball.
The values of $\lambda$ provide a sensitivity analysis rather than a prescribed calibration rule.

\begin{figure}[H]
\centering
\includegraphics[width=\textwidth]{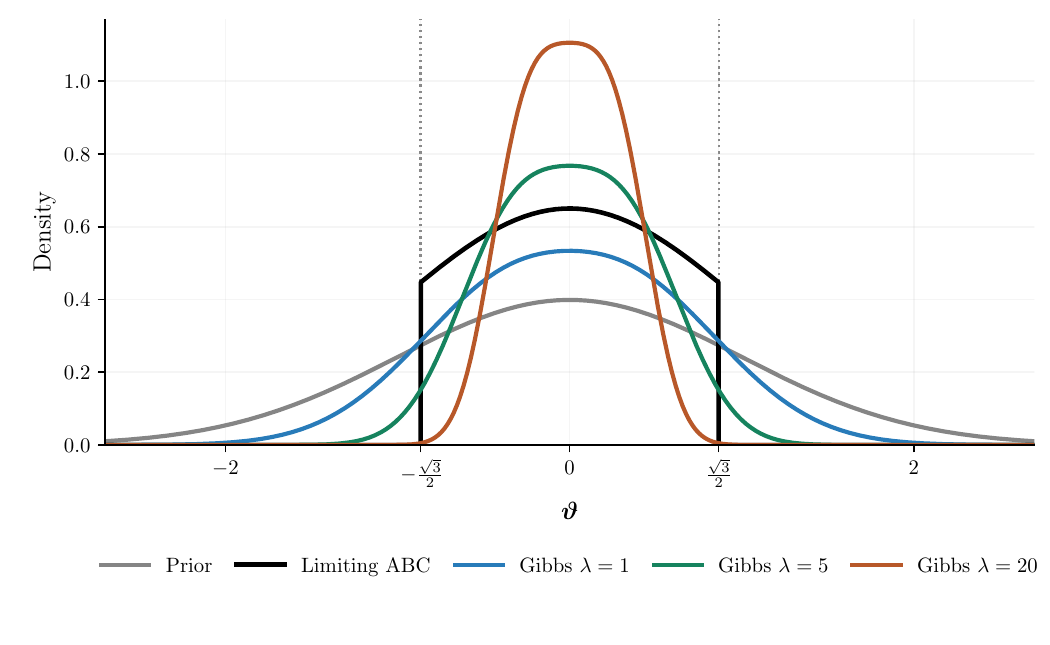}
\caption{Comparison of the prior, the limiting ABC posterior, and the penalty-based Gibbs posteriors for $\lambda\in\{1,5,20\}$ in the misspecified Gaussian model. Dotted vertical lines mark the compatibility boundaries $\pm\sqrt{3}/2$.}
\label{fig:gaussian-gibbs}
\end{figure}

\begin{table}[H]
\centering
\begin{tabular}{lccc}
\hline
Distribution
& Standard deviation
& $\bP\bracks{\abs*{\bm{\theta}}\geq\sqrt{3}/2}$
& $\bE[D(\bm{\theta})]$ \\
\hline
Prior
& 1.000000 & 0.386476 & 0.497255 \\
Limiting ABC
& 0.475324 & 0.000000 & 0.173484 \\
Gibbs ($\lambda=1$)
& 0.694645 & 0.219139 & 0.298952 \\
Gibbs ($\lambda=5$)
& 0.454945 & 0.046365 & 0.156001 \\
Gibbs ($\lambda=20$)
& 0.304087 & 0.000515 & 0.079122 \\
\hline
\end{tabular}
\caption{Comparison of the prior, limiting ABC posterior, and penalty-based Gibbs posteriors in the Gaussian model. Here $D(\bm{\theta})=\sqrt{\bm{\theta}^{\:\! 2}+1/4}-1/2$. Values are obtained by deterministic numerical integration.}
\label{tab:gaussian-gibbs}
\end{table}

\end{example}

\begin{example}[Finite-state models]
\label{ex:finite.state-models} 

Consider the case where $\cY = \Set{\! \by_{\:\! 1} , \dots , \by_{\:\! m} \!} \subset \R^{m_{\;\!\! \cY}}$\:\!\!, $m \in \N^{\:\! *}$\:\!\!, and let
\[
\bm{\mu}^{\;\!\! \star} \;\!\!
= \;\!\!
\sum_{i \:\! = \:\! 1}^{m}
q_{\;\!\! i}^{\star}
\bm{\delta}_{\by_{\:\! i}}
\]
and, for every $\bm{\theta} \in \Theta$,
\[
\bm{\nu}_{\bm{\theta}}
= \;\!\!
\sum_{i \:\! = \:\! 1}^{m}
q_{\;\!\! \bm{\theta} \:\!\!, i} \:\!
\bm{\delta}_{\by_{\:\! i}} ,
\]
where $q_{\;\!\! i}^{\star}$\;\!\! and $q_{\;\!\! \bm{\theta} \:\!\!, i}$, $i = 1 , \dots , m$, denote the corresponding probability weights.
Then
\[
W_{\:\!\! p}^{\:\! p \;\!\!}
\big{(}\:\! \bm{\mu}^{\;\!\! \star} \:\!\!,\;\!\! \bm{\nu}_{\bm{\theta}} \big{)}
\:\! = \:\!\!
\min_{\substack{
\gamma_{\:\! ij} \ge \;\! 0
\\
\sum_{\:\! j \:\! = \:\! 1}^{\:\! m} \;\!\! \gamma_{\:\! ij} \:\! = \;\! q_{i}^{\star}
\\
\sum_{\:\! i \:\! = \:\! 1}^{\:\! m} \;\!\! \gamma_{\:\! ij} \:\! = \;\! q_{\bm{\theta} \:\!\!, j}
}}
\sum_{i,j \:\! = \:\! 1}^{m} \;\!\!
\gamma_{\:\! i \;\!\! j} \:\!
\bd_{\:\!\! \cY}^{\;\! p}
\big{(}\:\! \by_{\:\! i} \:\!,\;\!\! \by_{\;\!\! j} \big{)} .
\]
Hence, the asymptotic compatibility profile is determined by finite-dimensional optimal transport problems, while the associated functional is obtained by integrating the corresponding transformed costs with respect to the parameter distribution.
This covers, among others, categorical and finite-state models.

\end{example}

The formulation illustrated by Examples~\ref{ex:gaussian-models} and~\ref{ex:finite.state-models} also applies to mixed discrete--continuous laws.
For instance, a zero-inflated model may take the form
\[
\bm{\nu}_{\bm{\theta}} \;\!\!
=
a_{\bm{\theta}} \:\!
\bm{\delta}_{\by_{\:\! 0}} \;\!\!
+
\big{(} 1 - a_{\bm{\theta}} \big{)} \:\!
\bm{\eta}_{\:\! \bm{\theta}} ,
\quad
a_{\bm{\theta}} \;\!\!
\in
\mo{[} 0 \:\!,\:\!\! 1 \mc{]} \:\! ,
\]
where $\by_{\:\! 0} \in \cY$\;\!\!, and $\bm{\eta}_{\:\! \bm{\theta}} \;\!\! \in \ccP_{\;\!\! p}(\cY)$ may be, for example, a Gaussian distribution or a finite Gaussian mixture.
The corresponding Wasserstein discrepancy is defined directly between the resulting hybrid probability measures, without requiring separate transport criteria for their discrete and continuous components.


Deterministic models are recovered as an immediate special case.
Indeed, if
\[
\bm{\mu}^{\;\!\! \star} \;\!\!
=
\bm{\delta}_{\by^{\star}}
\]
and, for every $\bm{\theta} \in \Theta$,
\[
\bm{\nu}_{\bm{\theta}} \;\!\!
=
\bm{\delta}_{\by_{\bm{\theta}}} ,
\]
where $\by^{\star} \:\!\! , \by_{\:\! \bm{\theta}} \;\!\! \in \cY$\;\!\!, then
\[
W_{\:\!\! p} \big{(}\:\! \bm{\mu}^{\;\!\! \star} \:\!\!,\;\!\! \bm{\nu}_{\bm{\theta}} \big{)} \;\!\!
=
\bd_{\;\!\! \cY}
\big{(}\:\! \by^{\star} \:\!\!,\;\!\! \by_{\:\! \bm{\theta}} \big{)} .
\]
Therefore, in the deterministic setting, the Wasserstein discrepancy reduces exactly to the underlying distance between the observational and model outputs.



\section{Conclusions}
\label{sec:conclusions} 

We have formulated compatibility-weight inference so that the normalization argument, the sampling assumptions, and the discrepancy geometry can be examined separately.
For vanishing tolerances, prior-integrated convergence of rescaled weights gives total variation convergence of the corresponding posterior distributions.
For fixed supercritical tolerances and increasing sample size, we recover the conditional-prior mechanism established in coarsened inference and Wasserstein ABC~\cite{miller-dunson_2019,bernton-jacob-gerber-robert_2019}, here with expected total variation convergence under general representation and stabilization assumptions.
Under partial stabilization, the persistence results provide lower bounds on observation-averaged acceptance, without identifying a limiting posterior.

The Gibbs identity gives convergence of minimum values and stability of approximate minimizers as the tolerance vanishes.
The fixed-tolerance limit has its own exact variational characterization: relative-entropy minimization subject to full mass on the compatibility region.
Additional finite penalties on the limiting discrepancy define different generalized Bayesian updates.
For Wasserstein discrepancies, their entropy-regularized problems reduce exactly to distributions over model laws, and power penalties admit the distance-to-ball representation under a geodesic assumption.
The Gaussian example illustrates fixed-tolerance filtering with quantified total variation errors and compares the resulting conditional prior with penalty-based Gibbs updates across different penalty scales.

Open questions include joint sample-size and tolerance regimes, local analysis at critical tolerances, and variational convergence of the finite-sample objectives.
Extensions to other regularizers require corresponding control of approximate minimizers.
Conditional or sequential models, uncertain limiting observational laws, and connections with distributionally robust optimization and multivariate risk assessment also merit investigation; see, e.g.,~\cite{esfahani-kuhn_2018,gao-kleywegt_2023,mastrogiacomo-rocca-tarsia_2026}.


\begin{thebibliography}{99}

\bibitem{agueh-carlier_2011}
Agueh, M., \& Carlier, G. (2011).
\newblock \textit{Barycenters in the Wasserstein space}.
\newblock SIAM Journal on Mathematical Analysis, 43(2), 904--924.

\bibitem{ambrosio-gigli-savare_2008}
Ambrosio, L., Gigli, N., \& Savar\'e, G. (2008).
\newblock \textit{Gradient Flows: In Metric Spaces and in the Space of Probability Measures}.
\newblock Birkh\"auser, Basel.


\bibitem{beaumont-zhang-balding_2002}
Beaumont, M.~A., Zhang, W., \& Balding, D.~J. (2002).
\newblock \textit{Approximate Bayesian computation in population genetics}.
\newblock Genetics, 162, 2025--2035.

\bibitem{bernton-jacob-gerber-robert_2019}
Bernton, E., Jacob, P.~E., Gerber, M., \& Robert, C.~P. (2019).
\newblock \textit{Approximate Bayesian computation with the Wasserstein distance}.
\newblock Journal of the Royal Statistical Society: Series B (Statistical Methodology), 81(2), 235--269.

\bibitem{bissiri-holmes-walker_2016}
Bissiri, P.~G., Holmes, C.~C., \& Walker, S.~G. (2016).
\newblock \textit{A General Framework for Updating Belief Distributions}.
\newblock Journal of the Royal Statistical Society: Series B (Statistical Methodology), 78(5), 1103--1130.


\bibitem{boissard-le.gouic_2014}
Boissard, E., \& Le~Gouic, T. (2014).
\newblock \textit{On the mean speed of convergence of empirical and occupation measures in Wasserstein distance}.
\newblock Annales de l'Institut Henri Poincar\'e, Probabilit\'es et Statistiques, 50(2), 539--563.

\bibitem{chen-mira-onnela_2020}
Chen, S., Mira, A., \& Onnela, J.-P. (2020).
\newblock \textit{Flexible model selection for mechanistic network models}.
\newblock Journal of Complex Networks, 8(2), cnz024.

\bibitem{dal.maso_1993}
Dal Maso, G. (1993).
\newblock \textit{An Introduction to $\Gamma$\;\!\!-\:\!Convergence}.
\newblock Progress in Nonlinear Differential Equations and Their Applications, Vol.~8.
\newblock Birkh\"auser, Boston.

\bibitem{dembo-zeitouni_1998}
Dembo, A., \& Zeitouni, O. (1998).
\newblock \textit{Large Deviations Techniques and Applications}.
\newblock 2nd~ed.
\newblock Applications of Mathematics, Vol.~38.
\newblock Springer, New York.

\bibitem{donsker-varadhan_1975}
Donsker, M.~D., \& Varadhan, S.~R.~S. (1975).
\newblock \textit{Asymptotic evaluation of certain Markov process expectations for large time, I}.
\newblock Communications on Pure and Applied Mathematics, 28(1), 1--47.

\bibitem{dupuis-ellis_1997}
Dupuis, P., \& Ellis, R.~S. (1997).
\newblock \textit{A Weak Convergence Approach to the Theory of Large Deviations}.
\newblock Wiley Series in Probability and Statistics.
\newblock John Wiley \& Sons, New York.

\bibitem{dutta-mira-onnela_2018}
Dutta, R., Mira, A., \& Onnela, J.-P. (2018).
\newblock \textit{Bayesian inference of spreading processes on networks}.
\newblock Proceedings of the Royal Society A: Mathematical, Physical and Engineering Sciences, 474(2215), 20180129.

\bibitem{ebert-dutta-mengersen-mira-ruggeri-wu_2018}
Ebert, A., Dutta, R., Mengersen, K., Mira, A., Ruggeri, F., \& Wu, P. (2018).
\newblock \textit{Likelihood-free parameter estimation for dynamic queueing networks: Case study of passenger flow in an international airport terminal}.
\newblock arXiv preprint, \href{https://arxiv.org/abs/1804.02526}{arXiv:1804.02526}.


\bibitem{fearnhead-prangle_2012}
Fearnhead, P., \& Prangle, D. (2012).
\newblock \textit{Constructing summary statistics for approximate Bayesian computation: semi-automatic approximate Bayesian computation}.
\newblock Journal of the Royal Statistical Society: Series B (Statistical Methodology), 74(3), 419--474.

\bibitem{fournier-guillin_2015}
Fournier, N., \& Guillin, A. (2015).
\newblock \textit{On the rate of convergence in Wasserstein distance of the empirical measure}.
\newblock Probability Theory and Related Fields, 162(3--4), 707--738.

\bibitem{frazier-martin-robert-rousseau_2018}
Frazier, D.~T., Martin, G.~M., Robert, C.~P., \& Rousseau, J. (2018).
\newblock \textit{Asymptotic properties of approximate Bayesian computation}.
\newblock Biometrika, 105(3), 593--607.

\bibitem{frazier-robert-rousseau_2020}
Frazier, D.~T., Robert, C.~P., \& Rousseau, J. (2020).
\newblock \textit{Model misspecification in approximate Bayesian computation: Consequences and diagnostics}.
\newblock Journal of the Royal Statistical Society: Series B (Statistical Methodology), 82(2), 421--444.

\bibitem{gao-kleywegt_2023}
Gao, R., \& Kleywegt, A.~J. (2023).
\newblock \textit{Distributionally robust stochastic optimization with Wasserstein distance}.
\newblock Mathematics of Operations Research, 48(2), 603--655.

\bibitem{gelbrich_1990}
Gelbrich, M. (1990).
\newblock \textit{On a formula for the $L^2$\;\!\! Wasserstein metric between measures on Euclidean and Hilbert spaces}.
\newblock Mathematische Nachrichten, 147(1), 185--203.

\bibitem{goyal-onnela_2023}
Goyal, R., \& Onnela, J.-P. (2023).
\newblock \textit{Framework for converting mechanistic network models to probabilistic models}.
\newblock Journal of Complex Networks, 11(5), cnad034.

\bibitem{jiang-wu-wong_2018}
Jiang, B., Wu, T.-Y., \& Wong, W.~H. (2018).
\newblock \textit{Approximate Bayesian computation with Kullback--Leibler divergence as data discrepancy}.
\newblock Proceedings of Machine Learning Research, 84, 1711--1721.

\bibitem{kantorovich_1942}
Kantorovich, L.~V. (1942).
\newblock \textit{On the translocation of masses}.
\newblock Doklady Akademii Nauk SSSR, 37(7--8), 227--229.

\bibitem{lin-ho-chen-cuturi-jordan_2020}
Lin, T., Ho, N., Chen, X., Cuturi, M., \& Jordan, M.~I. (2020).
\newblock \textit{Fixed-Support Wasserstein Barycenters: Computational Hardness and Fast Algorithm}.
\newblock Advances in Neural Information Processing Systems, 33, 5368--5380.

\bibitem{marin-pudlo-robert-ryder_2012}
Marin, J.-M., Pudlo, P., Robert, C.~P., \& Ryder, R.~J. (2012).
\newblock \textit{Approximate Bayesian computational methods}.
\newblock Statistics and Computing, 22, 1167--1180.

\bibitem{mastrogiacomo-rocca-tarsia_2026}
Mastrogiacomo, E., Rocca, M., \& Tarsia, M. (2026).
\newblock \textit{Multi-Objective Optimization and its Connection to Multivariate Risk Measures}.
\newblock Journal of Optimization Theory and Applications, 209(1), 30.

\bibitem{matsubara-knoblauch-briol-oates_2022}
Matsubara, T., Knoblauch, J., Briol, F.-X., \& Oates, C.~J. (2022).
\newblock \textit{Robust generalised Bayesian inference for intractable likelihoods}.
\newblock Journal of the Royal Statistical Society: Series B (Statistical Methodology), 84(3), 997--1022.

\bibitem{mccann_1997}
McCann, R.~J. (1997).
\newblock \textit{A convexity principle for interacting gases}.
\newblock Advances in Mathematics, 128(1), 153--179.

\bibitem{miller-dunson_2019}
Miller, J.~W., \& Dunson, D.~B. (2019).
\newblock \textit{Robust Bayesian inference via coarsening}.
\newblock Journal of the American Statistical Association, 114(527), 1113--1125.

\bibitem{esfahani-kuhn_2018}
Mohajerin~Esfahani, P., \& Kuhn, D. (2018).
\newblock \textit{Data-driven distributionally robust optimization using the Wasserstein metric: Performance guarantees and tractable reformulations}.
\newblock Mathematical Programming, 171(1--2), 115--166.

\bibitem{nguyen-arbel-lu-forbes_2019}
Nguyen, H.~D., Arbel, J., L\"u, H., \& Forbes, F. (2019).
\newblock \textit{Approximate Bayesian computation via the energy statistic}.
\newblock IEEE Transactions on Signal Processing, 67(19), 5005--5019.

\bibitem{pacchiardi-khoo-dutta_2024}
Pacchiardi, L., Khoo, S., \& Dutta, R. (2024).
\newblock \textit{Generalized Bayesian likelihood-free inference}.
\newblock Electronic Journal of Statistics, 18(2), 3628--3686.

\bibitem{panaretos-zemel_2019}
Panaretos, V.~M., \& Zemel, Y. (2019).
\newblock \textit{Statistical aspects of Wasserstein distances}.
\newblock Annual Review of Statistics and Its Application, 6, 405--431.

\bibitem{parthasarathy_1967}
Parthasarathy, K.~R. (1967).
\newblock \textit{Probability Measures on Metric Spaces}.
\newblock Academic Press, New York--London.

\bibitem{park-jitkrittum-sejdinovic_2016}
Park, M., Jitkrittum, W., \& Sejdinovic, D. (2016).
\newblock \textit{K2-ABC: Approximate Bayesian computation with kernel embeddings}.
\newblock Proceedings of Machine Learning Research, 51, 398--407.

\bibitem{peyre-cuturi_2019}
Peyr\'e, G., \& Cuturi, M. (2019).
\newblock \textit{Computational Optimal Transport}.
\newblock Foundations and Trends in Machine Learning, 11(5--6), 355--607.

\bibitem{santambrogio_2015}
Santambrogio, F. (2015).
\newblock \textit{Optimal Transport for Applied Mathematicians: Calculus of Variations, PDEs, and Modeling}.
\newblock Birkh\"auser, Cham.

\bibitem{sisson-fan-beaumont_2019}
Sisson, S.~A., Fan, Y., \& Beaumont, M.~A. (2019).
\newblock \textit{Handbook of Approximate Bayesian Computation}.
\newblock Chapman \& Hall/CRC, Boca Raton.

\bibitem{sriperumbudur-gretton-fukumizu-scholkopf-lanckriet_2010}
Sriperumbudur, B.~K., Gretton, A., Fukumizu, K., Sch\"olkopf, B., \& Lanckriet, G.~R.~G. (2010).
\newblock \textit{Hilbert space embeddings and metrics on probability measures}.
\newblock Journal of Machine Learning Research, 11, 1517--1561.


\bibitem{villani_2003}
Villani, C. (2003).
\newblock \textit{Topics in Optimal Transportation}.
\newblock American Mathematical Society, Providence, RI.

\bibitem{villani_2009}
Villani, C. (2009).
\newblock \textit{Optimal Transport: Old and New}.
\newblock Springer, Berlin.


\bibitem{weed-bach_2019}
Weed, J., \& Bach, F. (2019).
\newblock \textit{Sharp asymptotic and finite-sample rates of convergence of empirical measures in Wasserstein distance}.
\newblock Bernoulli, 25(4A), 2620--2648.

\bibitem{wilkinson_2013}
Wilkinson, R.~D. (2013).
\newblock \textit{Approximate Bayesian computation (ABC) gives exact results under the assumption of model error}.
\newblock Statistical Applications in Genetics and Molecular Biology, 12(2), 129--141.

\end{thebibliography}
\end{document}